\documentclass{emsprocart}

\usepackage{srcltx}

\usepackage[all]{xy}
\usepackage{lscape}
\usepackage{tikz}
\usetikzlibrary{arrows}
\usetikzlibrary{snakes}
\usetikzlibrary{shapes}

\contact[amiot@math.unistra.fr]{Claire Amiot, Institut de Recherche Math\'ematique Avanc\'ee, 7 rue Descartes 67000 Strasbourg, France}

\newcommand{\Hom}{{\sf Hom }}\newcommand{\RHom}{\mathbf{R}{\sf Hom }}
\newcommand{\End}{{\sf End }}
\newcommand{\Ext}{{\sf Ext }}

\renewcommand{\mod}{{\sf mod \hspace{.02in}  }}

\newcommand{\Mod}{{\sf Mod \hspace{.02in} }}

\newcommand{\ind}{{\sf ind \hspace{.02in} }}
\newcommand{\per}{{\sf per \hspace{.02in}  }}
\newcommand{\proj}{{\sf proj \hspace{.02in} }}

\newcommand{\Sub}{{\sf Sub \hspace{.02in} }}

\newcommand{\add}{{\sf add \hspace{.02in} }}
\newcommand{\gr}{{\sf gr \hspace{.02in} }}

\newcommand{\Gr}{{\sf Gr \hspace{.02in} }}

\newcommand{\ten}{\otimes}
\newcommand{\lten}{\overset{\mathbf{L}}{\ten}}

\newcommand{\Talg}{{\sf T}}
\newcommand{\SL}{{\sf SL}}

\newcommand{\Jac}{{\sf Jac}}

\newcommand{\Cc}{\mathcal{C}}
\newcommand{\Dd}{\mathcal{D}}
\newcommand{\Ee}{\mathcal{E}}
\newcommand{\Ff}{\mathcal{F}}
\newcommand{\Hh}{\mathcal{H}}

\newcommand{\Tt}{\mathcal{T}}

\newcommand{\Rr}{\mathcal{R}}

\newcommand{\Vv}{\mathcal{V}}
\newcommand{\Uu}{\mathcal{U}}
\newcommand{\Xx}{\mathcal{X}}

\newcommand{\SSS}{\mathbb{S}}

\newcommand{\ZZ}{\mathbb{Z}}
\newcommand{\CM}{{\sf CM \hspace{.02in}}}
\newcommand{\Db}{\mathcal{D}^{\rm b}}
\newcommand{\BPi}{{\mathbf\Pi}}

\newcommand{\bsm}{\begin{smallmatrix}}
\newcommand{\esm}{\end{smallmatrix}}

\numberwithin{figure}{subsection}
\numberwithin{equation}{subsection}

\newtheorem{theorem}{Theorem}[section]

\newtheorem{corollary}[theorem]{Corollary}
\newtheorem{proposition}[theorem]{Proposition}

\newtheorem{question}[theorem]{Question}
\newtheorem{conjecture}[theorem]{Conjecture}

\theoremstyle{definition}
\newtheorem{definition}[theorem]{Definition}
\newtheorem{example}[theorem]{Example}
\newtheorem{remark}[theorem]{Remark}

\author[C. Amiot]{Claire Amiot}

\title{On generalized cluster categories}

\begin{document}

\maketitle

\begin{abstract}Cluster categories have been introduced by Buan, Marsh, Reineke, Reiten and Todorov in order to categorify Fomin-Zelevinsky cluster algebras. This survey motivates and outlines the construction of  a generalization of cluster categories, and explains different applications of these new categories in representation theory.
\end{abstract}

\medskip

\section*{Introduction}
In 2003, Marsh, Reineke and Zelevinsky attempted in \cite{MRZ} to understand Fomin-Zelevinsky cluster algebras (defined in \cite{FZ1}) in terms of representations of quivers. This article was immediately followed by the fundamental paper \cite{BMR+06} of Buan, Marsh, Reineke, Reiten and Todorov which stated the definition of cluster categories. In this paper, the authors associate to each finite dimensional hereditary algebra a triangulated category endowed with a special set of objects called cluster-tilting. The combinatorics of these objects is closely related to the combinatorics of acyclic cluster algebras, and especially with the mutation of quivers.

The same kind of phenomena appear naturally in the stable categories of modules over preprojective algebras of Dynkin type, and have been studied by Geiss, Leclerc and Schr\"oer in \cite{GLS06}, \cite{GLS07}. Actually most of the results on cluster-tilting theory of \cite{BMR+06}, \cite{BMR1}, \cite{BMR2}, \cite{GLS06}, \cite{GLS07} have been proved in the more general setting of $\Hom$-finite, triangulated, $2$-Calabi-Yau categories with cluster-tilting objects (see \cite{KR07}, \cite{IY08}, \cite{BIRS09}).

Since then, other categories with the above properties have been constructed and investigated. One finds stable subcategories of modules over a preprojective algebra of type $Q$ associated with any element in the Coxeter  group of $Q$ (see \cite{BIRS09}, \cite{BIRSm}, \cite{GLS}, \cite{GLS10}); stable Cohen-Macaulay modules over isolated singularities (see \cite{Aus86}, \cite{Iya07a}, \cite{BIKR08}); and generalized cluster categories associated with finite dimensional algebras of global dimension at most two, or with Jacobi-finite quivers with potential (introduced in \cite{Ami08}, \cite{Ami09}). 

The context of this survey is this last family of triangulated $2$-Calabi-Yau categories: the generalized cluster categories. The aim, here, is first to give some motivation for enlarging the family of cluster categories (Sections 1 and 2). Then we will explain the general construction of these new cluster categories (Section 3), and link these new categories with the categories cited above (Section 4). Finally, we give some applications in representation theory of the cluster categories associated with algebras of global dimension at most 2 (Section 5). 
\section*{Contents}
\begin{enumerate}
    \setlength{\topsep}{0pt}%
      \setlength{\parskip}{0pt}
    \setlength{\itemsep}{0pt}%
    \settowidth{\labelwidth}{9}%
    \setlength{\leftmargin}{0pt}%
    \setlength{\itemindent}{0pt}%
  \item[1\ \ ] 
    Motivation
    \begin{enumerate}
      \setlength{\topsep}{0pt}
      \setlength{\parskip}{0pt}
      \setlength{\itemsep}{0pt}
      \setlength{\leftmargin}{0pt}
      \setlength{\partopsep}{0pt}
      \item[1.1 ] 
        Mutation of quivers
      \item[1.2 ] 
        Cluster category of type $A_3$
      \item[1.3 ] 
        Stable module category of a preprojective algebra of type $A_3$
    \end{enumerate}
  \item[2\ \ ] 
    2-Calabi-Yau categories and cluster-tilting theory
    \begin{enumerate}
      \setlength{\topsep}{0pt}
      \setlength{\parskip}{0pt}
      \setlength{\itemsep}{0pt}
      \setlength{\leftmargin}{0pt}
      \item[2.1 ] 
        Iyama-Yoshino mutation
      \item[2.2 ] 
        2-Calabi-Yau-tilted algebras
      \item[2.3 ] 
        Two fundamental families of examples
      \item[2.4 ] 
        Quivers with potential
    \end{enumerate}
  \item[3\ \ ] 
    From $3$-Calabi-Yau DG-algebras to 2-Calabi-Yau categories
    \begin{enumerate}
      \setlength{\topsep}{0pt}
      \setlength{\parskip}{0pt}
      \setlength{\itemsep}{0pt}
      \setlength{\leftmargin}{0pt}
      \item[3.1 ] 
        Graded algebras and DG algebras
      \item[3.2 ] 
        General construction
      \item[3.3 ] 
        Ginzburg DG algebras
      \item[3.4 ] 
        Application to surfaces with marked points
      \item[3.5 ] 
        Derived preprojective algebras
      \item[3.6 ] 
        Link between the two constructions
      \item[3.7 ] 
        Higher generalized cluster categories
    \end{enumerate}
  \item[4\ \ ] 
    Stable categories as generalized cluster categories
    \begin{enumerate}
      \setlength{\topsep}{0pt}
      \setlength{\parskip}{0pt}
      \setlength{\itemsep}{0pt}
      \setlength{\leftmargin}{0pt}
      \item[4.1 ] 
        Preprojective algebras of Dynkin type
      \item[4.2 ] 
        Stable categories associated to words
      \item[4.3 ] 
        Cohen Macaulay modules over isolated singularities
    \end{enumerate}
  \item[5\ \ ] 
    On the $\mathbb {Z}$-grading on the $3$-preprojective algebra
    \begin{enumerate}
      \setlength{\topsep}{0pt}
      \setlength{\parskip}{0pt}
      \setlength{\itemsep}{0pt}
      \setlength{\leftmargin}{0pt}
      \item[5.1 ] 
        Grading of $\Pi _3(\Lambda )$
      \item[5.2 ] 
        Mutation of graded QPs
      \item[5.3 ] 
        Application to cluster equivalence
      \item[5.4 ] 
        Algebras of acyclic cluster type
      \item[5.5 ] 
        Gradable $\Pi _3{\Lambda }$-modules and piecewise hereditary algebras
    \end{enumerate}
\end{enumerate}

\subsection*{Acknowledgments}
I deeply thank Andrzej Skowro\'nski for encouraging me to write this survey. Some of the results presented here are part of my Ph.D. thesis under the supervision of Bernhard Keller. I thank him for his constant support and helpful comments and suggestions on a previous version of this survey. I am also grateful to Steffen Oppermann for helpful comments. I thank him and my other coauthors Osamu Iyama, Gordana Todorov and Idun Reiten for pleasant collaborations.

\subsection*{Notation}
$k$ is an algebraically closed field. All categories considered in this paper are $k$-linear additive Krull-Schmidt categories. By $\Hom$-finite categories, we mean, categories such that $\Hom(X,Y)$ is finite dimensional for any objects $X$ and $Y$.  We denote by $D=\Hom_k(-,k)$ the $k$-dual. The tensor products are over the field $k$ when not specified. For an object $T$ in a category $\Cc$, we denote by $\add(T)$ the additive closure of $T$, that is the smallest full subcategory of $\Cc$ containing $T$ and stable under taking direct summands. 

All modules considered here are right modules. For a $k$-algebra $A$, we denote by $\Mod A$ the category of right modules and by $\mod A$ the category of finitely presented right $A$-modules.

A quiver $Q=(Q_0,Q_1,s,t)$ is given by a set of vertices $Q_0$, a set of arrows $Q_1$, a source map $s:Q_1\to Q_0$ and a target map $t:Q_1\to Q_0$. We denote by $e_i, i\in Q_0$ the set of primitive idempotents of the path algebra $kQ$. By a Dynkin quiver, we mean a quiver whose underlying graph is of simply laced Dynkin type.

\section{Motivation}\label{section motivation}
\subsection{Mutation of quivers}
In the fundamental paper \cite{FZ1}, Fomin and Zelevinsky introduced the notion of mutation of quivers as follows.
\begin{definition}[Fomin-Zelevinsky \cite{FZ1}]
Let $Q$ be a finite quiver without loops and oriented cycles of length $2$ ($2$-cycles for short). Let $i$ be a vertex of $Q$. The \emph{mutation of the quiver $Q$ at the vertex $i$} is a quiver denoted by $\mu_i(Q)$ and constructed from $Q$ using the following rule:
\begin{itemize}
\item[(M1)] for any couple of arrows $j\to i\to k$, add an arrow $j\to k$;
\item[(M2)] reverse the arrows incident with $i$;
\item[(M3)] remove a maximal collection of $2$-cycles.
\end{itemize}
\end{definition}

This definition is one of the key steps in the definition of \emph{Cluster Algebras}. Even if the initial motivation of Fomin and Zelevinsky to define cluster algebras was to give a combinatorial and algebraic setup to understand canonical basis and total positivity in algebraic groups (see \cite{Lus}, \cite{Lus2}, \cite{GLS}), these algebras have led to many applications in very different domains of mathematics:
\begin{itemize}
\item representations of  groups of surfaces  (higher Teichm\"uller theory \cite{FG}, \cite{FG2}, \cite{FST08});
\item  discrete dynamical systems, (Y-systems, integrable systems \cite{FZ5}, \cite{Kel10});
\item non commutative algebraic geometry (Donaldson-Thomas invariants \cite{KS},  Calabi-Yau algebras \cite{G06});
\item Poisson geometry \cite{GSV};
\item quiver and finite dimensional algebras representations. 
\end{itemize}
This paper is focused on this last connection and on the new insights that cluster algebras have brought in representation theory. However, in order to make this survey not too long, we will not discuss cluster algebras and their precise links with representation theory, but we will concentrate on the link between certain categories and quiver mutation. We refer to \cite{Kellersurvey}, \cite{BM06}, \cite{Reitensurvey} for nice overviews of categorifications of cluster algebras.

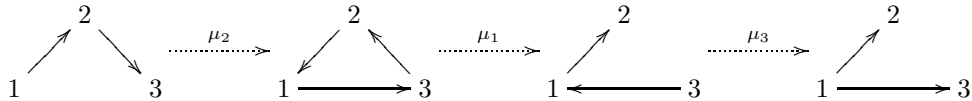
\begin{example}
$$\xymatrix@-.3cm{ & 2\ar[dr] &&&& 2\ar[dl] &&&& 2 &&&& 2 &\\
1\ar[ur] && 3\ar@<.5cm>@{..>}[rr]^{\mu_2} && 1\ar[rr] && 3\ar@<.5cm>@{..>}[rr]^{\mu_1}\ar[ul] && 1\ar[ur] && 3 \ar[ll]\ar@<.5cm>@{..>}[rr]^{\mu_3}&&1\ar[rr]\ar[ur]&&3}.$$
\end{example}
 One easily checks that the mutation at a given vertex is an involution. If $i$ is a source (i.e. there are no arrows with target $i$) or a sink (i.e. there are no arrows with source $i$), the mutation $\mu_i$ consists only in reversing arrows incident with $i$ (step (M2)). Hence it coincides with the reflection introduced by Bernstein, Gelfand, and Ponomarov \cite{BGP}. Therefore mutation of quivers can be seen as a generalization of reflections (which are just defined if $i$ is a source or a sink). 
 
The BGP reflections  have really nice applications in representation theory, and are closely related to tilting theory of hereditary algebras. For instance, combining the functorial interpretation of reflections functors by Brenner and Butler \cite{BB}, the interpretation of tilting modules for derived categories by Happel \cite[Thm~1.6]{Hap87} and its description of the derived categories of hereditary algebras \cite[Cor.~4.8]{Hap87}, we obtain that the reflections characterize combinatorially derived equivalence between hereditary algebras.  
 \begin{theorem}[Happel]\label{Happel}
 Let $Q$ and $Q'$ be two acyclic quivers. Then the algebras $kQ$ and $kQ'$ are derived equivalent if and only if one can pass from $Q$ to $Q'$ using a finite sequence of reflections.
\end{theorem} 

One hopes therefore that the notion of quiver mutation has also nice applications and meaning in representation theory. 

In the rest of the section, we give two fundamental examples of categories where the combinatorics of the quiver mutation of $1\to 2\to 3$ appear naturally. These categories come from representation theory, they are very close to categories of finite dimensional modules over finite dimensional algebras.

\subsection{Cluster category of type $A_3$}\label{subsection CA3}

Let $Q$ be the following quiver $1\to 2\to 3$. We consider the category $\mod kQ$ of finite dimensional (right) modules over the path algebra $kQ$. We refer to the books \cite{ASS}, \cite{ARS}, \cite{GR}, \cite{Rin} for a wealth of information on the representation theory of quivers and finite dimensional algebras.

Since $A_3$ is a Dynkin quiver, this category has finitely many indecomposable modules. The Auslander-Reiten quiver is the following:
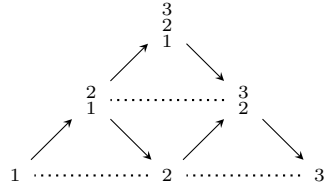
\[\scalebox{1}{
\begin{tikzpicture}[>=stealth,scale=1]
\node (P2) at (2,0) {${\bsm 1\esm}$};
\node (P3) at (4,0) {${\bsm 2\esm}$};
\node (P4) at (6,0){${\bsm 3\esm}$};
\node (B2) at (3,1) {${\bsm 2\\1\esm}$};
\node (B3) at (5,1) {${\bsm 3\\2\esm}$};
\node (A3) at (4,2) {${\bsm 3\\2\\1\esm}$};

\draw [->] (P2)  --  (B2);
\draw[->](P3)--(B3);
\draw[->](B2)--(P3);
\draw [->] (B3)  --  (P4);
\draw [->] (B2)  --  (A3);
\draw [dotted, thick] (B3)  --  (B2);
\draw [dotted, thick] (P4)  --  (P3);
\draw[dotted, thick] (P3)--(P2);
\draw[->] (A3)--(B3);

\end{tikzpicture}}
\]

Here the simples associated with the vertices are symbolized by ${\bsm1\esm}$, ${\bsm2\esm}$, ${\bsm3\esm}$. The module ${\bsm3\\2\\1\esm}$ is an indecomposable module $M$ such that there exists a filtration $M_1\subset M_2\subset M_3=M$ with $M_1\simeq{\bsm 1\esm}$, $M_2/M_1\simeq {\bsm 2 \esm}$ and $M_3/M_2\simeq {\bsm 3 \esm}$. In this situation, it determines the module $M$ up to isomorphism. There are three projective indecomposable objects corresponding to the vertices of $Q$ and which are $P_1=e_1kQ={\bsm 1 \esm}$, $P_2=e_2kQ={\bsm 2\\ 1\esm}$ and $P_3=e_3kQ={\bsm3\\2\\1\esm}$. The dotted lines describe the AR-translation that we denote by $\tau$; it induces a bijection between indecomposable non projective $kA_3$-modules and indecomposable non injective $kA_3$-modules. Dotted lines also correspond to minimal relations in the AR-quiver of $\mod kQ$.

\bigskip
The bounded derived category $\Db(A)$, where $A$ is a finite dimensional algebra of finite global dimension, is a $k$-category whose objects are bounded complexes of finite dimensional (right) $A$-modules. Its morphisms are obtained from morphisms of complexes by inverting formally quasi-isomorphisms. We refer to \cite{Hap} (see also \cite{Kel07}) for more precise description. For an object $M\in \Db(A)$ we denote by $M[1]$ the shift complex defined by $M[1]^n:=M^{n+1}$ and $d_{M[1]}=-d_M$. This category is a triangulated category, with suspension functor $M\to M[1]$ corresponding to the shift. The functor $\mod A\to \Db(A)$ which sends a $A$-module $M$ on the complex $\ldots \to 0\to M\to 0\ldots$ concentrated in degree $0$ is fully faithful. Moreover for any $M$ and $N$ in $\mod(A)$ and $i\in \ZZ$ we have:
$$\Hom_{\Db(A)}(M,N[i])\simeq \Ext^{i}_{A}(M,N).$$ The category $\Db(A)$  has a Serre functor $\SSS=-\lten_A DA$ (which sends the projective $A$-module $e_i A$ on the injective $A$-module $I_i=e_i DA$) and an AR-translation $\tau=-\lten_A DA[-1]$ which extends the AR-translation of $\mod A$ and that we still denote by $\tau$.

\bigskip
In the case $A=kQ$, the indecomposable objects are isomorphic to stalk complexes, that is, are of the form $X[i]$, where $i\in\mathbb{Z}$ and $X$ is an indecomposable $kQ$-module. The AR-quiver of $\Db(kQ)$ is the following infinite quiver (cf. \cite{Hap87}):

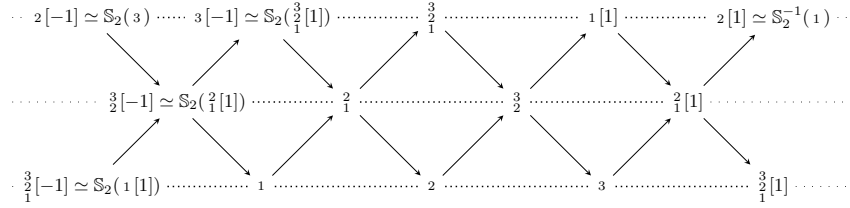
\[\scalebox{.7}{
\begin{tikzpicture}[>=stealth,scale=1.6]
\node (P0) at (-1,0) {};
\node (P1) at (0,0) {${\bsm 3\\2\\1\esm}[-1]\simeq\SSS_2({\bsm 1\esm}[1])$};
\node (P2) at (2,0) {${\bsm 1\esm}$};
\node (P3) at (4,0) {${\bsm 2\esm}$};
\node (P4) at (6,0){${\bsm 3\esm}$};
\node (P5) at (8,0){${\bsm 3\\2\\1\esm}[1]$};
\node (P6) at (9,0) {};
\node (B0) at (-1,1) {};
\node (B1) at (1,1) {${\bsm 3\\2\esm}[-1]\simeq\SSS_2({\bsm 2\\1\esm}[1])$};
\node (B2) at (3,1) {${\bsm 2\\1\esm}$};
\node (B3) at (5,1) {${\bsm 3\\2\esm}$};
\node (B4) at (7,1) {${\bsm 2\\1\esm}[1]$};
\node (B5) at (9,1) {};
\node (A0) at (-1,2) {};
\node (A1) at (0,2) {${\bsm 2\esm}[-1]\simeq\SSS_2({\bsm 3\esm})$};
\node (A2) at (2,2) {${\bsm 3\esm}[-1]\simeq\SSS_2({\bsm 3\\2\\1\esm}[1])$};
\node (A3) at (4,2) {${\bsm 3\\2\\1\esm}$};
\node (A4) at (6,2){${\bsm 1\esm}[1]$};
\node (A5) at (8,2){${\bsm 2\esm}[1]\simeq\SSS_2^{-1}({\bsm 1\esm})$};
\node (A6) at (9,2) {};
\draw [->] (P1)  --  (B1);
\draw [->] (P2)  --  (B2);
\draw [->] (P3)  --  (B3);
\draw [->] (P4)  --  (B4);
\draw [->] (B1)  --  (P2);
\draw [->] (B2)  --  (P3);
\draw [->] (B3)  --  (P4);
\draw [->] (B4)  --  (P5);
\draw [->] (A1)  --  (B1);
\draw [->] (A2)  --  (B2);
\draw [->] (A3)  --  (B3);
\draw [->] (A4)  --  (B4);
\draw [->] (B1)  --  (A2);
\draw [->] (B2)  --  (A3);
\draw [->] (B3)  --  (A4);
\draw [->] (B4)  --  (A5);
\draw [loosely dotted] (P0) -- (P1);\draw [loosely dotted] (P5) -- (P6);\draw [loosely dotted] (B0) -- (B1);\draw [loosely dotted] (B4) -- (B5);\draw [loosely dotted] (A0) -- (A1);\draw [loosely dotted] (A5) -- (A6);

\draw [dotted, thick] (P1) -- (P2);
\draw [dotted, thick] (P2) --(P3);
\draw [dotted, thick] (P3) -- (P4);
\draw [dotted, thick] (P4) --(P5);\draw [dotted, thick] (B1) -- (B2);\draw [dotted, thick] (B2) -- (B3);
\draw [dotted, thick] (B3) --(B4);\draw [dotted, thick] (A1) -- (A2);\draw [dotted, thick] (A2)--(A3);
\draw [dotted, thick] (A3) --(A4);\draw [dotted, thick] (A4) -- (A5);

\end{tikzpicture}}
\]
 
The functor $\SSS_2:=\SSS[-2]=\tau[-1]$ acts bijectively on the indecomposable objects of $\Db(kQ)$. It is an auto-equivalence of $\Db(kQ)$.

\bigskip
The cluster category $\Cc_{Q}$ is defined to be the orbit category of $\Db(kQ)$ by the functor $\SSS_2$: the indecomposable objects are the indecomposable objects of $\Db(kQ)$. The space of morphisms  between two objects in $\Db(kQ)$ is given by $$\Hom_{\Cc}(X,Y)=\bigoplus_{i\in \ZZ}\Hom_{\Dd}(X,\SSS_2^i Y).$$
Hence the objects $X$ and $\SSS_2^{i}X$ become isomorphic in $\Cc_Q$ for any $i\in \ZZ$. We denote by $\pi(X)$ the $\SSS_2$-orbit of an indecomposable $X$. One can see that there are finitely many indecomposable objects in $\Cc_{Q}$. They are of the form $\pi(X)$, where $X$ is an indecomposable $kQ$-module, or $X\simeq P[1]$ where $P$ is an indecomposable projective $kQ$-module. Hence the AR-quiver of the category $\Cc_{Q}$ is the following:

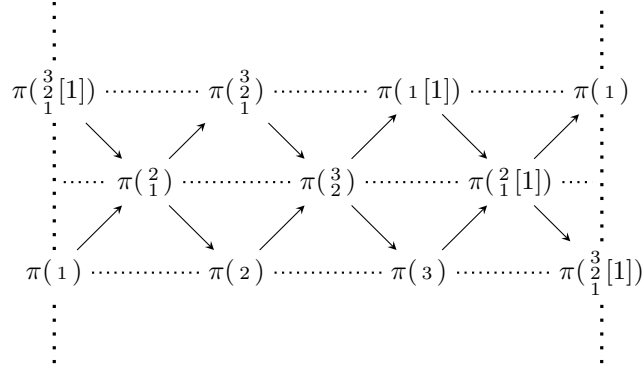
\[\scalebox{1}{
\begin{tikzpicture}[>=stealth,scale=1.2]
\node (P2) at (2,0) {$\pi({\bsm 1\esm})$};
\node (P3) at (4,0) {$\pi({\bsm 2\esm})$};
\node (P4) at (6,0){$\pi({\bsm 3\esm})$};
\node (P5) at (8,0){$\pi({\bsm 3\\2\\1\esm}[1])$};

\node (B1) at (2,1) {};
\node (B2) at (3,1) {$\pi({\bsm 2\\1\esm})$};
\node (B3) at (5,1) {$\pi({\bsm 3\\2\esm})$};
\node (B4) at (7,1) {$\pi({\bsm 2\\1\esm}[1])$};
\node (B5) at (8,1) {};
\node (A2) at (2,2) {$\pi({\bsm 3\\2\\1\esm}[1])$};
\node (A3) at (4,2) {$\pi({\bsm 3\\2\\1\esm})$};
\node (A4) at (6,2){$\pi({\bsm 1\esm}[1])$};
\node (A5) at (8,2){$\pi({\bsm 1\esm})$};

\draw [->] (P2)  --  (B2);
\draw [->] (P3)  --  (B3);
\draw [->] (P4)  --  (B4);
\draw [->] (B2)  --  (P3);
\draw [->] (B3)  --  (P4);
\draw [->] (B4)  --  (P5);
\draw [->] (A2)  --  (B2);
\draw [->] (A3)  --  (B3);
\draw [->] (A4)  --  (B4);
\draw [->] (B2)  --  (A3);
\draw [->] (B3)  --  (A4);
\draw [->] (B4)  --  (A5);

\draw [loosely dotted, very thick] (2,-1) -- (P2); \draw [loosely dotted, very thick] (P2) -- (A2);\draw [loosely dotted, very thick] (A2) -- (2,3);
 \draw [loosely dotted, very thick] (8,-1)--(P5);\draw [loosely dotted, very thick] (P5) -- (A5);\draw [loosely dotted, very thick] (A5) -- (8,3);
\draw [dotted, thick] (P2) --(P3);
\draw [dotted, thick] (P3) -- (P4);
\draw [dotted, thick] (P4) --(P5); \draw [dotted, thick] (B1) -- (B2);\draw [dotted, thick] (B2) -- (B3);
\draw [dotted, thick] (B4) -- (B5);
\draw [dotted, thick] (B3) --(B4);\draw [dotted, thick] (A2)--(A3);
\draw [dotted, thick] (A3) --(A4);\draw [dotted, thick] (A4) -- (A5);

\end{tikzpicture}}
\]
where the two vertical lines are identified. Therefore one can view the category $\Cc_Q$ as the module category $\mod kQ$ with extra objects ($\pi({\bsm1\esm}[1])$, $\pi({\bsm2\\1\esm}[1])$, $\pi({\bsm3\\2\\1\esm}[1])$) and extra morphisms: if $M$ and $N$ are $kQ$-modules, then
$$\Hom_\Cc(M,N)\simeq \Hom_{kQ}(M,N)\oplus \Ext^1_{kQ}(M,\tau^- N).$$
Note that in the category $\Cc_Q$, since the functor $\SSS_2=\SSS[-2]=\tau[-1]$ is isomorphic to the identity, then the functors $\tau$ and $[1]$ are isomorphic.

We are interested in the objects of $\Cc_Q$ which are rigid (i.e. satisfying $\Hom_\Cc(X,X[1])=0$), basic (i.e. with pairwise non-isomorphic direct summands), and maximal for these properties. An easy computation yields 14 maximal basic rigid objects which are :

\[ \tau^i(\pi({\bsm1\esm})\oplus\pi({\bsm2\\1\esm})\oplus\pi({\bsm3\\2\\1\esm})), i=0,\ldots, 5;  \quad \tau^i(\pi({\bsm2\esm})\oplus\pi({\bsm2\\1\esm})\oplus\pi({\bsm3\\2\\1\esm})), i=0,\ldots, 2;\]
\[ \tau^i(\pi({\bsm2\esm})\oplus\pi({\bsm3\\2\esm})\oplus\pi({\bsm3\\2\\1\esm})), i=0,\ldots, 2;\quad \tau^i(\pi({\bsm1\esm})\oplus\pi({\bsm 3\esm})\oplus\pi({\bsm3\\2\\1\esm})), i=0,1.\]

All these maximal rigid objects have the same number of indecomposable summands. The Gabriel quivers of their endomorphism algebras are the following:
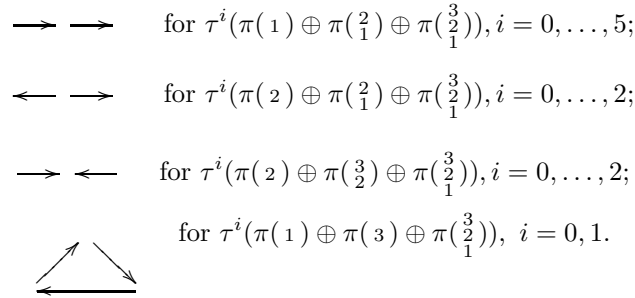
\[ \xymatrix@-.3cm{\ar[r] &  \ar[r] & } \quad \textrm{ for }\tau^i(\pi({\bsm1\esm})\oplus\pi({\bsm2\\1\esm})\oplus\pi({\bsm3\\2\\1\esm})), i=0,\ldots, 5; \]
\[ \xymatrix@-.3cm{ & \ar[l] \ar[r] & } \quad \textrm{ for }\tau^i(\pi({\bsm2\esm})\oplus\pi({\bsm2\\1\esm})\oplus\pi({\bsm3\\2\\1\esm})), i=0,\ldots, 2;\]

\[ \xymatrix@-.3cm{\ar[r] &   & \ar[l]} \quad \textrm{for }\tau^i(\pi({\bsm2\esm})\oplus\pi({\bsm3\\2\esm})\oplus\pi({\bsm3\\2\\1\esm})), i=0,\ldots, 2;\]
\[ \xymatrix@-.3cm{& \ar[dr]  & \\ \ar[ur]&&\ar[ll]} \quad \textrm{for }\tau^i(\pi({\bsm1\esm})\oplus\pi({\bsm 3\esm})\oplus\pi({\bsm3\\2\\1\esm})),\ i=0,1.\]

 These objects satisfy a remarkable property: given a maximal basic rigid object, one can replace an indecomposable summand by another one in a unique way to find another maximal basic rigid object. 
 This process can be understood as a mutation, and under this process the quiver of the endomorphism algebra changes according to the mutation rule. For instance, if we write the maximal rigid object together with its quiver we obtain:

$$
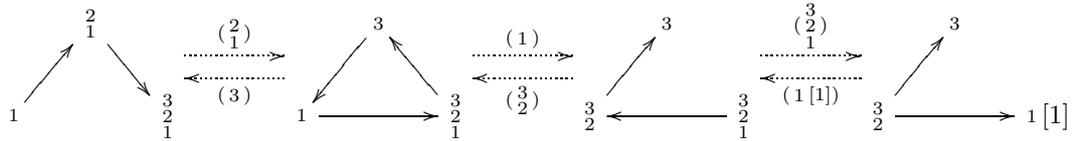
\xymatrix@-.3cm{ & {\bsm2\\1\esm} \ar[dr] &&&& {\bsm 3\esm}\ar[dl] &&&& {\bsm3\esm} &&&&  {\bsm3\esm} &\\
{\bsm1\esm} \ar[ur] && {\bsm3\\2\\1\esm} \ar@<.8cm>@{..>}[rr]^{({\bsm2\\1\esm})} && {\bsm1\esm}\ar[rr]\ar@<-.5cm>@{..>}[ll]^{({\bsm 3\esm})} && {\bsm3\\2\\1\esm}\ar@<.8cm>@{..>}[rr]^{({\bsm1\esm})}\ar[ul] && \ar@<-.5cm>@{..>}[ll]^{({\bsm3\\2\esm})}{\bsm3\\2\esm}\ar[ur] && {\bsm3\\2\\1\esm} \ar[ll]\ar@<.8cm>@{..>}[rr]^{ ({\bsm3\\2\\1\esm})}&&\ar@<-.5cm>@{..>}[ll]^{({\bsm1\esm}[1])}{\bsm3\\2\esm}\ar[rr]\ar[ur]&&{\bsm1\esm}[1]}.$$
 
\subsection{Stable module category of a preprojective algebra of type $A_3$}\label{subsection PiA3}

Let $Q$ be the quiver $1\to2\to 3$. Then the preprojective algebra $\Pi_2(kQ)$ (see \cite{Rin98}) is presented by the quiver
$$\xymatrix{1\ar@/^/[r]^a & 2\ar@/^/[r]^b\ar@/^/[l]^{a^*} & 3\ar@/^/[l]^{b^*}}\quad \textrm{with relations }a^*a=0,\ aa^*-b^*b=0, \  bb^*=0.$$
It is a finite dimensional algebra. 
The projective indecomposable $\Pi_2(kQ)$-modules are (up to isomorphism) $P_1=I_3={\bsm  1\\2\\3\esm}$, $P_2=I_2={\bsm  &2&\\3&&1\\ &2&\esm}$ and $P_3=I_1={\bsm  3\\2\\1\esm}$. They are also injective, since the algebra $\Pi_2(kQ)$ is self-injective.

The category $\mod \Pi_2(kQ)$ has finitely many indecomposable objects (up to isomorphism) and its AR-quiver is the following:

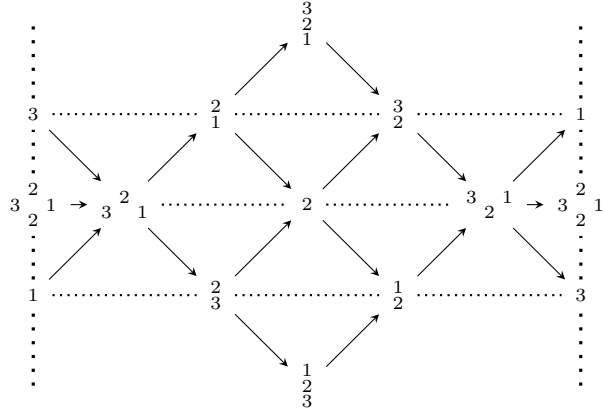
\[\scalebox{1}{
\begin{tikzpicture}[>=stealth,scale=1.2]
\node (P2) at (2,0) {${\bsm 1\esm}$};
\node (P3) at (4,0) {${\bsm 2\\3\esm}$};
\node (P4) at (6,0){${\bsm 1\\2\esm}$};
\node (P5) at (8,0){${\bsm 3\esm}$};
\node (P6) at (5,-1){${\bsm 1\\2\\3\esm}$};
\node (B1) at (2,1) {${\bsm &2&\\3&&1\\&2&\esm}$};
\node (B2) at (3,1) {${\bsm &2&\\3&&1\esm}$};
\node (B3) at (5,1) {${\bsm 2\esm}$};
\node (B4) at (7,1) {${\bsm 3&&1\\&2&\esm}$};
\node (B5) at (8,1) {${\bsm &2&\\3&&1\\&2&\esm}$};
\node (A2) at (2,2) {${\bsm 3\esm}$};
\node (A3) at (4,2) {${\bsm 2\\1\esm}$};
\node (A4) at (6,2){${\bsm 3\\2\esm}$};
\node (A5) at (8,2){${\bsm 1\esm}$};
\node (A6) at (5,3){${\bsm 3\\2\\1\esm}$};

\draw [->] (P2)  --  (B2);
\draw [->] (P3)  --  (B3);
\draw [->] (P4)  --  (B4);
\draw [->] (B2)  --  (P3);
\draw [->] (B3)  --  (P4);
\draw [->] (B4)  --  (P5);
\draw [->] (A2)  --  (B2);
\draw [->] (A3)  --  (B3);
\draw [->] (A4)  --  (B4);
\draw [->] (B2)  --  (A3);
\draw [->] (B3)  --  (A4);
\draw [->] (B4)  --  (A5);
\draw [->] (A3) -- (A6);
\draw [->] (A6) -- (A4);
\draw[->] (P3)--(P6);\draw[->] (P6)--(P4);

\draw [loosely dotted, very thick] (2,-1) --(P2);\draw [loosely dotted, very thick] (P2) -- (B1);\draw [loosely dotted, very thick] (B1) --(A2); \draw [loosely dotted, very thick] (A2) -- (2,3); 

\draw [loosely dotted, very thick] (8,-1)-- (P5);\draw [loosely dotted, very thick] (P5) --(B5);\draw [loosely dotted, very thick] (B5) --(A5);\draw [loosely dotted, very thick] (A5)--(8,3);
\draw [dotted, thick] (P2) --(P3);
\draw [dotted, thick] (P3) -- (P4);
\draw [dotted, thick] (P4) --(P5); \draw [->] (B1) -- (B2);\draw [dotted, thick] (B2) -- (B3);
\draw [->] (B4) -- (B5);
\draw [dotted, thick] (B3) --(B4);\draw [dotted, thick] (A2)--(A3);
\draw [dotted, thick] (A3) --(A4);\draw [dotted, thick] (A4) -- (A5);
\end{tikzpicture}}
\]
where the vertical dotted lines are identified.
There are 14 different basic maximal rigid objects (that is objects $M$ satisfying $\Ext^1_{\Pi_2(kQ)}(M,M)=0$). They all have 6 direct summands, among which the three projective-injectives.

Therefore, if we work in the stable category $\underline{\mod}\Pi_2(kQ)$ (see \cite{Hap}), where all morphisms factorizing through a projective-injective vanish, we obtain 14 maximal rigid objects with 3 direct summands (indeed the projective-injective objects become isomorphic to zero in $\underline{\mod}\Pi_2(kQ)$). Moreover, the category  $\underline{\mod}\Pi_2(kQ)$ satisfies also the nice property that, given a maximal rigid object, one can replace any of its indecomposable summands by another one in a unique way to obtain another maximal rigid object. 

The category $\underline{\mod}\Pi_2(kQ)$ has the same AR-quiver as the cluster category $\Cc_{Q}$, moreover if we look at the combinatorics of the quivers of the maximal rigid objects in $\underline{\mod}\Pi_2(kQ)$, we obtain:

$$
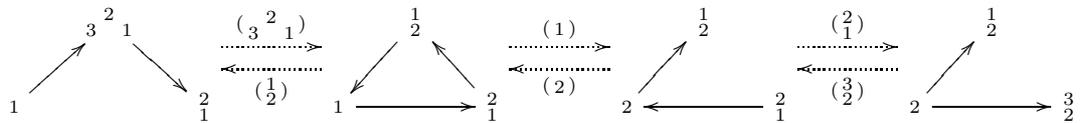
\xymatrix@-.3cm{ & {\bsm&2&\\3&&1\esm} \ar[dr] &&&& {\bsm 1\\2\esm}\ar[dl] &&&& {\bsm1\\2\esm} &&&&  {\bsm1\\2\esm} &\\
{\bsm1\esm} \ar[ur] && {\bsm2\\1\esm} \ar@<.8cm>@{..>}[rr]^{({\bsm&2&\\3&&1\esm})} && {\bsm1\esm}\ar[rr]\ar@<-.5cm>@{..>}[ll]^{({\bsm 1\\2\esm})} && {\bsm2\\1\esm}\ar@<.8cm>@{..>}[rr]^{({\bsm1\esm})}\ar[ul] && \ar@<-.5cm>@{..>}[ll]^{({\bsm2\esm})}{\bsm2\esm}\ar[ur] && {\bsm2\\1\esm} \ar[ll]\ar@<.8cm>@{..>}[rr]^{ ({\bsm2\\1\esm})}&&\ar@<-.5cm>@{..>}[ll]^{({\bsm3\\2\esm})}{\bsm2\esm}\ar[rr]\ar[ur]&&{\bsm3\\2\esm}}.$$

In conclusion, the categories $\Cc_{Q}$ and $\underline{\mod}\Pi_2(kQ)$ seem to be equivalent, and the combinatorics of their maximal rigid objects is closely related to the Fomin-Zelevinsky mutation of quivers.

Our aim in this survey is to show how general these phenomena are, and to construct a large class of categories in which similar phenomena occur.
 
\section{2-Calabi-Yau categories and cluster-tilting theory}
The previous examples are not isolated. These two categories have the same properties, namely they are triangulated $2$-Calabi-Yau categories, and have cluster-tilting objects (which in these cases are the same as maximal rigid). We will see in this section that most of the phenomena that appear in the two cases described in Section \ref{section motivation} still hold in the general setup of $2$-Calabi-Yau categories with cluster-tilting objects. 

\subsection{Iyama-Yoshino mutation}

A triangulated category which is $\Hom$-finite is called \emph{$d$-Calabi-Yau} ($d$-CY for short) if there is a bifunctorial isomorphism $$\Hom_{\Cc}(X,Y)\simeq D\Hom_{\Cc}(Y,X[d])$$ where $D=\Hom_k(-,k)$ is the usual duality over $k$.

\begin{definition}
Let $\Cc$ be a $\Hom$-finite triangulated category. An object $T\in\Cc$ is called \emph{cluster-tilting} (or $2$-cluster-tilting) if $T$ is basic and if we have
$$\add(T)=\{ X\in \Cc \mid \Hom_{\Cc}(X,T[1])=0\}=\{X\in \Cc\mid \Hom_{\Cc}(T,X[1])=0\}.$$
\end{definition}
\noindent
Note that a cluster-tilting object is maximal rigid (the converse is not always true cf. \cite{BIKR08}), and that the second equality in the definition always holds when $\Cc$ is $2$-Calabi-Yau.

\medskip
If there exists a cluster-tilting object in a $2$-CY category $\Cc$, then it is possible to construct others by a recursive process resumed in the following:
\begin{theorem}[Iyama-Yoshino, \cite{IY08}]\label{theorem IY mutation}
Let $\Cc$ be a $\Hom$-finite 2-CY triangulated category with a cluster-tilting object $T$. Let $T_i$ be an indecomposable direct summand of $T\simeq T_i\oplus T_0$. Then there exists a unique indecomposable $T^*_i$ non isomorphic to $ T_i$ such that $T_0\oplus T_i^*$ is cluster-tilting. Morever $T_i$ and $T_i^*$ are linked by the existence of triangles
$$\xymatrix{T_i\ar[r]^u& B\ar[r]^v& T_i^*\ar[r]^w & T_i[1] } \quad\textrm{and}\quad \xymatrix{T_i^*\ar[r]^{u'}& B'\ar[r]^{v'}& T_i\ar[r]^{w'} & T_i^*[1] }$$
where $u$ and $u'$ are minimal left $\add(T_0)$-approximations and $v$ and $v'$ are minimal right $\add(T_0)$-approximations.
\end{theorem}

These exchange triangles permit to mutate cluster-tilting objects and have been first described by Buan, Marsh, Reineke, Reiten and Todorov (see \cite[Proposition~6.9]{BMR+06}) for cluster categories.  The corresponding exchange short exact sequences in module categories over a preprojective algebra of Dynkin type appeared also in the work of Geiss, Leclerc and Schr\"oer (see \cite[Lemma 5.1]{GLS06}). The general statement is due to Iyama and Yoshino \cite[Theorem 5.3]{IY08}. This is why we decided to refer to this mutation as the \emph{IY-mutation} of cluster-tilting objects in this article.

\medskip

This recursive process of mutation of cluster-tilting objects is closely related to the notion of mutation of quivers defined by Fomin and Zelevinsky \cite{FZ1} in the following sense.

\begin{theorem}[Buan-Iyama-Reiten-Scott \cite{BIRS09}]\label{BIRSc}
Let $\Cc$ be a $\Hom$-finite 2-CY triangulated category with cluster-tilting object $T$. Let $T_i$ be an indecomposable direct summand of $T$, and denote by $T'$ the cluster-tilting object $\mu_{T_i}(T)$. Denote by $Q_T$ (resp. $Q_{T'}$) the Gabriel quiver of the endomorphism algebra $\End_{\Cc}(T)$ (resp. $Q_{T'}$).Assume that there are no loops and no $2$-cycles at the vertex $i$ of $Q_T$ (resp. $Q_{T'}$) corresponding to the indecomposable $T_i$ (resp. $T_i^*$). Then we have
$$Q_{T'}=\mu_i(Q_T),$$ where $\mu_i$ is the Fomin-Zelevinsky quiver mutation.
\end{theorem} 
We illustrate this result by the following diagram
\[\xymatrix@C=1.7cm@R=0.7cm{\underset{\tiny{\rm cluster-tilting}}{T}\ar@{<..>}[rr]^{\tiny{\rm IY-mutation}}\ar[d] && \underset{\tiny{\rm cluster-tilting}}{T'}\ar[d]\\ Q_T\ar@{<..>}[rr]^{\tiny{\rm FZ-mutation}} && Q_{T'}  \,.\!\!\!}\]

The corresponding results have been first shown in the setting of cluster categories in \cite{BMR1} and in the setting of preprojective algebras of Dynkin type in \cite{GLS06}.

\subsection{2-Calabi-Yau-tilted algebras}

The endomorphism algebras $\End_\Cc(T)$ where $T$ is a cluster-tilting object are of importance, they are called \emph{$2$-CY tilted algebras} (or \emph{cluster-tilted algebras} if $\Cc$ is a cluster category). The following result says that a $2$-CY category $\Cc$ with cluster-tilting objects is very close to a module category over a $2$-CY-tilted algebra. Such a $2$-CY category can be seen as a ``recollement'' of module categories over $2$-CY tilted algebras.
\begin{proposition}[Keller-Reiten \cite{KR07}]
Let $\Cc$ be a 2-CY triangulated category with a cluster-tilting object $T$. Then the functor $$F_T=\Hom_\Cc(T,-):\Cc\to \mod \End_\Cc(T)$$ induces an equivalence $\Cc/\add(T[1])\simeq \mod \End_\Cc(T)$.

If the objects $T$ and $T'$ are linked by an IY-mutation, then the categories $\mod \End_\Cc(T)$ and $\mod \End_\Cc(T')$ are nearly Morita equivalent, that is there exists a simple $\End_\Cc(T)$-module $S$ and a simple $\End_\Cc(T')$-module $S'$, and an equivalence of categories $\mod \End_\Cc(T)/\add(S)\simeq \mod \End_\Cc(T')/\add(S')$.
\end{proposition}
These results have been first studied in \cite{BMR1} for cluster categories. 
The 2-CY-tilted algebras are Gorenstein (i.e. $\End_\Cc(T)$ has finite injective dimension as a module over itself) and are either of infinite global dimension, or hereditary (see \cite{KR07}). Other nice properties of the functor $\Hom_\Cc(T,-):\Cc\to \mod \End_\Cc(T)$ are studied in \cite{KZ}.

\begin{example} Let $Q$ be the quiver $1\to2\to3$.
 The category $\Cc_Q/\add(T[1])$ where $T=\pi({\bsm 1\esm})\oplus \pi({\bsm 2\\1\esm})\oplus \pi({\bsm 3\\2\\1\esm})$ has the same AR-quiver as $\mod kQ$. (Note that $kQ$ is an hereditary algebra.)

If we take $T=\pi({\bsm 1\esm})\oplus \pi({\bsm 3\esm})\oplus \pi({\bsm 3\\2\\1\esm})$, then the associated 2-CY-titled algebra is given by the quiver 
$
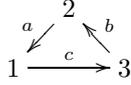
\xymatrix@-0.5cm{ &2\ar[dl]_a& \\1\ar[rr]^c && 3\ar[ul]_b} \quad \textrm{with relations} \quad ab=bc=ca=0.$ It is an algebra of infinite global dimension.
The module category of this algebra has finitely many indecomposables and its AR-quiver is of the form

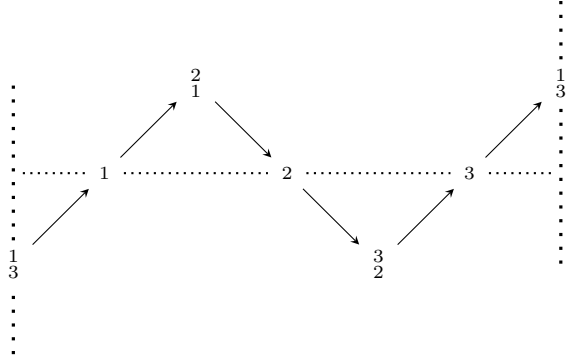
\[\scalebox{1}{
\begin{tikzpicture}[>=stealth,scale=1.2]
\node (P2) at (2,0) {${\bsm 1\\3\esm}$};
\node (P3) at (4,0) {};
\node (P4) at (6,0){${\bsm 3\\2\esm}$};
\node (P5) at (8,0){};

\node (B1) at (2,1) {};
\node (B2) at (3,1) {${\bsm 1\esm}$};
\node (B3) at (5,1) {${\bsm 2\esm}$};
\node (B4) at (7,1) {${\bsm 3\esm}$};
\node (B5) at (8,1) {};
\node (A2) at (2,2) {};
\node (A3) at (4,2) {${\bsm 2\\1\esm}$};
\node (A4) at (6,2){};
\node (A5) at (8,2){${\bsm 1\\3\esm}$};

\draw [->] (P2)  --  (B2);
\draw [->] (P4)  --  (B4);
\draw [->] (B3)  --  (P4);
\draw [->] (A3)  --  (B3);
\draw [->] (B2)  --  (A3);
\draw [->] (B4)  --  (A5);
\draw [loosely dotted, very thick] (2,-1) -- (P2);\draw [loosely dotted, very thick] (P2) -- (2,2); \draw [loosely dotted, very thick] (8,0)--(A5); \draw [loosely dotted, very thick] (A5) -- (8,3);
 \draw [dotted, thick] (B1) -- (B2);\draw [dotted, thick] (B2) -- (B3);
\draw [dotted, thick] (B4) -- (B5);
\draw [dotted, thick] (B3) --(B4);

\end{tikzpicture}}
\]
It is the same quiver as the AR-quiver of $\Cc_Q/\add(T[1])$.
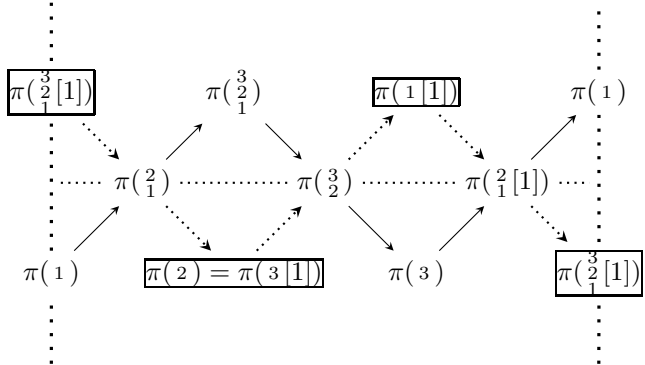
\[\scalebox{1}{
\begin{tikzpicture}[>=stealth,scale=1.2]
\node (P2) at (2,0) {$\pi({\bsm 1\esm})$};
\node (P3) at (4,0) {\frame{$\pi({\bsm 2\esm})=\pi({\bsm3\esm}[1])$}};
\node (P4) at (6,0){$\pi({\bsm 3\esm})$};
\node (P5) at (8,0){\frame{$\pi({\bsm 3\\2\\1\esm}[1])$}};

\node (B1) at (2,1) {};
\node (B2) at (3,1) {$\pi({\bsm 2\\1\esm})$};
\node (B3) at (5,1) {$\pi({\bsm 3\\2\esm})$};
\node (B4) at (7,1) {$\pi({\bsm 2\\1\esm}[1])$};
\node (B5) at (8,1) {};
\node (A2) at (2,2) {\frame{$\pi({\bsm 3\\2\\1\esm}[1])$}};
\node (A3) at (4,2) {$\pi({\bsm 3\\2\\1\esm})$};
\node (A4) at (6,2){\frame{$\pi({\bsm 1\esm}[1])$}};
\node (A5) at (8,2){$\pi({\bsm 1\esm})$};

\draw [->] (P2)  --  (B2);
\draw [->, dotted, thick] (P3)  --  (B3);
\draw [->] (P4)  --  (B4);
\draw [->, dotted, thick] (B2)  --  (P3);
\draw [->] (B3)  --  (P4);
\draw [->, dotted, thick] (B4)  --  (P5);
\draw [->, dotted, thick] (A2)  --  (B2);
\draw [->] (A3)  --  (B3);
\draw [->, dotted, thick] (A4)  --  (B4);
\draw [->] (B2)  --  (A3);
\draw [->, dotted, thick] (B3)  --  (A4);
\draw [->] (B4)  --  (A5);

\draw [loosely dotted, very thick] (2,-1) -- (P2); \draw [loosely dotted, very thick] (P2) -- (A2);\draw [loosely dotted, very thick] (A2) -- (2,3);
 \draw [loosely dotted, very thick] (8,-1)--(P5);\draw [loosely dotted, very thick] (P5) -- (A5);\draw [loosely dotted, very thick] (A5) -- (8,3);

 \draw [dotted, thick] (B1) -- (B2);\draw [dotted, thick] (B2) -- (B3);
\draw [dotted, thick] (B4) -- (B5);
\draw [dotted, thick] (B3) --(B4);

\end{tikzpicture}}
\]

\end{example}

\subsection{Two fundamental families of examples}

The two examples studied in section \ref{section motivation} are part of two fundamental families of 2-CY categories with cluster-tilting objects. We give here the general construction of these families.

\bigskip
\noindent
\textbf{The acyclic cluster category}

The first one is given by the cluster category $\Cc_Q$ associated with an acyclic quiver $Q$. This category was first defined in \cite{BMR+06}, following the decorated representations of acyclic quivers introduced in \cite{MRZ}.  

The \emph{acyclic cluster category} is defined in \cite{BMR+06} as the orbit category of the bounded derived category $\Db(kQ)$ of finite dimensional modules over $kQ$ by the autoequivalence $\SSS_2=\SSS[-2]\simeq \tau[-1]$, where $\SSS$ is the Serre functor $-\lten_{kQ}DkQ$, and $\tau$ is the AR-translation of $\Db(kQ)$.
The objects of this category are the same as those of $\Db(kQ)$ and the spaces of morphisms are given by $$\Hom_{\Cc_Q}(X,Y):=\bigoplus_{i\in \mathbb{Z}}\Hom_{\Db(kQ)}(X,\SSS_2^iY).$$ The canonical functor $\pi:\Db(kQ)\rightarrow \Cc_Q$ satisfies $\pi\circ \SSS_2\simeq \pi $.  Moreover the $k$-category $\Cc_Q$ satisfies the property that, for any $k$-category $\Tt$ and functor $F:\Db(kQ)\to \Tt$ with $F\circ\SSS_2\simeq F$, the functor $F$ factors through $\pi$.   $$\xymatrix@-.3cm{\Db(kQ)\ar[dr]_{\pi}\ar[dd]_{\SSS_2}\ar@/^/[drr]^F &&\\ & \Cc_Q\ar@{..>}[r]& \Tt \\ \Db(kQ)\ar[ur]^{\pi}\ar@/_/[urr]_F &}$$ 
For any $X$ and $Y$ in $\Db(kQ)$, the inifinite sum $\bigoplus_{i\in\ZZ}\Hom_{\Dd}(X,\SSS_2^{i}Y)$ has finitely many non zero summands, hence the category $\Cc_Q$ is $\Hom$-finite. Moreover the shift functor, and the Serre functor of $\Db(kQ)$ yield a shift functor and a Serre functor for $\Cc_Q$. The acyclic cluster category $\Cc_Q$ satisfies the $2$-CY property by construction since the functor $\SSS_2=\SSS[-2]$ becomes isomorphic to the identity in the orbit category $\Db(kQ)/\SSS_2$.

Furthermore, the category $\Cc_Q$ has even more structure as shown in the following fundamental result.

\begin{theorem}[Keller \cite{Kel05}]\label{CQ triangulated}
Let $Q$ be an acyclic quiver. The acyclic cluster category has a natural structure of triangulated category making the functor $\pi:\Db(kQ)\rightarrow \Cc_Q$ a triangle functor. 
\end{theorem}

The triangles in $\Db(kQ)$ yield natural candidates for triangles in $\Cc_Q$, but checking that they satisfy the axioms of triangulated categories (see e.g. \cite{Hap}) is a difficult task. A direct proof (by axioms checking) that $\Cc_Q$ is triangulated is especially difficult for axiom $(TR1)$: Given a morphism in $\Cc_Q$, how to define its cone if it is not liftable into a morphism in $\Db(kQ)$? Keller obtains Theorem \ref{CQ triangulated} using different techniques. He first embeds $\Cc_Q$ in a triangulated category (called triangulated hull), and then shows that this embedding is dense.

\medskip\noindent
\textbf{Link with tilting theory}

If $T$ is a tilting $kQ$-module, then $\pi(T)$ is a cluster-tilting object as shown in \cite{BMR+06}. Moreover the mutation of the cluster-tilting objects is closely related to mutation of tilting $kQ$-modules in the following sense (see \cite{BMR+06}): If $T\simeq T_i\oplus T_0$ is a tilting $kQ$-module with $T_i$ indecomposable and if there exists an indecomposable module $T_i^*$ such that  $T_0\oplus T_i^*$ is a tilting module, then $\pi(T_i^*\oplus T_0)$ is the IY-mutation of $\pi(T)$ at $\pi(T_i)$. The exchange triangles are the images of exchange sequences in the module category. Furthermore, the 2-CY-tilted algebra $\End_\Cc(\pi(T))$ (called \emph{cluster-tilted} in this case) is isomorphic to the trivial extension of the algebra $B=\End_{kQ}(T)$ by the $B$-$B$-bimodule $\Ext^2_B(D(B),B)$ as shown in \cite{ABS08}.

An advantage of cluster-tilting theory over tilting theory is that in cluster-tilting theory it is always possible to mutate. In a  hereditary module category,  an almost complete tilting module does not always have 2 complements (cf.  \cite{HU89}, \cite{RS90}). 

Using strong results of tilting theory due to Happel and Unger \cite{HU03}, one has the following property.

\begin{theorem}[Prop. 3.5 \cite{ BMR+06}]\label{HU theorem}
Let $Q$ be an acyclic quiver. Any cluster-tilting objects $T$ and $T'$ in $\Cc_Q$ are IY-mutation-equivalent, that is one can pass from $T$ to $T'$ using a finite sequence of IY-mutations.
\end{theorem}
This fact is unfortunately not known to be true or false in general 2-CY categories (see aslo Remark \ref{Remark Bruestle}).

The category $\Cc_Q$ provides a first  and very important example which permits to categorify the acyclic cluster algebras. We refer to \cite{BMR+06}, \cite{CC06}, \cite{BMR1}, \cite{BMRT}, \cite{CK1}, \cite{CK2} for results on the categorification of acyclic cluster algebras. Detailed overviews of this subject can be found in \cite{BM06}, \cite{Reitensurvey} and in the first sections of \cite{Kellersurvey}.

\begin{remark}
The construction of $\Cc_Q$ can be generalized  to any hereditary category. Therefore, it is possible to consider $\Cc_{\Hh}=\Db(\Hh)/\SSS_2$, where $\Hh={\sf coh} \mathbb{X}$, and $\mathbb{X}$ is a weighted projective line. This orbit category is still triangulated by \cite{Kel05} (the generalization from acyclic quivers to hereditary categories was first suggested by Asashiba). These cluster categories are studied in details by Barot, Kussin and Lenzing in~\cite{BKL10}.
\end{remark}

\medskip
\noindent
\textbf{Recognition theorem for acyclic cluster categories.}

The following result due to Keller and Reiten ensures that the acyclic cluster categories are the only categories satisfying the properties needed to categorify acyclic cluster algebras.

\begin{theorem}[Keller-Reiten \cite{KR08}]\label{KR recognition theorem}
Let $\Cc$ be an algebraic triangulated $2$-CY category with a cluster-tilting object $T$. If the quiver of the endomorphism algebra $\End_{\Cc}(T)$ is acyclic, then there exists a triangle equivalence $\Cc\simeq \Cc_Q$.
\end{theorem}

Note that this result implies in particular that the cluster-tilting objects of the category $\Cc$ are all mutation equivalent. The equivalence of the two categories described in subsections \ref{subsection CA3} and \ref{subsection PiA3} is a consequence of Theorem \ref{KR recognition theorem}.

\medskip
\noindent
\textbf{Geometric description of the acyclic cluster category.}

Another description of the cluster category has been given independently by Caldero, Chapoton and Schiffler in~\cite{CCS}, \cite{CCS2} in the case where the quiver $Q$ is an orientation of the graph $A_n$. The description is geometric: in this situation, indecomposable objects in $\Cc_Q$ correspond to homotopy classes of arcs joining two non consecutive vertices of an $(n+3)$-gon, and cluster-tilting objects correspond to ideal triangulations of the $(n+3)$-gon. One can associate a quiver to any ideal triangulation $\tau$: the vertices are in bijection with the inner arcs of $\tau$, and two vertices are linked by an arrow if and only if the two corresponding arcs are part of the same triangle. Mutation of cluster-tilting objects corresponds to the flip of an arc in an ideal triangulation.

A similar description of $\Cc_Q$ exists in the case where $Q$ is some (acyclic) orientation of the graphs $D_n$, $\tilde{A}_n$ and $\tilde{D}_n$. In these cases, the surface with marked points is respectively an $n$-gon with a puncture, an annulus without puncture, and an $n$-gon with two punctures (see \cite{Schi1}, \cite{ST}, \cite{Schi2}).

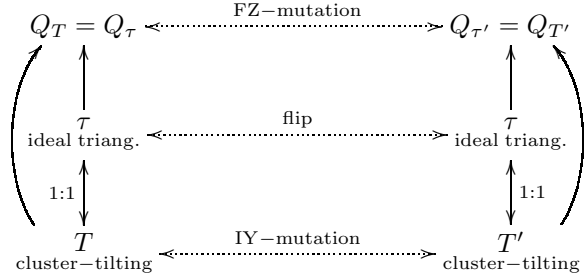
\[\xymatrix@C=1.7cm{Q_T=Q_\tau\ar@{<..>}[rr]^{\tiny{\rm FZ-mutation}} && Q_{\tau'}=Q_{T'}\\ \underset{\tiny{\rm ideal\ triang.}}{\tau} \ar@{<..>}[rr]^{\tiny{\rm flip}}\ar[u]\ar@{<->}[d]_-{1:1}&&\underset{\tiny{\rm ideal\ triang.}}{\tau}\ar[u]\ar@{<->}[d]^-{1:1} \\ \underset{\tiny{\rm cluster-tilting}}{T}\ar@{<..>}[rr]^{\tiny{\rm IY-mutation}} \ar@(l,l)[uu]&& \underset{\tiny{\rm cluster-tilting}}{T'}\ar@(r,r)[uu]   }\]

\begin{example}
Here is the correspondence between cluster-tilting objects of the example in subsection \ref{subsection CA3} and ideal triangulations of the hexagon.
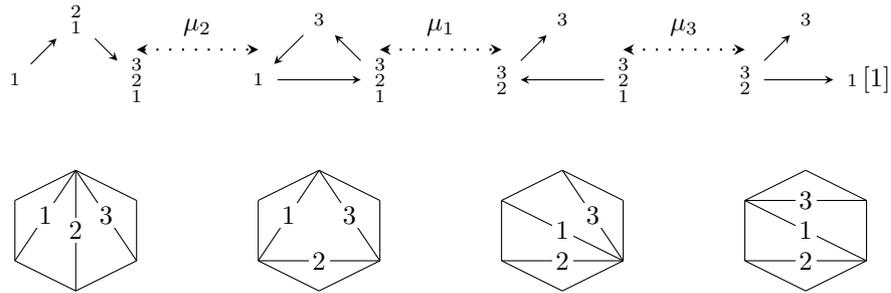
\[\scalebox{1}{
\begin{tikzpicture}[>=stealth,scale=.8]
\node (A1) at (0,0) {${\bsm 1\esm}$};
\node(A2) at (1,1){${\bsm 2\\1\esm}$};
\node (A3) at (2,0){${\bsm 3\\2\\1\esm}$};
\node (B1) at (4,0) {${\bsm 1\esm}$};
\node(B2) at (5,1){${\bsm 3\esm}$};
\node (B3) at (6,0){${\bsm 3\\2\\1\esm}$};
\node (C1) at (8,0) {${\bsm 3\\2\esm}$};
\node(C2) at (9,1){${\bsm 3\esm}$};
\node (C3) at (10,0){${\bsm 3\\2\\1\esm}$};
\node (D1) at (12,0) {${\bsm 3\\2\esm}$};
\node(D2) at (13,1){${\bsm 3\esm}$};
\node (D3) at (14,0){${\bsm 1\esm}[1]$};
\draw[loosely dotted, thick, <->] (2,0.5)--node[swap,yshift=3mm] {$\mu_2$}(4,0.5);
\draw[loosely dotted, thick, <->] (6,0.5)--node[swap,yshift=3mm] {$\mu_1$}(8,0.5);\draw[loosely dotted, thick, <->] (10,0.5)--node[swap,yshift=3mm] {$\mu_3$}(12,0.5);

\draw[->] (A1)--(A2);
\draw[->] (A2)--(A3);
\draw[->] (B2)--(B1);
\draw[->] (B1)--(B3);
\draw[->] (B3)--(B2);
\draw[->] (C3)--(C1);
\draw[->] (C1)--(C2);
\draw[->](D1)--(D2);
\draw[->](D1)--(D3);

\draw (0,-2)--(0,-3);\draw (0,-3)--(1,-3.5);\draw (1,-3.5)--(2,-3);\draw (2,-3)--(2,-2);\draw (2,-2)--(1,-1.5);\draw (1,-1.5)--(0,-2);
\draw (1,-1.5)--node [fill=white,inner sep=.5mm]{$2$}(1,-3.5); \draw (1,-1.5)--node [fill=white,inner sep=.5mm]{$1$}(0,-3); \draw (1,-1.5)--node [fill=white,inner sep=.5mm]{$3$}(2,-3);

\draw (4,-2)--(4,-3);\draw (4,-3)--(5,-3.5);\draw (5,-3.5)--(6,-3);\draw (6,-3)--(6,-2);\draw (6,-2)--(5,-1.5);\draw (5,-1.5)--(4,-2);
\draw (6,-3)--node [fill=white,inner sep=.5mm]{$2$}(4,-3); \draw (5,-1.5)--node [fill=white,inner sep=.5mm]{$1$}(4,-3); \draw (5,-1.5)--node [fill=white,inner sep=.5mm]{$3$}(6,-3);

\draw (8,-2)--(8,-3);\draw (8,-3)--(9,-3.5);\draw (9,-3.5)--(10,-3);\draw (10,-3)--(10,-2);\draw (10,-2)--(9,-1.5);\draw (9,-1.5)--(8,-2);
\draw (10,-3)--node [fill=white,inner sep=.5mm]{$2$}(8,-3); \draw (8,-2)--node [fill=white,inner sep=.5mm]{$1$}(10,-3); \draw (9,-1.5)--node [fill=white,inner sep=.5mm]{$3$}(10,-3);

\draw (12,-2)--(12,-3);\draw (12,-3)--(13,-3.5);\draw (13,-3.5)--(14,-3);\draw (14,-3)--(14,-2);\draw (14,-2)--(13,-1.5);\draw (13,-1.5)--(12,-2);
\draw (14,-3)--node [fill=white,inner sep=.5mm]{$2$}(12,-3); \draw (12,-2)--node [fill=white,inner sep=.5mm]{$1$}(14,-3); \draw (12,-2)--node [fill=white,inner sep=.5mm]{$3$}(14,-2);

\end{tikzpicture}}\]

\end{example}
\bigskip
\noindent
\textbf{Preprojective algebras of Dynkin type}

A very different approach is given by module categories over preprojective algebras. Let $\Delta$ be a finite graph, and $Q$ be an acyclic orientation of $\Delta$. The double quiver $\overline{Q}$ is obtained from $Q$ by adding to each arrow $a:i\rightarrow j$ of $Q$ an arrow $a^*:j\rightarrow i$. The preprojective algebra $\Pi_2(kQ)$ is then defined to be the quotient of the path algebra $kQ$ by the ideal of relations $\sum_{a\in Q_1}aa^*-a^*a$. This algebra is finite dimensional if and only if $Q$ is of Dynkin type, and its module category is closely related to the theory of the Lie algebra of type $Q$. When $Q$ is Dynkin, the algebra $\Pi_2(kQ)$ is selfinjective and the stable category $\underline{\mod}\Pi_2(kQ)$ is triangulated (see \cite{Hap}) and $2$-Calabi-Yau. Moreover it does not depend on the choice of the orientation of $Q$, so we denote it by $\Pi_2(\Delta)$.

In their work \cite{GLS06}, \cite{GLS07a} about cluster algebras arising in Lie theory, Geiss, Leclerc and Schr\"oer have constructed special cluster-tilting objects in the module category $\mod \Pi_2(\Delta)$, which helped them to categorify certain cluster algebras. Further developments of this link between cluster algebras and preprojective algebras can be found in \cite{GLS07}, \cite{GLS08}, \cite{GLS}, \cite{GLS10} (see also \cite{GLSsurvey} for an overview).

\begin{theorem}[Geiss-Leclerc-Schr\"oer \cite{GLS07a}]
Let $\Delta$ be a simply laced Dynkin diagram. Then for any orientation $Q$ of $\Delta$, there exists a cluster-tilting object $T_Q$ in the category $\mod\Pi_2(\Delta)$ so that the Gabriel quiver of the cluster-tilting object $T_Q$ is obtained from the AR quiver of $\mod kQ$ by adding arrows corresponding to the AR-translation. 
\end{theorem}

\begin{example}
In the example of the first section with $\Delta=A_3$, with the three orientations $1\rightarrow 2\rightarrow 3$, $1\leftarrow 2\rightarrow 3$ and $1\rightarrow 2\leftarrow 3$ of $\Delta$, we obtain respectively the following three cluster-tilting objects

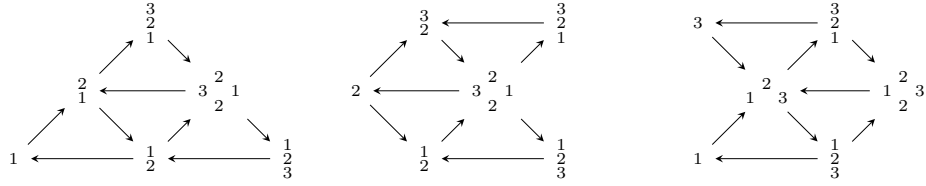
\[\scalebox{.9}{
\begin{tikzpicture}[>=stealth,scale=1]
\node (P2) at (2,0) {${\bsm 1\esm}$};
\node (P3) at (4,0) {${\bsm 1\\2\esm}$};
\node (P4) at (6,0){${\bsm 1\\2\\3\esm}$};
\node (B2) at (3,1) {${\bsm 2\\1\esm}$};
\node (B3) at (5,1) {${\bsm &2&\\3&&1\\&2&\esm}$};
\node (A3) at (4,2) {${\bsm 3\\2\\1\esm}$};

\node (P2') at (8,0) {${\bsm 1\\2\esm}$};
\node (P3') at (10,0) {${\bsm 1\\2\\3\esm}$};
\node (B2') at (7,1) {${\bsm 2\esm}$};
\node (B3') at (9,1) {${\bsm &2&\\3&&1\\&2&\esm}$};
\node (A2') at (8,2) {${\bsm 3\\2\esm}$};
\node (A3') at (10,2){${\bsm 3\\2\\1\esm}$};

\node (P2'') at (12,0) {${\bsm 1\esm}$};
\node (P3'') at (14,0) {${\bsm 1\\2\\3\esm}$};
\node (B2'') at (13,1) {${\bsm &2&\\1&&3\esm}$};
\node (B3'') at (15,1) {${\bsm &2&\\1&&3\\&2&\esm}$};
\node (A2'') at (12,2) {${\bsm 3\esm}$};
\node (A3'') at (14,2){${\bsm 3\\2\\1\esm}$};

\draw [->] (B2')--(P2');\draw [->] (P2')--(B3');
\draw [->] (B2')--(A2');\draw [->] (A2')--(B3');
\draw [->] (B3')--(P3');\draw [->] (B3')--(B2');
\draw [->] (B3')--(A3');\draw [->] (P3')--(P2');\draw [->] (A3')--(A2');

\draw [->] (P2'')--(B2'');\draw [->] (B2'')--(P3'');
\draw [->] (A2'')--(B2'');\draw [->] (B2'')--(A3'');
\draw [->] (P3'')--(B3'');\draw [->] (B3'')--(B2'');
\draw [->] (A3'')--(B3'');\draw [->] (P3'')--(P2'');\draw [->] (A3'')--(A2'');

\draw [->] (P2)  --  (B2);
\draw[->](P3)--(B3);
\draw[->](B2)--(P3);
\draw [->] (B3)  --  (P4);
\draw [->] (B2)  --  (A3);
\draw [->] (B3)  --  (B2);
\draw [->] (P4)  --  (P3);
\draw[->] (P3)--(P2);
\draw[->] (A3)--(B3);

\end{tikzpicture}}
\]

\end{example}

\bigskip\noindent
\textbf{Correspondence}

The link between cluster categories and stable module categories of preprojective algebras of Dynkin type is given in the following.
\begin{theorem} Let $\Delta$ be a simply laced Dynkin diagram and $Q$ be an acyclic quiver. There is a triangle equivalence $\underline{\mod}\Pi_2(\Delta)\simeq \Cc_Q$ if and only if  we are in one of the cases 
  $\Delta=A_2$ and $Q$ is of type $A_1$,
  $\Delta=A_3$ and $Q$ is of type $A_3$,
 $\Delta=A_4$ and $Q$ is of type $D_6$.

\end{theorem}

\begin{remark}
In the case $\Delta=A_5$, the category $\underline{\mod}\Pi_2(A_5)$ is equivalent to the cluster category $\Cc_{\Hh}$ where $\Hh$ is of tubular type $E^{(1,1)}_8$.
\end{remark}

Even though the categories $\underline{\mod}\Pi_2(\Delta)$ and $\Cc_Q$ are constructed in a very distinct way, this theorem shows that these categories are not so different in a certain sense. This observation allows us to ask the following.

\begin{question}\label{problem 1}
Is it possible to generalize the constructions above so that  the $2$-Calabi-Yau categories $\underline{\mod}\Pi_2(\Delta)$ and $\Cc_Q$ are part of the same family of categories?
\end{question}

There are two different approaches to this problem. One consists in generalizing the construction $\underline{\mod}\Pi_2(\Delta)$ and the other consists in generalizing the construction of the cluster category. A construction of the first type is given in \cite{GLS07} and in \cite{BIRS09} (see subsection \ref{subsection words})  and a construction of the second type is given in \cite{Ami08} and \cite{Ami09} (see section \ref{section general construction}).

\subsection{Quivers with potential}

Theorem \ref{BIRSc} links the quivers of the $2$-CY-tilted algebras appearing in the same IY-mutation class of cluster-tilting object. This leads to the natural question: 

\medskip
\noindent
\emph{Is there a combinatorial way to deduce the relations of the $2$-CY-tilted algebra after mutation?}
\medskip

This question (among others) brought Derksen Weyman and Zelevinsky to introduce in \cite{DWZ} the notions of quiver with potential (QP for short),  Jacobian algebra and mutation of QP.

\begin{definition}
A \emph{potential} $W$ on a quiver $Q$ is an element in $\widehat{kQ}/\widehat{[kQ,kQ]}$ where $\widehat{kQ}$ is the completion of the path algebra $kQ$ for the $J$-adic topology  ($J$ being the ideal of $kQ$ generated by the arrows) where $\widehat{[kQ,kQ]}$ is the subspace of $\widehat{kQ}$ generated by the commutators of the algebra $kQ$. 
\end{definition}
\noindent
In other words, a potential is a (possibly infinite) linear combination of cycles of $Q$, up to cyclic equivalence ($a_1a_2\ldots a_n\sim a_2a_3\ldots a_na_1$).

\begin{definition}
Let $Q$ be a quiver. The partial derivative $\partial:\widehat{kQ}/\widehat{[kQ,kQ]}\rightarrow \widehat{kQ}$ is defined to be the unique continuous linear map which sends the class of a path $p$ to the sum $\sum_{p=uav}vu$ taken over all decompositions of the path $p$.

Let $(Q,W)$ be a quiver with potential. The \emph{Jacobian algebra} of $(Q,W)$ is defined to be $\Jac(Q,W):=\widehat{kQ}/\langle \partial_a W, a\in Q_1\rangle$, where $\langle \partial_a W, a \in Q_1\rangle$ is the ideal of $\widehat{kQ}$ generated by  $\partial_aW$ for all $a\in Q_1$.

A QP $(Q,W)$ is  \emph{Jacobi-finite} if its Jacobian algebra is finite dimensional.
\end{definition}

In \cite{DWZ} the authors introduced the notion of \emph{reduction} of a QP. The reduction of $(Q,W)$ consists in finding a QP $(Q^{\rm r},W^{\rm r})$ whose key properties are that the Jacobian algebras $\Jac(Q,W)$ and $\Jac(Q^{\rm r},W^{\rm r})$ are isomorphic and $W^{\rm r}$ has no 2-cycles as summands. In the case of a Jacobi-finite QP, the quiver $Q^{\rm r}$ is the Gabriel quiver of the Jacobian algebra and is then uniquely determined.  The potential $W^{\rm r}$ is not uniquely determined. The notion of reduction is defined up to an equivalence relation of QPs called \emph{right equivalence} (see \cite{DWZ}).

Derksen, Weyman and Zelevinsky also refined the notion of the FZ-mutation of a quiver into the notion of mutation of a QP (that we will call DWZ-mutation) at a vertex of $Q$ without loop. Steps (M1) and (M2) are the same as in FZ-mutation. The new potential is defined to be a sum $[W]+W^*$, constructed from $W$ using the new arrows. Step (M3') consists in reducing the QP we obtained. For generic QPs, step (M3') of the DWZ-mutation coincides with step (M3) of the FZ-mutation when restricted to the quiver. Let us illustrate this process in an example. 
\begin{example}
Let us define $(Q,W)$  as follows: 
$$\xymatrix@-0.3cm{& &2\ar[dr]^b && \\ Q=&1\ar[rr]^c\ar[ur]^a&&3\ar[dl]^d&\quad \textrm{ with }W=edc.\\ &&4\ar[ul]^e&& } $$ The Jacobian ideal is defined by the relations $ed=ce=dc=0$. One easily checks that $(Q,W)$ is Jacobi-finite.

Let us mutate $(Q,W)$ at the vertex $3$. After steps $(M1)$ and $(M2)$ we obtain the QP $(Q',W')$ defined as follows:
$$
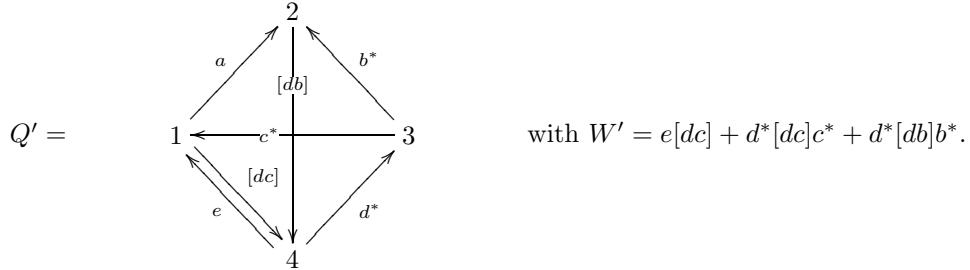
\xymatrix@+0.3cm{ &&2\ar[dd]|(.3){[db]} & &\\ Q'=&1\ar@<.5ex>[dr]^{[dc]}\ar[ur]^a&&3\ar[ll]|(.6){c^*}\ar[ul]_{b^*}&\textrm{ with }W'=e[dc]+d^*[dc]c^*+d^*[db]b^*.\\ &&4\ar@<.5ex>[ul]^e\ar[ur]_{d^*}& &} 
$$
 Ê
After Step (M3') we obtain the QP:
$$
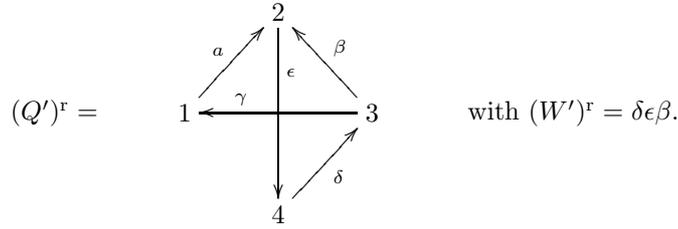
\xymatrix{& &2\ar[dd]^(.3){\epsilon} & &\\ (Q')^{\rm r}=&1\ar[ur]^a&&3\ar[ll]_(.7){\gamma}\ar[ul]_{\beta}&\textrm{ with }(W')^{\rm r}=\delta\epsilon\beta.\\ &&4\ar[ur]_{\delta}&& } 
$$ The map $kQ'\to k(Q')^{\rm r}$ sending $a$, $b^*$, $c^*$, $d^*$ $e$, $[dc]$, $[db]$ on respectively $\alpha$, $\beta$, $\gamma$, $\delta$, $-\gamma\delta$, $0$, $\epsilon$, sends $W'$ on $(W')^{\rm r}$ and induces an isomorphism between $\Jac(Q',W')$ and $\Jac((Q')^{\rm r},(W')^{\rm r})$.

The mutation of $(Q,W)$ at $3$ is defined to be $((Q')^{\rm r},(W')^{\rm r})$ (up to right equivalence).

\end{example}

If there do not occur $2$-cycles at any iterate mutation of $(Q,W)$, then the potential is \emph{non-degenerate}  \cite{DWZ} and the DWZ-mutation (when restricted to the quiver) coincides with the FZ-mutation of quivers.  This happens in particular when the potential is  \emph{rigid}, that is when all cycles of $Q$ are cyclically equivalent to an element in the Jacobian ideal $\langle \partial_aW, a\in Q_1\rangle$. This notion of rigidity is stable under mutation \cite{DWZ}. 

Moreover the following theorem gives an answer to the previous question in the case where the $2$-CY-tilted algebra is Jacobian.

\begin{theorem}[Buan-Iyama-Reiten-Smith \cite{BIRSm}]\label{BIRSm}
Let $T$ be a cluster-tilting object in an $\Hom$-finite $2$-CY category $\Cc$ and let $T_i$ be an indecomposable direct summand of $T$. Assume that there is an algebra isomorphism $\End_{\Cc}(T)\simeq \Jac(Q,W)$, for some QP $(Q,W)$ and assume that there are no loops nor $2$-cycles at vertex $i$ (corresponding to $T_i$) in the quiver of $\End_{\Cc}(T)$ and a technical assumption called \emph{glueing condition}. Then there is an isomorphism of algebras $$\Jac(\mu_i(Q,W))\simeq \End_{\Cc}(\mu_{T_i}(T)).$$
\end{theorem}

\[\xymatrix@C=1.7cm{\underset{\tiny{\rm cluster-tilting}}{T}\ar@{<..>}[rr]^{\tiny{\rm IY-mutation}}\ar@{|->}[d] && \underset{\tiny{\rm cluster-tilting}}{T'}\ar@{|->}[d]\\ \End_\Cc(T)\simeq \Jac(Q,W) && \End_\Cc(T')\simeq \Jac(Q',W')\\ (Q,W)\ar@{<..>}[rr]^{\tiny{\rm DWZ-mutation}}\ar@{|->}[u] && (Q',W')\ar@{|->}[u]     }\]

This result applies in particular for acyclic cluster categories. If $Q$ is an acyclic quiver, the endomorphism algebra of the canonical cluster-tilting object $\pi(kQ)\in\Cc_Q$ is isomorphic to the Jacobian algebra $\Jac(Q,0)$ and $0$ is a rigid potential for the quiver $Q$. Moreover the categories $\Cc_Q$ satisfy the glueing condition \cite{BIRSm}. Therefore, combining Theorems \ref{BIRSm} and \ref{HU theorem}, all cluster-tilted algebras are Jacobian algebras associated with a rigid QP and can be deduced one from each other by mutation.

The endomorphism algebras of the canonical cluster-tilting object of the category $\underline{\mod }\Pi_2(\Delta)$, where $\Delta$ is a Dynkin graph, are also Jacobian algebras with rigid QP and satisfy the glueing condition \cite[Thm 6.5]{BIRSm}.

The two natural questions are now the following:

\begin{question}\label{problem 2}
\begin{itemize}
\item[1.] Are all $2$-CY-tilted algebras Jacobian algebras?
\item[2.] Are all Jacobian algebras $2$-CY-tilted algebras?
\end{itemize}
\end{question}

The first question is still open, but the answer is expected to be negative by recent results of Davison \cite{Davison} and Van den Bergh \cite{VdB2}, which pointed out the existence of algebras which are bimodule 3-CY but which do not come from a potential. 
The second one has a positive answer in the Jacobi-finite case (see subsection \ref{subsection ginzburg}).

\medskip
Derksen Weyman and Zelevinsky also describe the mutation of (decorated) representation of Jacobian algebras in \cite{DWZ}. Their aim was to generalize the construction of mutations of decorated representations \cite{MRZ}. Given a $\Jac(Q,W)$-module $M$, they associate a $\Jac(\mu_i(Q,W))$-module $\mu_i(M)$. This mutation preserves indecomposability. 

\section{From $3$-Calabi-Yau DG-algebras to 2-Calabi-Yau categories}\label{section general construction}

\subsection{Graded algebras and DG algebras}
In this section, we recall some definitions concerning differential graded algebras, differential graded modules, and derived categories. We refer to \cite{Kel94} for definitions and properties on differential graded algebras and associated triangulated categories.

\medskip
\noindent
\textbf{Graded algebras}

We denote by $\Gr k$ the category of graded $k$-vector spaces. Recall that given two $\ZZ$-graded vector spaces $M=\bigoplus_{p\in \ZZ} M_p$ and $N=\bigoplus_{p\in\ZZ}N_p$, a morphism $f\in\Hom_{\Gr k}(M,N)$ is of the form $f=\bigoplus_{p\in \ZZ}(f_p)$ where $f_p\in \Hom_k(M_p,N_p)$ for any $p\in \ZZ$.

Let $A=\bigoplus_{i\in \ZZ}A_i$ be a $\mathbb{Z}$-graded algebra. A $\ZZ$-graded $A$-module $M$ is a $\ZZ$-graded $k$-vector space $M=\bigoplus_{i\in \ZZ}M_i$ with a morphism of $\mathbb{Z}$-graded algebras $A\to \bigoplus_{p\in \ZZ}\Hom_{\Gr k}(M,M(p))$, where $M(p)$ is the $\ZZ$-graded $k$-module such that $M(p)_i=M_{p+i}$ for any $i\in \ZZ$. In other words, for any $n,p\in\mathbb{Z}$ there is a $k$-linear map $M_n\times A_p\to M_{n+p}$ sending $(m,a)$ to $m.a$ such that $m.1=m$ for all $m\in M$ and such that the following diagram commutes
$$\xymatrix{M_n\times A_p\times A_q \ar[r]\ar[d] & M_n\times A_{p+q}\ar[d]\\M_{n+p}\times A_q \ar[r] & M_{n+p+q},}$$
where the maps are induced by the multiplication in $A$ for the first row and by the action of $A$ on $M$ respectively for the others. Then the \emph{degree shift} $M(1)$ of a graded module $M$ is still a graded $A$-module. 
A \emph{morphism} of graded $A$-modules is a morphism $f:M\to N$ homogeneous of degree $0$, that is $f=\bigoplus_{n\in \ZZ}f_n$ where $f_n\in\Hom_k(M_n, N_n)$ and which commutes with the action of $A$, that is for any $a\in A_p$ there is a commutative diagram
$$\xymatrix{M_n\ar[r]^{f_n}\ar[d]_{.a} & N_n\ar[d]^{.a}\\ M_{n+p}\ar[r]^{f_{n+p}} & N_{n+p}.}$$

The category of graded $A$-modules $\Gr A$ is an abelian category. Therefore we can define the derived category of graded $A$-modules $\Dd(\Gr A)$ as usual (see \cite{Kel07}). For a complex of graded $A$-modules, we denote by $M[1]$ the complex such that $M[1]^n=M^{n+1}$ and $d_{M[1]}=-d_M$.

\begin{definition}
Let $d\geq 2$. We say that $A$ is \emph{bimodule $d$-Calabi-Yau of Gorenstein parameter 1} if 
there exists an isomorphism $$\RHom_{A^{\rm e}}(A,A^{\rm e})[d](-1)\simeq A \quad\rm{in }\ \Dd(\Gr A^{\rm e})$$ where $A^{\rm e}$ is the graded algebra $A^{\rm op}\otimes A$. 
\end{definition}

\medskip
\noindent
\textbf{DG algebras}
\begin{definition}
A \emph{differential graded algebra} (=DG algebra for short) $A$ is a $\mathbb{Z}$-graded $k$-algebra $A=\bigoplus_{n\in \mathbb{Z}}A^n$ with a differential $d_A$, that is a $k$-endomorphism of $A$ homogeneous of degree $1$ satisfying the Leibniz rule: for any $a\in A^p$ and any $b\in A$ we have $d_A(ab)=d_A(a)b+(-1)^pad_A(b).$
\end{definition} 

A \emph{DG $A$-module} $M$ is a graded $A$-module $M=\bigoplus_{n\in\ZZ}M^n$, endowed with a differential $d_M$, that is    $d_M\in \Hom_{\Gr A}(M,M(1))$ such that $d_M\circ d_M=0$.
Moreover the differential $d_M$ of the complex $M$ is compatible with the differential of $A$ in the following sense
$$ \forall m\in M^n,\ \forall a\in A^p,\ d_M(m.a)=m.d_A(a) +(-1)^n d_M(m).a$$
$$\xymatrix@-.5cm{M^n\times A^p\ar[rr]^{(d_M,1)}\ar[dd]_{(1,d_A)}\ar[dr] && M^{n+1}\times A^p\ar[dd]\\ & M^{n+p}\ar[dr]^-{d_M}& \\ M^n\times A^{p+1}\ar[rr]&& M^{n+p+1}}.$$
Note that this equality for $p=0$ implies that $M$ has a structure of complex of $Z^0(A)$-modules.
A morphism of DG $A$-modules is a morphism of graded $A$-modules which is a morphism of complexes. Note that if we endow the graded $A$-module $M(1)$ with the differential $d_{M(1)}=-d_M$, it is a DG $A$-module and there is a canonical isomorphism of DG $A$-modules $M(1)\simeq M[1]$.

The \emph{derived category} $\Dd(A)$ is defined as follows. The objects are DG $A$-modules, and morphisms are equivalence classes  of diagrams  
$$s^{-1}f:\xymatrix@-.5cm{M\ar[dr]_-f&&N\ar[dl]^-s \\ & N' &}$$ where $f$ is a morphism of DG $A$-modules, and $s$ is a morphism of DG $A$-modules such that for any $n\in \mathbb{Z}$, the morphism $H^n(s):H^n(N)\to H^n(N')$ is an isomorphism of $H^0(A)$-module ($s$ is a quasi-isomorphism).
This is a triangulated category. We denote by $\per A$ the thick subcategory (= the smallest triangulated category stable under direct summands) generated by $A$. We denote by $\Db(A)$ the subcategory of $\Dd(A)$ whose objects are the DG $A$-modules with finite dimensional total cohomology.
\begin{remark}
If $A$ is a $k$-algebra, we can view it as a DG algebra concentrated in degree $0$, and with differential $0$. Then we recover the usual derived categories $\Dd(A)$, $\per A$ and $\Db(A)$.
\end{remark}

\begin{definition}
A DG algebra is \emph{bimodule $d$-Calabi-Yau}, if there exists an isomorphism $$\RHom_{A^{\rm e}}(A,A^{\rm e})\simeq A[d] \quad\rm{in}\ \Dd(A^{\rm e})$$ where $A^{\rm e}$ is the DG algebra $A^{\rm op}\otimes A$.
\end{definition}
The next proposition says that the bimodule $d$-Calabi-Yau property for $A$ implies the $d$-CY property for the bounded derived category $\Db(A)$.
\begin{proposition}[\cite{G06, Kellersurvey}]
If $A$ is a DG algebra which is bimodule $d$-Calabi-Yau, then there exists a functorial isomorphism
$$\Hom_{\Dd A}(X,Y)\simeq D\Hom_{\Dd A}(Y, X[d]), \quad \textrm{for any } X\in \Dd A\ \textrm{and } Y\in \Db(A).$$
\end{proposition}
\begin{remark}
If $A$ is a graded algebra, we can view it as a DG algebra with differential $0$. Then if $X$ is a complex of graded $A$-modules, for any $n\in \ZZ$, the differential $d_X:X^n=\bigoplus_{i\in\mathbb{Z}}X^n_i\to X^{n+1}=\bigoplus_{i\in\ZZ}X^{n+1}_i$ is homogeneous of degree $0$, it goes from $X^n_i$ to $X^{n+1}_i$. If we set $M^n:=\bigoplus_{i\in \ZZ}X^{n-i}_i$, then the complex $(M, d_X)$ has naturally a structure of DG $A$-module. This  yields a canonical functor $\Dd (\Gr A)\to \Dd A$, which induces a fully faithful functor $\per\Gr A/(1)\to\per A$. The category $\per(A)$ is generated, as a triangulated category, by the image of the orbit category $\per (\Gr A)/(1)$ through this functor. More precisely, $\per A$ is the triangulated hull of the orbit category $\per (\Gr A)/(1)$.

If the graded algebra $A$ is bimodule $d$-Calabi-Yau of Gorenstein parameter $1$, then the DG algebra $A$ endowed with the zero differential is bimodule $d$-Calabi-Yau.
\end{remark}

\subsection{General construction}
The next result is the main step in the construction of new $2$-CY categories with cluster-tilting object which generalize the acyclic cluster categories.
\begin{theorem}[Thm 2.1 \cite{Ami09}]\label{theorem general}
Let $\BPi$ be a DG-algebra with the following properties:
\begin{itemize}
\item[(a)] $\BPi$ is homologically smooth (i.e. $\BPi\in\per(\BPi^{\rm e})$),
\item[(b)] $H^p(\BPi)=0$ for all $p\geq 1$,
\item[(c)] $H^0(\BPi)$ is finite dimensional as $k$-vector space,
\item[(d)] $\BPi$ is bimodule $3$-CY.
\end{itemize}
Then the triangulated category $\Cc(\BPi)=\per\BPi/\Db(\BPi)$ is $\Hom$-finite, $2$-CY and the object $\BPi$ is a cluster-tilting object with $\End_{\Cc}(\BPi)\simeq H^0(\BPi)$.
\end{theorem}

Let us make a few comments on the hypotheses and on the proof of this theorem.
Hypothesis (a) implies that $\Db(\BPi)\subset \per \BPi$ (\cite{Kellersurvey}), so that the category $\Cc(\BPi)=\per\BPi/\Db(\BPi)$ exists. Hypothesis (b) implies the existence of a natural $t$-structure (coming from the usual truncation) on $\per\BPi$, which is an essential ingredient for the proof. Moreover, hypotheses (c) and (d) imply that $H^p(\BPi)$ is finite dimensional over $k$ for all $p\in \ZZ$. 
Thus for any $X\in \per\BPi$, there is a triangle $$\xymatrix{\tau_{\leq n}X \ar[r] & X\ar[r] & \tau_{>n}X\ar[r] &\tau_{\leq n}X[1]},$$ where $\tau_{> n}X$ is in $\Db(\BPi)$. Moreover we have $$\Hom_{\Cc}(X,Y)\simeq \lim_{n\rightarrow \infty}\Hom_{\per\BPi}(\tau_{\le n}X, \tau_{\le n}Y),$$ which implies the $\Hom$-finiteness of the category $\Cc$.

Now we define a full subcategory $\Ff(\BPi)$ of $\per\BPi$ by $$\Ff(\BPi)=(\per\BPi)_{\leq 0}\cap ((\per\BPi)_{\geq -2})^\perp,$$ where $(\per\BPi)_{\leq p}$ (resp. $(\per\BPi)_{\geq p}$) is the full subcategory of $\per \BPi$ consisting of objects having their homology concentrated in degrees $\leq p$ (resp. $\geq p$). An important step of the proof consists in showing that the composition 
$$\xymatrix{\Ff(\BPi)\ar@{^(->}[r]& \per\BPi\ar[r] & \per\BPi/\Db\BPi =\Cc(\BPi)}$$ is an equivalence. Notice that the subcategory $\Ff(\BPi)$ is not stable under the shift functor.

This equivalence implies in particular the following. 
\begin{proposition}\label{prop these}
Let $\BPi$ be as in Theorem \ref{theorem general}.  Then the following diagram is commutative:
$$\xymatrix{\per\BPi\supset\Ff (\BPi)\ar[dr]_{H^0}\ar[rr]^{\sim} && \Cc(\BPi)\ar[dl]^{\Hom_\Cc(\BPi,-)=F_{\BPi}}\\ & \mod H^0(\BPi)& }$$
\end{proposition}

\subsection{Ginzburg DG algebras}\label{subsection ginzburg}
The general theorem above applies to Ginzburg DG algebras associated with a Jacobi-finite QPs.

\begin{definition}[Ginzburg \cite{G06}]
Let $(Q,W)$ be a QP. Let $Q^G$ be the graded quiver with the same set of vertices as $Q$ and whose arrows are:
\begin{itemize}
\item the arrows of $Q$ (of degree $0$);
\item an arrow $a^*:j\rightarrow i$ of degree $-1$ for each arrow $a:i\rightarrow j$ of $Q$;
\item a loop $t_i:i\rightarrow i$ of degree $-2$ for each vertex $i\in Q_0$.
\end{itemize}
The \emph{completed Ginzburg DG algebra} $\widehat{\Gamma}(Q,W)$ is the DG algebra whose underlying graded algebra is the completion (for the $J$-adic topology) of the graded path algebra $kQ^G$. The differential of $\widehat{\Gamma}(Q,W)$ is the unique continuous linear endomorphism homogeneous of degree $1$ which satisfies the Leibniz rule
(i.e. $d(uv) = (du)v + (-1)^p udv$  for all homogeneous $u$ of degree $p$ and all $v$), and takes the following values on the arrows of $Q^G$:
$$d(a)=0 \textrm{ and }d(a^*)=\partial_a W \ \forall a\in Q_1;\quad d(t_i)=e_i(\sum_{a\in Q_1}[a,a^*])e_i \ \forall i\in Q_0.$$
\end{definition} 

\begin{theorem}[Keller \cite{Kel09}]
The completed Ginzburg DG algebra $\widehat{\Gamma}(Q,W)$ is homologically smooth and bimodule $3$-Calabi-Yau.
\end{theorem}
It is immediate to see that $\widehat{\Gamma}(Q,W)$ is non zero only in negative degrees, and that $H^0(\widehat{\Gamma}(Q,W))\simeq \Jac(Q,W)$. Therefore by Theorem \ref{theorem general} we get the following.
\begin{corollary}\label{Jacobi finite cluster categories}
Let $(Q,W)$ be a Jacobi-finite QP. Then the category $$\Cc_{(Q,W)}:=\per \widehat{\Gamma}(Q,W)/\Db(\widehat{\Gamma}(Q,W))$$ is $\Hom$-finite, $2$-Calabi-Yau, and has a canonical cluster-tilting object whose endomorphism algebra is isomorphic to $\Jac(Q,W)$.
\end{corollary}

This category $\Cc_{(Q,W)}$ is called the \emph{cluster category associated with a QP}. This corollary gives in particular an answer to Question \ref{problem 2} (2) in the Jacobi-finite case.

\begin{remark}

If $(Q,W)$ is not Jacobi-finite, a generalization of the category $\Cc_{(Q,W)}$, which is not $\Hom$-finite, is constructed in \cite{Pla10}.

Note that in a recent paper \cite{VdB}, Van den Bergh has shown that complete DG algebras in negative degrees which are bimodule $3$-Calabi-Yau are quasi-isomorphic to deformations of Ginzburg DG algebras in the sense of \cite{Kel09}. Hence a complete DG algebra $\BPi$ as in Theorem \ref{theorem general} is a (deformation) of a Ginzburg DG algebra associated with a Jacobi-finite QP. 
\end{remark}

The following two results give a link between Ginzburg DG algebras associated with QPs linked by mutation.

\begin{theorem}\label{KY theorem}
Let $(Q,W)$ be a QP without loops and $i\in Q_0$ not on a $2$-cycle in $Q$. Denote by $\Gamma:= \widehat{\Gamma}(Q,W)$ and $\Gamma':= \widehat{\Gamma}(\mu_i(Q,W))$ the completed Ginzburg DG algebras.
\begin{itemize}
\item[(a)]\cite{KY09}, \cite{Kel09} 
There are triangle equivalences
$$\xymatrix{\per \Gamma\ar[r]^-\sim& \per \Gamma'\\ \Db \Gamma\ar[r]^-\sim\ar@{^(->}[u] & \Db\Gamma . \ar@{^(->}[u]}$$
Hence we have a triangle equivalence $\Cc(Q,W)\simeq \Cc(\mu_i(Q,W))$.
\item[(b)] \cite{Pla10}
We have a diagram
$$\xymatrix@C=1.7cm{\per \Gamma \ar[rr]^\sim \ar[d]^{H^0}&& \per \Gamma'\ar[d]^{H^0}\\ \mod \Jac(Q,W) \ar@{<..>}[rr]^{\rm DWZ-mutation}_{\rm for \ representations} && \mod\Jac(\mu_i(Q,W)) .}$$
\end{itemize}
\end{theorem}
Combining (b) together with Proposition \ref{prop these}, we obtain that in the Jacobi-finite case, for any cluster-tilting object $T\in \Cc_{(Q,W)}$ which is IY-mutation equivalent to the canonical one, we have:
\begin{equation}\label{equation diagram}\xymatrix@C=1.7cm{ T\in \ar@(ur,ul)@{<..>}[rr]^{\rm IY-mutation}&\Cc_{(Q,W)}\ar[dl]_{F_{T'}}\ar[dr]^{F_{T}} &\ni T'\\ \mod \End_\Cc(T) \ar@{<..>}[rr]^{\rm DWZ-mutation}_{\rm for \ representations} && \mod\End_\Cc(T').}\end{equation}

This is not clear from the definition wether the cluster categories $\Cc_{(Q,W)}$ satisfy the glueing condition of Theorem \ref{BIRSm}. However Theorem \ref{KY theorem}(a) ensures that if $(Q,W)$ is a non degenerate Jacobi-finite QP, then  IY-mutation is compatible with DWZ-mutation for cluster-tilting objects mutation equivalent to the canonical one $\Gamma (Q,W)$.

\subsection{Application to surfaces with marked points}
Let $(\Sigma,M)$ be a pair consisting of a compact Riemann surface $\Sigma$ with non-empty boundary and a set M of marked points of $\Sigma$, with at least one marked point on each boundary component.  By an arc on $\Sigma$, we mean the homotopy class of a non crossing simple curve on $\Sigma$ with endpoints in $M$ (which may coincide), which does not cut out an unpunctured monogon, or an unpunctured digon.  

To each ideal triangulation $\tau$ of the surface $(\Sigma,M)$, Labardini associates in \cite{Lab1} a QP $(Q_\tau,W_\tau)$ which is rigid and Jacobi-finite. Moreover he shows that the flip of the triangulation coincides with the DWZ-mutation of the quiver with potential. 
Since any two triangulations are linked by a finite sequence of flips, the generalized cluster categories obtained from the QPs associated with the triangulations are all equivalent. More precisely, combining results in \cite{Lab1} with Corollary~\ref{Jacobi finite cluster categories} and Theorem~\ref{KY theorem} (a) we get the following.

\begin{corollary}\label{corollary labardini}
Let $(\Sigma, M)$ be a surface with marked points with non empty boundary. Then there exists a $\Hom$-finite triangulated 2-CY category $\Cc_{(\Sigma,M)}$ with a cluster-tilting object $T_\tau$ corresponding to each ideal triangulation $\tau$ such that we have the following commutative diagram 
\[\xymatrix@C=1.7cm@R=0.7cm{\underset{\tiny{\rm triangulation}}{\tau}\ar@{<..>}[rr]^{\tiny{\rm flip}}\ar@{|->}[d] && \underset{\tiny{\rm triangulation}}{\tau'}\ar@{|->}[d]\\ \underset{\tiny{\rm cluster-tilting}}{T_{\tau}}\ar@{<..>}[rr]^{\tiny{\rm IY-mutation}} && \underset{\tiny{\rm cluster-tilting}\,.\!\!\!}{T_{\tau'}} }\]
All the cluster-tilting objects $T_\tau$ for $\tau$ an ideal triangulation are IY-mutation equivalent.
\end{corollary} 

In \cite{Lab2} Labardini moreover associates to each arc $j$ on the surface  a  module over  the Jacobian algebra $\Jac(Q_\tau,W_\tau)$ for each triangulation $\tau$, in a way compatible with DWZ-mutation for representations. More precisely, if $X_\tau(j)$ (resp. $X_{\tau'}(j))$) is the $\Jac(Q_\tau,W_\tau)$-module (resp. $\Jac(Q_{\tau'},W_{\tau'})$-module), where $\tau'$ is a flip of $\tau$, then one can pass from $X_\tau(j)$ to $X_{\tau'}(j)$ using the DWZ-mutation for representations on the corresponding vertex:
\[\xymatrix@C=1.7cm@R=0.7cm{\underset{\tiny{\rm triangulation}}{\tau}\ar@{<..>}[rr]^{\tiny{\rm flip}}\ar@{|->}[d] && \underset{\tiny{\rm triangulation}}{\tau'}\ar@{|->}[d]\\ X_{\tau}(j)\ar@{<..>}[rr]^{\tiny{\rm DWZ-mutation}}_{\tiny{\rm for \ representations}} && X_{\tau'}(j) .\!  }\]

Denote by $IY(\Sigma,M)\subset \Cc_{(\Sigma,M)}$ the set of all cluster-tilting objects of $\Cc_{(\Sigma,M)}$ which are in the IY-mutation class of a cluster-tilting object  $T_\tau$, where $\tau$ is an ideal triangulation (it does not depend on the choice of $\tau$ by Corollary \ref{corollary labardini}). If a triangulation $\tau$ contains a self-folded triangle, it is not possible to flip $\tau$ at the inside radius of the self-folded triangle. To avoid this problem, Fomin, Shapiro and Thurston introduced in \cite{FST08} the notions of \emph{tagged arcs} and \emph{tagged triangulations} which generalize the notions of arcs and triangulation. Then it is possible to mutate any tagged arc of any tagged triangulations. Hence we obtain a bijection:
$$\xymatrix{IY(\Sigma,M)\ar@{<..>}[rr]^-{1:1} && \{ \textrm{tagged triangulations of }(\Sigma, M) \}}.$$

Let $\ind IY(\Sigma,M)$ be the set of all indecomposable summands of $IY(\Sigma,M)$. Then since any (tagged) arc on $(\Sigma,M)$ can be completed into a (tagged) triangulation, there is a bijection:
$$\xymatrix{\ind IY(\Sigma,M)\ar@{<..>}[rr]^-{1:1} && \{ \textrm{tagged arcs on }(\Sigma, M) \}}.$$

Now fix an ideal triangulation $\tau$ of $(\Sigma, M)$ and denote by $T_\tau$ the associated cluster-tilting of $\Cc_{(\Sigma,M)}$. Using results of \cite{Lab2} together with the diagram \ref{equation diagram} we obtain the following commutative diagram:
$$\xymatrix{\ind IY(\Sigma,M)\setminus \ind(T_\tau)\ar@{<..}[rr]\ar@{^(->}[d] && \{ \textrm{arcs  on }(\Sigma, M) \textrm{ not in } \tau \}\ar[ddl]^{X_\tau}
\\ \Cc_{(\Sigma,M)}\ar[dr]_{F_{T_\tau[-1]}} && 
\\ &\mod \Jac (Q_\tau,W_\tau) &
},$$
where  $F_{T_\tau[-1]}=\Hom_\Cc(T_\tau[-1],-)\simeq \Hom_\Cc(T_\tau, -[1])$.

\begin{remark}\label{Remark Bruestle}
In the case of an unpunctured surface $(\Sigma,M)$, Br\"ustle and Zhang study the category $\Cc_{(\Sigma,M)}$ very precisely in  \cite{BZ10} (see also \cite{ABCP}), and show an analogue of Theorem \ref{HU theorem}, that is any cluster-tilting object of $\Cc_{(\Sigma,M)}$ is in $IY(\Sigma,M)$.
\end{remark}

\begin{example}\label{example surface}
Let $\Sigma$ be a surface of genus $0$ with  boundary having two connected components, and let $M$ be a set of marked points on $\Sigma$ consisting of three marked points on the boundary, and one puncture. Let $\tau$ be the following ideal triangulation
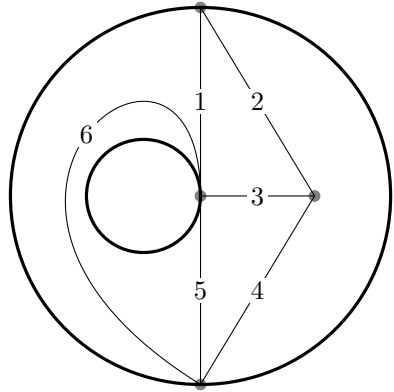
\[\scalebox{1}{
\begin{tikzpicture}[>=stealth,scale=1]
\filldraw [gray] (0,0) circle (2pt);
\filldraw [gray] (1.5,0) circle (2pt);
\filldraw [gray] (0,2.5) circle (2pt);
\filldraw [gray] (0,-2.5) circle (2pt);

\draw (0,0) [very thick] circle (2.5cm);
\draw (-0.75,0) [very thick] circle (0.75cm);
\draw (0,0)-- node [fill=white,inner sep=.5mm]{$1$}(0,2.5);
\draw (0,0)--node [fill=white,inner sep=.5mm]{$5$}(0,-2.5);
\draw (0,0)--node [fill=white,inner sep=.5mm]{$3$}(1.5,0);
\draw (1.5,0)--node [fill=white,inner sep=.5mm]{$2$}(0,2.5);
\draw (1.5,0)--node [fill=white,inner sep=.5mm]{$4$}(0,-2.5);
\draw (0,0)..controls (0,3) and (-4,0)..node [fill=white,inner sep=.5mm]{$6$}(0,-2.5);

\end{tikzpicture}}\]

The quiver with potential $(Q_\tau,W_\tau)$ is the following
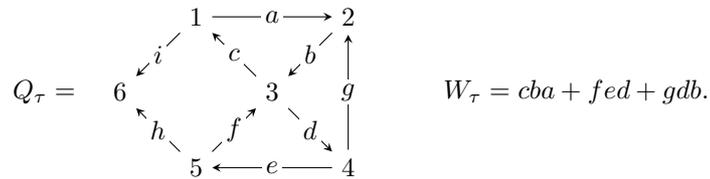
\[\scalebox{1}{
\begin{tikzpicture}[>=stealth,scale=1]
\node (P0) at (-1,0) {$Q_\tau=$};
\node (P1) at (1,1){$1$};
\node (P2) at (3,1){$2$};
\node (P3) at (2,0){$3$};
\node (P4) at (3,-1){$4$};
\node (P5) at (1,-1){$5$};
\node (P6) at (0,0){$6$};
\draw [->] (P1)--node[fill=white,inner sep=.5mm]{$a$} (P2);\draw [->] (P2)--node[fill=white,inner sep=.5mm]{$b$}(P3);\draw [->] (P3)--node[fill=white,inner sep=.5mm]{$c$}(P1);\draw [->] (P3)--node[fill=white,inner sep=.5mm]{$d$}(P4);\draw [->] (P4)--node[fill=white,inner sep=.5mm]{$e$}(P5);\draw [->] (P5)--node[fill=white,inner sep=.5mm]{$f$}(P3);\draw [->] (P4)--node[fill=white,inner sep=.5mm]{$g$}(P2);\draw [->] (P5)--node[fill=white,inner sep=.5mm]{$h$}(P6);\draw [->] (P1)--node[fill=white,inner sep=.5mm]{$i$}(P6);

\node (P7) at (6,0){$W_\tau=cba+fed+gdb$.};
 
\end{tikzpicture}}\]

Let $j$ be an arc of $(\Sigma,M)$ which is not in $\tau$.  In the most simple cases, the composition series of the module $X_\tau(j)$ correspond to the arcs of $\tau$ crossed transversally by $j$. For instance, if $j$ is the following arc
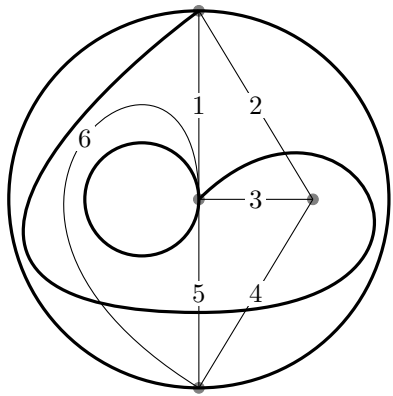
\[\scalebox{1}{
\begin{tikzpicture}[>=stealth,scale=1]
\filldraw [gray] (0,0) circle (2pt);
\filldraw [gray] (1.5,0) circle (2pt);
\filldraw [gray] (0,2.5) circle (2pt);
\filldraw [gray] (0,-2.5) circle (2pt);

\draw (0,0) [very thick] circle (2.5cm);
\draw (-0.75,0) [very thick] circle (0.75cm);
\draw (0,0)-- node [fill=white,inner sep=.5mm]{$1$}(0,2.5);
\draw (0,0)--node [fill=white,inner sep=.5mm]{$5$}(0,-2.5);
\draw (0,0)--node [fill=white,inner sep=.5mm]{$3$}(1.5,0);
\draw (1.5,0)--node [fill=white,inner sep=.5mm]{$2$}(0,2.5);
\draw (1.5,0)--node [fill=white,inner sep=.5mm]{$4$}(0,-2.5);
\draw (0,0)..controls (0,3) and (-4,0)..node [fill=white,inner sep=.5mm]{$6$}(0,-2.5);

\draw [very thick](0,2.5).. controls (-2,1) and (-4,-1.5)..(0,-1.5);
\draw [very thick] (0,-1.5)..controls (4,-1.5) and (2,2)..(0,0);
\end{tikzpicture}}\]
the module $X_\tau(j)$ has the following composition series $X_\tau(j)={\bsm &&&6\\2&&5&\\&4&&&\esm}$. 

Let $j$ be an arc of $\tau$ and $j'$ be the flip of this arc with respect to $\tau$. Then the arc $j'$ crosses the triangulation $\tau$ uniquely  through the arc $j$, thus the corresponding $\Jac(Q_\tau,W_\tau)$-module $X_\tau(j')$ is the simple module associated with the vertex of $Q_\tau$ corresponding to the arc $j$ of $\tau$.  

\end{example}

\begin{remark}
Labardini also associates in \cite{Lab1} QPs to ideal triangulations of surfaces without boundary.  For instance if $(\Sigma,M)$ is a once punctured torus, then all triangulations give the same QP:
$$\xymatrix{&&2\ar@<.5ex>[dl]^a\ar@<-.5ex>[dl]_{a'}&&\\ Q=&1\ar@<.5ex>[rr]^c\ar@<-.5ex>[rr]_{c'}&&3\ar@<.5ex>[ul]^b\ar@<-.5ex>[ul]_{b'}&W=abc+a'b'c'+ab'ca'bc'.}$$
This QP is not rigid, but is non-degenerate and Jacobi-finite (cf. \cite[Section 8]{Lab2}).  Using Corollary~\ref{Jacobi finite cluster categories}, it is possible to associate a generalized cluster category with cluster-tilting objects to this surface. However the non-degeneracy and Jacobi-finiteness are not known to be true or false for general surfaces without boundary. 
\end{remark}

\subsection{Derived preprojective algebras}
Derived preprojective algebras give another application of Theorem \ref{theorem general}.

\begin{definition}[Keller \cite{Kel09}]
Let $\Lambda$ be a finite dimensional algebra of global dimension at most $2$. 
Denote by $\Theta_2$ a cofibrant resolution of 
$\RHom_{\Lambda^{\rm e}}(\Lambda,\Lambda^{\rm e})[2]\in\Dd(\Lambda^{\rm e})$.
Then the  \emph{derived $3$-preprojective algebra} is
defined as the tensor DG algebra
\[ \BPi_{3}(\Lambda):= \Talg_\Lambda \Theta_2=\Lambda\oplus \Theta_2\oplus (\Theta_2\ten_\Lambda \Theta_2)\oplus \ldots\]
and the algebra $\Pi_{3}(\Lambda):=H^0(\BPi_{3}(\Lambda)) $ is called the \emph{$3$-preprojective algebra}.
\end{definition}

Since the algebra $\Lambda$ is a finite dimensional algebra of finite global dimension, there is an isomorphism of $\Lambda$-$\Lambda$-bimodules $\RHom_{\Lambda^{\rm e}}(\Lambda, \Lambda^{\rm e})\simeq \RHom_{\Lambda}(D\Lambda,\Lambda)$. Moreover the functor 
$$-\lten_{\Lambda}\RHom_{\Lambda}(D\Lambda, \Lambda):\Db(\Lambda)\longrightarrow \Db(\Lambda)$$ is a quasi inverse for the Serre functor of $\Db(\Lambda)$ (see \cite{Kel94} or \cite[Chapter 1]{Ami08}). Hence we have an isomorphism $\Pi_3(\Lambda)\simeq \bigoplus_{p\geq 0}\Hom_{\Db\Lambda}(\Lambda, \SSS_2^{-p}\Lambda).$ Since $\Lambda$ is of global dimension at most 2, then 
$\Hom_{\Db\Lambda}(\Lambda, \SSS_2^p\Lambda)$ vanishes for $p\geq 1$. Therefore we have isomorphisms of algebras:
\begin{equation}\label{equation isomorphism}\Pi_3(\Lambda)\simeq \bigoplus_{p\in \ZZ}\Hom_{\Db\Lambda}(\Lambda, \SSS_2^{-p}\Lambda)\simeq \Talg_\Lambda\Ext_\Lambda^2(D\Lambda,\Lambda) \quad \textrm{(cf. \cite[Prop. 4.7]{Ami09}}).\end{equation}

\begin{definition}[Iyama]
An algebra of global dimension at most $2$ is called \emph{$\tau_2$-finite} if $\Pi_3(\Lambda)$ is finite dimensional. This is equivalent to the fact that the endofunctor $\tau_2:=H^0(\SSS_2)$ of $\mod \Lambda$ is nilpotent.
\end{definition}

Note that if $\Lambda$ is an hereditary algebra, then $\Pi_3(\Lambda)\simeq \Lambda$ since $\Ext^2_\Lambda(D\Lambda,\Lambda)$ vanishes. Hence for any acyclic quiver $Q$, the algebra $kQ$ is $\tau_2$-finite.

\begin{theorem}[Keller, \cite{Kel09}]\label{theorem derived preproj}
Let $\Lambda$ be an algebra of global dimension at most 2. Then the derived preprojective algebra $\BPi_3(\Lambda)$ is homologically smooth and bimodule $3$-CY.
\end{theorem} 
Applying then Theorem \ref{theorem general} we get the following construction.
\begin{corollary}\label{corollary existence cluster cat}
Let $\Lambda$ be as in Theorem \ref{theorem derived preproj}. If $\Lambda$ is moreover $\tau_2$-finite, then the category $$\Cc_2(\Lambda):=\per \BPi_3(\Lambda)/\Db(\BPi_3(\Lambda))$$ is $\Hom$-finite and $2$-Calabi-Yau. 

Moreover the object $\BPi_{3}(\Lambda)$ is cluster-tilting in $\Cc_2(\Lambda)$ with endomorphism algebra $\Pi_3(\Lambda)$.
\end{corollary}
This category is called \emph{the cluster category associated with an algebra of global dimension at most 2}. 
The next result gives the link between the category $\Cc_2(\Lambda)$ and the orbit category $\Db(\Lambda)/\SSS_2$.

\begin{theorem}\label{theorem nice generalization}
Let $\Lambda$ be a finite dimensional algebra of global dimension at most~2. 
\begin{itemize}
\item[(1)] For any $X$ in $\Db(\Lambda)$, there is an isomorphism $\SSS_2(X)\lten_{\Lambda}\BPi_3(\Lambda)\simeq X\lten_{\Lambda}\BPi_3(\Lambda)$ in $\Cc_2(\Lambda)$  functorial in $X$. Thus there is a commutative diagram
$$\xymatrix{\Db(\Lambda)\ar[rr]^{-\lten_{\Lambda}\BPi_3(\Lambda)}\ar[d]^\pi &&  \per{\BPi_3(\Lambda)}\ar[d]\\ 
\Db(\Lambda)/\SSS_2\ar@{..>}[rr] && \Cc_2(\Lambda)\,.\!\!\!}$$
Moreover, the factorization $\Db(\Lambda)/\SSS_2\rightarrow \Cc_2(\Lambda)$ is fully faithful.

\item[(2)] The smallest triangulated subcategory of $\Cc_2(\Lambda)$ containing the orbit category $\Db(\Lambda)/\SSS_2$ is $\Cc_{2}(\Lambda)$.
\end{itemize}
\end{theorem}

Let $J$ be the ideal $J=\Theta_2\oplus (\Theta_2\ten_\Lambda \Theta_2)\oplus \ldots$ of $\BPi_3(\Lambda)$. Then using the isomorphism \ref{equation isomorphism}, we obtain an isomorphism $\SSS_2(X)\lten_{\Lambda} J\simeq X\lten_{\Lambda}\BPi_3(\Lambda)$. The morphism in (1) of the above theorem is induced by the inclusion $J\to \BPi_3(\Lambda)$ in $\per (\BPi_3(\Lambda)^{\rm e})$.  The cone of  $J\to \BPi_3(\Lambda)$ is in $\Db(\BPi_3(\Lambda)^{\rm e})$ since $\Lambda$ is finite dimensional. Therefore the morphism   $$\SSS_2(X)\lten_{\Lambda}\BPi_3(\Lambda)\simeq X\lten_{\Lambda}J\rightarrow X\lten_{\Lambda}\BPi_3(\Lambda) \quad \textrm{in } \per \BPi_3(\Lambda)$$ has its cone in $\Db(\BPi_3(\Lambda))$, so is an isomorphism in $\Cc_2(\Lambda).$

Point $(2)$ of this theorem says that $\Cc_2(\Lambda)$ can be understood as the triangulated hull of the orbit category $\Db(\Lambda)/\SSS_2$.  Keller introduced this notion with DG categories. A DG category is a category where morphisms have a structure of $k$-complexes. We refer to \cite{Kel94}, \cite{Kel05} for precise definitions and constructions (see also \cite[Appendix]{IO09}, and \cite[Section 6]{AO10a}).

The philosophy is the following: any algebraic triangulated category is triangle equivalent to $H^0(\Xx)$, where $\Xx$ is a DG category. But for a given DG category $\Xx$, the category $H^0(\Xx)$ is not always triangulated. However the category $H^0(\Xx)$  can be viewed (via the Yoneda functor) as a full subcategory of $H^0({\sf DGMod }\Xx)$ which is a triangulated category.  Here ${\sf DGMod} \Xx $ is the DG category of the DG $\Xx$-modules.  The triangulated hull of $H^0(\Xx)$ is defined to be the smallest triangulated subcategory of  $H^0({\sf DGMod }\Xx)$ containing $H^0(\Xx)$. In our situation, the category $\Db(\Lambda)$ has a canonical enhancement in a DG category $\Xx$ (that is $\Db(\Lambda)$ is equivalent to some category $H^0(\Xx)$ for a canonical $\Xx$ ). The functor $\SSS_2$ can be canonically lifted to a DG functor $S$ of $\Xx$. One can define the orbit category $\Xx/S$; it is a DG category. The orbit category $\Db(\Lambda)/\SSS_2$ is equivalent to the category $H^0(\Xx/S)$. Its triangulated hull is defined to be the triangulated hull of $H^0(\Xx/S)$ following the previous construction. Since the enhancement $\Xx$ and the lift $S$ are canonical, one can speak of \emph{the} triangulated hull. But in general, triangulated hulls depend on the choice of the enhancement.

From Theorem \ref{theorem nice generalization} one deduces the following fact.
\begin{corollary}
Let $\Lambda=kQ$ for an acyclic quiver $Q$. Then $kQ$ is $\tau_2$-finite and there is a triangle equivalence $\Cc_2(kQ)\simeq \Cc_Q$.
\end{corollary}
This corollary associated with Theorem \ref{theorem nice generalization} shows that $\Cc_2(\Lambda)$ is the most natural generalization of the cluster category when replacing the algebra $kQ$ of global dimension 1 by the algebra $\Lambda$ of global dimension $2$.

\begin{remark}
If $T$ is a tilting $kQ$-module, where $Q$ is an acyclic quiver, then the algebra $B=\End_{kQ}(T)$ has global dimension at most 2 and is $\tau_2$-finite. The triangle functor $-\lten_BT:\Db(B)\to \Db(kQ)$ is an equivalence (see \cite{Hap87}) sending $B\in \Db(B)$ on $T\in \Db(kQ)$. This equivalence induces a triangle equivalence between the cluster categories $\Cc_2(B)$ and $\Cc_Q$ sending $\pi(B)$ on $\pi(T)$. Moreover the tensor product $\Ext^2(DB,B)\ten_B\Ext^2(DB,B)$ vanishes, hence the algebra $\Pi_3(B)$ is isomorphic to the trivial extension $B\oplus \Ext^2(DB,B)$.  The isomorphism $\Pi_3(B)\simeq \End_{\Cc_Q}(T)$  recovers then a result of \cite{ABS08}.

Note that in general, a $\tau_2$-finite algebra $B$ of global dimension $2$ can satisfy $\Ext^2(DB,B)\ten_B\Ext^2(DB,B)\neq 0$, even if $\Cc_2(B)$ is an acyclic cluster category, as shown in the following example.
Let $B$ be the Auslander algebra of $\mod kQ$, where $Q$ is the linear orientation of $A_3$

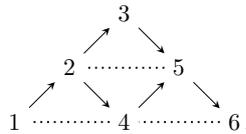
\[\scalebox{0.9}{
\begin{tikzpicture}[>=stealth,scale=0.8]
\node (A1) at (4,0) {$1$};
\node (A2) at (5,1) {$2$};
\node (A3) at (6,2) {$3$};
\node (A4) at (6,0) {$4$};
\node (A5) at (7,1) {$5$};
\node (A6) at (8,0) {$6$\,.\!\!\!};

\draw [->] (A1)--(A2);\draw [->] (A2)--(A3);\draw [->] (A2)--(A4);\draw [->] (A3)--(A5);\draw [->] (A5)--(A6);\draw [->] (A4)--(A5);
\draw[dotted, thick] (A4)--(A1);\draw[dotted, thick] (A6)--(A4);\draw[dotted, thick] (A5)--(A2);

\end{tikzpicture}}\]
Then an easy computation gives $\Ext^2_B(e_4DB,e_1B)\neq 0$ and $\Ext^2_B(e_6DB,e_4B)\neq 0$, hence we have
$$e_1\Ext^2_B(DB,B)\ten_B\Ext^2_B(DB,B)e_6=\Ext^2_B(DB,e_1B)\ten_B\Ext^2_B(e_6DB,B)\neq 0.$$

However, the category $\Cc_2(B)$ is equivalent to the acyclic cluster category $\Cc_{D_6}$, since the algebra $B$ is an iterated tilted algebra of type $D_6$ and of global dimension~2.

\end{remark} 

From Keller's results on triangulated hulls \cite{Kel05}, one deduces the following property for the category $\Cc_2(\Lambda)$ (see Appendix of \cite{IO09} for details).

\begin{theorem}\label{univ prop}
Let $\Lambda$ be a finite dimensional algebra of global dimension at most $2$.
Let $\Ee$ be a Frobenius category. Let $M$ be an object in $\Db(\Lambda^{\rm op}\ten\Ee)$. Assume that there is a morphism $f:\RHom_{\Lambda^{\rm e}}(\Lambda,\Lambda^{\rm e})\lten_{\Lambda}M [2]\rightarrow M$ in $\Db(\Lambda^{\rm op}\ten\Ee)$ such that the cone of $f$, viewed as an object in $\Db(\Ee)$ is perfect. Then there is triangle functor $F:\Cc_2(\Lambda)\rightarrow \underline{\Ee}$ making the following diagram commutative:
$$\xymatrix{\Db(\Lambda)\ar[rr]^{-\lten_{\Lambda}M} \ar[d]^{\pi}&& \Db(\Ee)\ar[d]\\ \Cc_2(\Lambda)\ar[rr]^F && \Db(\Ee)/\per(\Ee)\simeq \underline{\Ee} .}$$
\end{theorem}

\begin{remark}
A very powerful application of cluster categories associated with algebras of global dimension at most $2$ is made in \cite{Kel10}, and allows to prove the periodicity conjecture for pairs of Dynkin diagrams (see also last sections of \cite{Kellersurvey} for an overview of this result).
\end{remark}

\subsection{Link between the two constructions}\label{section link}

In the last two subsections we have generalized the notion of cluster category, first using the Ginzburg DG algebra, and then derived $3$-preprojective algebra. So a natural question is the following:

\medskip
\emph{
What is the link between the categories $\Cc_{(Q,W)}$, where $(Q,W)$ is a Jacobi-finite QP, and the categories $\Cc_2(\Lambda)$, where $\Lambda$ is a $\tau_2$-finite algebra of global dimension at most 2?}
\medskip

Let $\Lambda=kQ/I$ (where $I$ is an admissible ideal of $kQ$) be a finite dimensional algebra of global dimension at most $2$. Let $\Rr$ be a minimal set of relations, i.e. the lift to $I$ of a basis of $I/(IJ+JI)$ (where $J$ is the ideal of $kQ$ generated by the arrows) compatible with the decomposition $I=\bigoplus_{i,j}e_jIe_i$. Let $\bar{Q}$ be the quiver obtained from $Q$ by adding an arrow $a_r:t(r)\rightarrow s(r)$ for each relation $r:s(r)\rightarrow t(r)$. Let $\bar{W}_{\Lambda}$ be the potential $\bar{W}_{\Lambda}=\sum_{r\in \Rr}ra_r$.

\begin{example}\label{example 3cycle}
Let $\Lambda$ be the algebra presented by the quiver $$\xymatrix@-0.5cm{& 2\ar[dl]_a & \\ 1\ar[rr]^c&&3\ar[ul]_b} \quad \textrm{and  the relation }ab=0.$$
Then the associated quiver with potential $(\bar{Q}_\Lambda,\bar{W})$ is
 $$\xymatrix@-0.3cm{& 2\ar[dl]_a & \\ 1\ar@<.5ex>[rr]^c\ar@<-.5ex>[rr]_{r_{ab}}&&3\ar[ul]_b} \quad \textrm{and   }\bar{W}=r_{ab}ab.$$
\end{example}

The following result is due to Keller \cite[Theorem 6.12]{Kel09}.

\begin{theorem}[Keller]\label{theorem global dimension 2 = ginzburg} Let $\Lambda$ be a finite dimensional algebra of global dimension at most 2. 
There exists a isomorphism of DG algebras $f:\Gamma(\bar{Q},\bar{W}_{\Lambda})\rightarrow \BPi_{3}(\Lambda)$ such that $H^n(f)$ is an isomorphism for each $n\in \mathbb{Z}$.

This quasi-equivalence induces triangle equivalences $\per \Gamma(\bar{Q},\bar{W}_{\Lambda})\simeq\per \BPi_3(\Lambda)$ sending the object $\Gamma(\bar{Q},\bar{W}_{\Lambda})$ to the object $\BPi_{3}(\Lambda)$ and $\Db \Gamma(\bar{Q},\bar{W}_{\Lambda})\simeq\Db \BPi_3(\Lambda)$.
\end{theorem}
\noindent
This result implies in particular an isomorphism of algebras \begin{equation}\label{equation jacobian}\Pi_{3}(\Lambda)\simeq \Talg_{\Lambda}\Ext^2(D\Lambda,\Lambda)\simeq \Jac(\bar{Q}_{\Lambda},\bar{W}).\end{equation}
\noindent
 It also gives an answer to the above question.
\begin{corollary}\label{corollary cluster cat}
Let $\Lambda$ and $(\bar{Q}_\Lambda, \bar{W}_\Lambda)$ be as in Theorem~\ref{theorem global dimension 2 = ginzburg}. Assume moreover that $\Lambda$  is $\tau_2$-finite. Then there exists a triangle equivalence $\Cc_2(\Lambda) \simeq \Cc_{(\bar{Q}_\Lambda, \bar{W}_\Lambda)}$ sending the cluster-tilting object $\BPi_3(\Lambda)\in \Cc_2(\Lambda)$ on the cluster-tilting object $\BPi_{3}(\Lambda)\in \Cc_{(\bar{Q}_\Lambda, \bar{W}_\Lambda)}$.
\end{corollary}

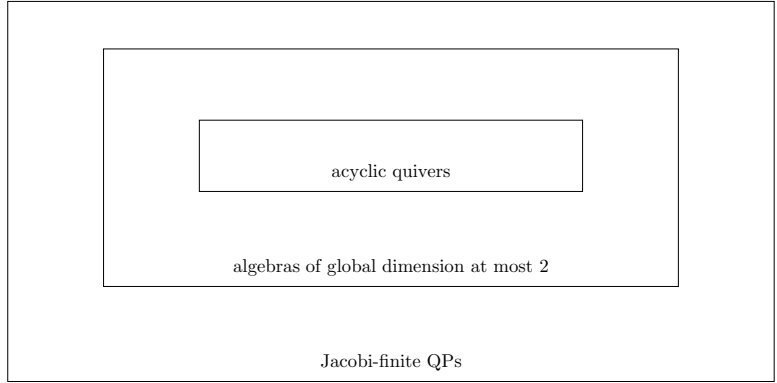
\[\scalebox{0.7}{
\begin{tikzpicture}[>=stealth,scale=1.8]
\node (A) at (4,0.2){Jacobi-finite QPs};
\node (B) at (4,1.2) {algebras of global dimension at most 2};
\node (C) at (4,2.2){acyclic quivers};
\draw (0,0) rectangle (8,4);
\draw (1,1) rectangle (7,3.5);
\draw (2,2) rectangle (6,2.75);
\end{tikzpicture}}\]

\begin{remark}
If $Q$ is an acyclic quiver, and $\Lambda=kQ$. Then  we have $\bar{Q}=Q$ and $\bar{W}=0$, hence we recover $\Pi_3(kQ)\simeq kQ\simeq \Jac(Q,0)$.

If $T$ is a tilting $kQ$-module (where $Q$ is an acyclic quiver), and if $\Lambda=\End_{kQ}(T)$, then the description of $\bar{Q}_\Lambda$ was already given in \cite{ABS08}. 
\end{remark}

\begin{remark}
One can show that the category $\Cc_{(Q,W)}$, where 
$$ Q=\xymatrix@-.5cm{&2\ar[dr]^b&\\ 1\ar[ur]^a && 3\ar[ll]_c} \quad \textrm{and} \quad W=cbacba ,$$ is not triangle equivalent to a generalized cluster category $\Cc_{2}(\Lambda)$ associated with an algebra of global dimension at most~2. Therefore the inclusion given in Corollary~\ref{corollary cluster cat} is strict.
\end{remark}

\begin{example}
Let $(\Sigma,M)$ be the surface with marked points as in Example \ref{example surface}.
Let $\Lambda$ be the algebra presented by the quiver 
\[\scalebox{1}{
\begin{tikzpicture}[>=stealth,scale=1]
\node (P0) at (-1,0) {$Q=$};
\node (P1) at (1,1){$1$};
\node (P2) at (3,1){$2$};
\node (P3) at (2,0){$3$};
\node (P4) at (3,-1){$4$};
\node (P5) at (1,-1){$5$};
\node (P6) at (0,0){$6$};
\draw [loosely dotted, thick] (P1)..controls (2,0.25).. (P2);\draw [->] (P2)--node[fill=white,inner sep=.5mm]{$b$}(P3);\draw [->] (P3)--node[fill=white,inner sep=.5mm]{$c$}(P1);\draw [->] (P3)--node[fill=white,inner sep=.5mm]{$d$}(P4);\draw [loosely dotted, thick] (P4)..controls (2,-0.25)..(P5);\draw [->] (P5)--node[fill=white,inner sep=.5mm]{$f$}(P3);\draw [loosely dotted, thick] (P4)..controls (2.5,0)..(P2);\draw [->] (P5)--node[fill=white,inner sep=.5mm]{$h$}(P6);\draw [->] (P1)--node[fill=white,inner sep=.5mm]{$i$}(P6);

\node (P7) at (6,0){with relations $cb=df=db=0$.};
 
\end{tikzpicture}}\]
This is a $\tau_2$-finite algebra of global dimension $2$. By Corollary \ref{corollary cluster cat}, there is an equivalence of triangulated categories  $\Cc_2(\Lambda)\simeq\Cc_{(\Sigma,M)}$ sending the cluster-tilting object $\pi(\Lambda)$ on the cluster-tilting object $T_\tau$ associated with the triangulation~$\tau$. 

We write $\pi(\Lambda)=\pi(e_2 \Lambda)\oplus T_0\in \Cc_2(\Lambda)$. One easily checks that there is an exchange triangle in $\Cc_2(\Lambda)$
$$\xymatrix{ \pi(e_2 \Lambda)\ar[r] & \pi (e_3 \Lambda)\ar[r] & \pi ({\bsm 3\\5 \esm})\ar[r] & \pi(e_2 \Lambda)[1]},$$
thus the object $\pi ({\bsm 3\\5 \esm})\oplus T_0$ is a cluster-tilting object in $\Cc_2(\Lambda)$.

Denote by $2'$ the flip of the arc $2$ with respect to the triangulation $\tau$. Then one checks that the image of the $\Lambda$-module ${\bsm 3\\5 \esm}$ through the functors
$$ \xymatrix@C=.2cm{\Db(\Lambda)\ar[d]^\pi & \\ \Cc_2(\Lambda)\ar@{-}[r]^\sim & \Cc_{(\Sigma,M)}\ar[d]^{F_{T_\tau[-1]}} \\ & \mod \Jac (Q_\tau,W_\tau)}$$ is isomorphic to the module $X_\tau(2')$ constructed by Labardini in \cite{Lab2}.  By the same process,  for  any arc $j$ of the surface, by making iterated IY-mutations and corresponding flips, one can construct an object $M$ in $\Db(\Lambda)$ such that $\pi(M)$ corresponds to $j$. 

For instance, the summand corresponding to the vertex $6$ of the object $\mu_6^L \circ\mu_5^L\circ\mu_4^L\circ\mu_2^L(\pi(\Lambda))$ of $\Db(\Lambda)$ (see Section 5 for a definition of the graded mutation $\mu^L)$ corresponds to the arc $j$ of Example \ref{example surface}.

\end{example}

\subsection{Higher generalized cluster categories}

In general, it is of interest to study $n$-Calabi-Yau categories $\Cc$ with $n$-cluster-tilting objects, i.e. objects $T$ such that $$\begin{array}{rcl}\add(T)& = &\{ X\in \Cc \mid \Ext^i_{\Cc}(X,T)=0, \ 1\leq i\leq n-1\}\\ &=&\{X\in \Cc\mid \Ext^i_{\Cc}(T,X)=0, \ 1\leq i \leq n-1\}.\end{array}$$

Theorem \ref{theorem general} generalizes when replacing the hypothesis $(d)$ by $\BPi$ being bimodule $n$-Calabi-Yau. The formal proof is given in \cite{Guo10}. Associated to a finite dimensional algebra $\Lambda$ of global dimension at most $n$, Keller also introduced in \cite{Kel09}  the \emph{derived $(n+1)$-preprojective algebra} $\BPi_{n+1}(\Lambda)$ as the DG algebra $\Talg_{\Lambda}\Theta_n$, where $\Theta_n$ is a cofibrant resolution of $\RHom_{\Lambda^{\rm e}}(\Lambda, \Lambda^{\rm e})[n]$. Then Corollary \ref{corollary existence cluster cat} generalizes when replacing $2$ by $n$.

In recent papers \cite{G06}, \cite{BSW10}, \cite{Kel09}, \cite{VdB}, generalizations of Ginzburg DG algebras are introduced and studied. 

\bigskip\noindent
\textbf{Special case $n=1$.}

The case $n=1$ is notable. A finite dimensional algebra of global dimension at most one is Morita equivalent to some path algebra of an acyclic quiver $Q$. The category $\Cc_1(kQ)$ is $\Hom$-finite if and only if the algebra $kQ$ is $\tau_1$-finite, that is when $Q$ is Dynkin, since $\tau_1$ is isomorphic to the usual Auslander-Reiten translation. In this case, the orbit category $\Db(kQ)/\SSS_1$ is equivalent to the category $\proj \Pi_2(kQ)$ and is naturally triangulated (\cite{Kel05} or \cite{Ami07}). They have finitely many indecomposables. A $1$-cluster-tilting object corresponds to a basic generator-cogenerator, therefore $\Pi_2(kQ)$ is the only $1$-cluster-tilting object in $\Cc_{1}(kQ)$. From results in \cite{BES} and \cite{Ami07}, we deduce the following.

\begin{theorem}
Let $Q$ be a Dynkin quiver.
\begin{itemize}
\item[(a)] \cite{Ami07} The deformed preprojective algebras defined in \cite{BES} are $1$-CY-tilted algebras.
\item[(b)] \cite{BES} If $char(k)=2$, then the deformed $2$-preprojective algebras are not isomorphic to $2$-preprojective algebras.
\end{itemize}
\end{theorem}
 Hence we may expect that there also exist some \emph{higher deformed preprojective algebras} that are $2$-CY-tilted algebras, which would yield a negative answer to Question~\ref{problem 2}~(1).

\section{Stable categories as generalized cluster categories}
In this section, we see how generalized cluster categories generalize other constructions of $2$-Calabi-Yau categories via Frobenius categories.

\subsection{Preprojective algebras of Dynkin type}
Once the notion of generalized cluster category is established, it is natural to see if it answers  Question \ref{problem 1}.

\begin{theorem}[\cite{Ami09}]\label{preproj Dynkin=cluster}
Let $Q$ be a Dynkin quiver.  Let $\underline{\Lambda}$ be the stable Auslander algebra of the category $\mod kQ$. Then there is a triangle equivalence
$$\underline{\mod}\Pi_2(kQ)\simeq \Cc_2(\underline{\Lambda}),$$ sending the canonical cluster-tilting object $\Lambda$ to the standard cluster-tilting object $M$ associated with the orientation of $Q$.
\end{theorem}

The proof consists in constructing a triangle functor $\Db(\underline{\Lambda})\to \Db(\Pi_2(kQ))$ and in using Theorem~\ref{univ prop} to deduce a triangle functor $\Cc_2(\underline{\Lambda})\to \underline{\mod}\Pi_2(kQ)$. We explain the construction of the functor in the example below.

\begin{example} Let $Q$ be the linear orientation of $A_3$. Let $X={\bsm1\esm}\oplus{\bsm2\\1\esm}\oplus{\bsm3\\2\\1\esm}\oplus{\bsm2\esm}\oplus{\bsm2\\1\esm}\oplus{\bsm3\\2\esm}\oplus{\bsm3\esm}$ be the $kQ$-generator-cogenerator. Define the algebras $\Lambda=\End_{kQ}(X)$, $\underline{\Lambda}:=\underline{\End}_{kQ}(X)$, and $\Gamma=\End_{\Pi}(M)$, where $M$ is the standard cluster-tilting object of the module category $\mod \Pi_2(Q)$. Then one easily checks that $\Lambda$ is a subalgebra of $\Gamma$ and that we have natural algebra morphisms:

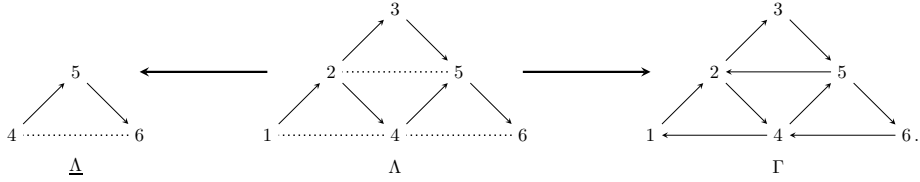
\[\scalebox{0.7}{
\begin{tikzpicture}[>=stealth,scale=1.2]
\node (L1) at (0,0) {$4$};
\node (L2) at (1,1) {$5$};
\node (L3) at (2,0) {$6$};
\node (Lbar) at (1,-0.5) {$\underline{\Lambda}$};
\node (A1) at (4,0) {$1$};
\node (A2) at (5,1) {$2$};
\node (A3) at (6,2) {$3$};
\node (A4) at (6,0) {$4$};
\node (A5) at (7,1) {$5$};
\node (A6) at (8,0) {$6$};
\node (L) at (6,-0.5) {$\Lambda$}; 
\node (B1) at (10,0) {$1$};
\node (B2) at (11,1) {$2$};
\node (B3) at (12,2) {$3$};
\node (B4) at (12,0) {$4$};
\node (B5) at (13,1) {$5$};
\node (B6) at (14,0) {$6$\,.\!\!\!};
\node (B) at (12,-0.5) {$\Gamma$};
\draw [->] (L1)--(L2);\draw [->] (L2)--(L3); \draw[dotted, thick] (L3)--(L1); 
\draw [->] (A1)--(A2);\draw [->] (A2)--(A3);\draw [->] (A2)--(A4);\draw [->] (A3)--(A5);\draw [->] (A5)--(A6);\draw [->] (A4)--(A5);
\draw[dotted, thick] (A4)--(A1);\draw[dotted, thick] (A6)--(A4);\draw[dotted, thick] (A5)--(A2);
\draw [->] (B1)--(B2);\draw [->] (B2)--(B3);\draw [->] (B2)--(B4);\draw [->] (B3)--(B5);\draw [->] (B5)--(B6);\draw [->] (B4)--(B5);
\draw[->] (B4)--(B1); \draw[->] (B6)--(B4); \draw[->] (B5)--(B2);
\draw [->, very thick] (4,1)--(2,1); \draw[->, very thick](8,1)--(10,1);
\end{tikzpicture}}\]
Then we consider the following triangle functor $$\xymatrix{F:\Db(\underline{\Lambda})\ar[rr]^-{Res}& & \Db(\Lambda)\ar[rr]^{-\lten_{\Lambda}M} && \Db(\Pi_2(kQ))}$$ induced by these algebra morphisms and check that it satisfies the hypotheses of Theorem \ref{univ prop}.
\end{example}

\begin{remark}
This theorem has been very nicely generalized in \cite{IO09} for higher cluster categories as follows: For $n\geq 2$, let $A$ be a $(n-1)$-representation-finite algebra (in particular it is of global dimension $n-1$ and $\tau_{n-1}$-finite). Then there is a triangle equivalence $\underline{\mod}\Pi_{n}(A)\simeq \Cc_n(\Lambda)$, where $\Lambda$ is the stable Auslander algebra $\Lambda=\underline{\End}_A(\Pi_{n}(A))$. 
\end{remark}

\subsection{Stable categories associated to words}\label{subsection words}
In this subsection, we describe the alternative answer to Question \ref{problem 1} developed in \cite{BIRS09} and the link with the generalized cluster categories associated to algebras of global dimension at most $2$.

Let $\Delta$ be a graph and $\Pi=\Pi_2(\Delta)$ be the associated preprojective algebra (up to isomorphism, it does not depend on the orientation of $\Delta$). Define the Coxeter group $C_\Delta$ associated to $\Delta$ by the generators $s_i$, for $i\in \Delta_0$ and the relations $s_i^2=1$, $s_is_j=s_js_i$ if there are no arrows linking $i$ with $j$, and $s_is_js_i=s_js_is_j$ if there is exactly one edge linking $i$ and $j$.

To any element $w$ of $C_\Delta$, Buan, Iyama, Reiten and Scott have associated an algebra $\Pi_w$ which is defined to be the quotient $\Pi/I_w$, where $$I_w:=\Pi (1-e_{i_1})\Pi(1-e_{i_2})\Pi\ldots \Pi(1-e_{i_l})\Pi$$ and $s_{i_1}s_{i_2}\ldots s_{i_l}$ is a reduced expression of the element $w$. The ideal $I_w$,  hence the algebra $\Pi_w$, does not depend on the choice of this reduced expression.

\begin{theorem}[Buan-Iyama-Reiten-Scott \cite{BIRS09}]
For any element $w$ in $C_\Delta$ the category $\Sub \Pi_w$ formed by all $\Pi$-submodules of modules in $\add\Pi_w$ is Frobenius and $\Hom$-finite.
The category $\underline{\Sub}\Pi_w$ is $2$-Calabi-Yau and contains a cluster-tilting object $M_{\mathbf{w}}$ for each reduced expression $\mathbf{w}$ of $w$.

All cluster-tilting objects of the form $M_\mathbf{w}$ are in the same IY-mutation class. More precisely, if $\mathbf{w}$ and $\mathbf{w'}$ are two reduced expressions of the element $w$ linked by a braid relation, then $M_\mathbf{w}$ and $M_{\mathbf{w}'}$ are linked by a IY-mutation:
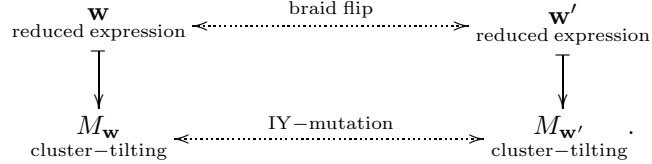
\[\xymatrix@C=1.7cm@R=0.7cm{\underset{\tiny{\rm reduced\ expression}}{\mathbf{w}}\ar@{<..>}[rr]^{\tiny{\rm braid \ flip}}\ar@{|->}[d] && \underset{\tiny{\rm reduced\ expression}}{\mathbf{w}'}\ar@{|->}[d]\\ \underset{\tiny{\rm cluster-tilting}}{M_{\mathbf{w}}}\ar@{<..>}[rr]^{\tiny{\rm IY-mutation}} && \underset{\tiny{\rm cluster-tilting}}{M_{\mathbf{w}'} }.}\]
\end{theorem}
Note that there are cluster-tilting objects in the IY-mutation class of $M_\mathbf{w}$ but which are not isomorphic to some $M_\mathbf{w'}$. 
 
These categories give an answer to Question \ref{problem 1}. 
Indeed, if $\Delta$ is a Dynkin graph and if $w$ is the element of maximal length in $C_\Delta$, then we have $\Sub\Pi_w\simeq \mod \Pi$. Moreover, there is a natural bijection between orientations of $\Delta$ and reduced expressions of $w$ (with commutativity relation). Using this bijection, the standard cluster-tilting objects of \cite{GLS06} correspond to the cluster-tilting objects $M_{\mathbf{w}}$.
For any acyclic quiver $Q$ (except $A_n$ with the linear orientation), there is a triangle equivalence $\underline{\Sub}\Pi_w\simeq \Cc_Q$ when $w=c^2$, where $c$ is the Coxeter element associated with the orientation of $Q$. See \cite{BIRSm}, \cite{AIRT} for further results on these categories. This construction and the constructions in Section 3 lead to the natural question.

\medskip
\emph{
What is the link between the categories $\Sub \Pi_w$ and $\Cc_2(\Lambda)$?}
\medskip

Let $\Delta$ be a graph. Then for any acyclic orientation $Q$ of $\Delta$, one can define a $\mathbb{Z}$-grading on $\Pi$, by assigning degree $0$ to the arrows of $Q$ and degree $1$ to the additional arrows. One easily checks that the preprojective relations are homogeneous of degree $1$. For any reduced expression $s_{i_1}\ldots s_{i_l}$ in $C_\Delta$, the ideal $\Pi (1-e_{i_1})\Pi\ldots \Pi (1-e_{i_l})\Pi$ is graded. The object $M_{\mathbf{w}}$ is defined as a sum of $e_i \Pi/I_w$, therefore it is a graded module.

\begin{theorem}[\cite{ART}]\label{theorem ART}
Let $\Delta$ be a graph and $\mathbf{w}$ be a reduced expression of an element $w$ in $C_\Delta$. Then there exists an orientation of $\Delta$ such that  there is a triangle equivalence $\underline{\Sub}\Pi_w\simeq \Cc_{2}(\underline{\Lambda})$ sending the canonical cluster-tilting object $\pi(\underline{\Lambda})$ on $M_{\mathbf{w}}$, where $\underline{\Lambda}$ is the stable endomorphism algebra $\underline{\End}_{\gr \Pi}(M_{\mathbf{w}})$ (in the category of graded $\Pi$-modules $\gr \Pi$).
\end{theorem}

As in Theorem \ref{preproj Dynkin=cluster}, the proof consists in constructing a triangle functor $\Db(\Lambda)\to \Db(\Sub\Pi_w)$ satisfying the hypotheses of Theorem~\ref{univ prop}.  
We may visualize the situation as follows

\[\scalebox{0.8}{
\begin{tikzpicture}[>=stealth,scale=1.9]
\node (A) at (4,0.2){cluster categories $\Cc_2(\Lambda)$};
\node (B) at (4,0.7) {categories $\underline{\Sub}\Pi_w$ associated with elements $w$ in Coxeter groups};
\node (C) at (2,1.2){acyclic quivers};
\node (D) at (5.5,2.8) {$\underline{\mod}\Pi_2(\Delta)$ with $\Delta$ Dynkin};
\draw (0,0) rectangle (8,4);
\draw (0.5,0.5) rectangle (7.5,3.5);
\draw (1,1) rectangle (4.1,2.1);
\draw (3.9,1.9) rectangle (7, 3);
\end{tikzpicture}}\]

\begin{example}
Let $\Delta$ be the graph $\widetilde{A}_2$, and $\mathbf{w}=s_2s_1s_2s_3s_2s_1s_3s_2\in C_\Delta$ be a reduced word.
The associated cluster-tilting object constructed in \cite{BIRS09} is the following:
$$M_{\mathbf{w}}={\bsm 2\esm}\oplus {\bsm 1\\2\esm}\oplus {\bsm 2\\1 \esm}\oplus {\bsm &&3&&\\&1&&2&\\2&&&&1\esm}\oplus {\bsm &&&2&\\&&3&&1\\&1&&2&\\2&&&&\esm}\oplus {\bsm &&&&1&&&\\&&&2&&3&&\\&&3&&1&&2&\\&1&&2&&&&1\\2&&&&&&&\esm}\oplus {\bsm &&&&&3&&\\&&&&1&&2&\\&&&2&&3&&1\\&&3&&1&&2&\\&1&&2&&&&\\2&&&&&&&\esm}\oplus {\bsm&&&&&&2&&&&\\&&&&&3&&1&&&\\&&&&1&&2&&3&\\&&&2&&3&&1&&2&\\&&3&&1&&2&&&&1\\&1&&2&&&&&&&\\2&&&&&&&&&&  \esm}$$ 
The last three summands (corresponding to the last $s_1$, $s_2$ and $s_3$ in the reduced word $\mathbf{w}$) are the projective-injectives of the category $\Sub\Pi_w$.
The endomorphism algebra $\End_\Pi(M_{\mathbf{w}})$ is the frozen Jacobian algebra (see \cite{BIRSm}) with  quiver $Q$
of the form 
\[\scalebox{0.7}{
\begin{tikzpicture}[>=stealth,scale=.8]
\node (B1) at (0,0){$1$};
\node (B2) at (2,0){$3$};
\node (B3) at (6,0){$5$};
\node (B4) at (10,0){$8$\,,\!\!\!};
\node (A1) at (1,2){$2$};
\node(A2) at (7,2){$6$};
\node (C1) at (4,-1){$4$};
\node (C2) at (8,-1){$7$};
\draw[->] (B2)--node[swap,yshift=-2mm] {$p_1$}(B1);\draw[->] (B4)--node[swap,yshift=2mm, xshift=2mm] {$p_3$}(B3);\draw[->](B3)--node[swap,yshift=2mm] {$p_2$}(B2);\draw[->] (A2)--node[swap,yshift=2mm] {$q$}(A1);\draw[->] (C2)--node[swap,yshift=-2mm] {$r$}(C1);
\draw[->] (B1)--node[swap,xshift=-2mm] {$\bar{a'}$}(A1);\draw[->] (A1)--node[swap,yshift=2mm] {$a'$}(B3);\draw[->] (B3)--node[swap,xshift=2mm] {$\bar{a}$}(A2);\draw[->] (A2)--node[swap,xshift=2mm] {$a$}(B4);
\draw[->] (B2)--node[swap,xshift=-2mm,yshift=-2mm] {$\bar{c'}$}(C1);\draw[->] (C1)--node[swap,yshift=-1mm,xshift=1mm] {$c'$}(B3);\draw[->] (B3)--node[swap,yshift=2mm] {$\bar{c}$}(C2);\draw[->] (C2)--node[swap,xshift=3mm] {$c$}(B4);
\draw[->] (A1)--node[swap,xshift=2mm] {$b'$}(C1);\draw[->] (C1)--node[swap,yshift=3mm] {$\bar{b}$}(A2);\draw[->] (A2)--node[swap,yshift=3mm] {$b$}(C2);
\end{tikzpicture}}\]
with potential $W=p_1p_2a'\bar{a'}-q\bar{a}a'+p_3a\bar{a}-q\bar{b}b'+rb\bar{b}+p_2c'\bar{c'}-r\bar{c}c'+p_3c\bar{c}$, and frozen arrows $a$, $b$ and $c$. That is we have an isomorphism
$$\End_{\Pi}(M_{\mathbf{w}})\simeq kQ/\langle \partial_xW, x\in Q_1, x\neq a,b,c\rangle.$$

The orientation of $\Delta$ is given by the quiver $Q$ restricted to vertices corresponding to projective-injectives, that is in this example $6$, $7$ and $8$. Hence we have 
$$\Delta=\xymatrix@=.3cm{&3\ar[dr]^c&\\1\ar[rr]_b\ar[ur]^{a}&&2} \textrm{ and } \bar{\Delta}=\xymatrix@=.8cm{&3\ar@<.5ex>[dr]|0\ar@<.5ex>[dl]|1&\\1\ar@<.5ex>[rr]|0\ar@<.5ex>[ur]|0&&2\ar@<.5ex>[ll]|1 \ar@<.5ex>[ul]|1 \,.\!\!\!}$$

With this grading on $\Pi$, the degree of the arrows of the graded algebra $\End_\Pi(M_{\mathbf{w}})$ are
$$d(p_1)=d(p_2)=d(p_3)=d(q)=d(r)=0,$$ 
$$d(a)=d(a')=d(b)=d(b')=d(c)=d(c')=0,$$
$$\textrm{and}\ d(\bar{a})=d(\bar{a'})=d(\bar{b})=d(\bar{c})=1.$$

Now we have natural algebra morphisms 

\[\scalebox{0.7}{
\begin{tikzpicture}[>=stealth,scale=.6]
\node (B1) at (0,0){$1$};
\node (B2) at (2,0){$3$};
\node (B3) at (6,0){$5$};
\node (B4) at (10,0){$8$};
\node (A1) at (1,2){$2$};
\node(A2) at (7,2){$6$};
\node (C1) at (4,-1){$4$};
\node (C2) at (8,-1){$7$};
\draw[->] (B2)--node[fill=white,inner sep=.2mm] {$0$}(B1);\draw[->] (B4)--node[fill=white,inner sep=.2mm] {$0$}(B3);\draw[->](B3)--node[fill=white,inner sep=.2mm] {$0$}(B2);\draw[->] (A2)--node[fill=white,inner sep=.2mm] {$0$}(A1);\draw[->] (C2)--node[fill=white,inner sep=.2mm] {$0$}(C1);
\draw[->] (B1)--node[fill=white,inner sep=.2mm] {$1$}(A1);\draw[->] (A1)--node[fill=white,inner sep=.2mm] {$0$}(B3);\draw[->] (B3)--node[fill=white,inner sep=.2mm] {$1$}(A2);\draw[->] (A2)--node[fill=white,inner sep=.2mm] {$0$}(B4);
\draw[->] (B2)--node[fill=white,inner sep=.2mm] {$1$}(C1);\draw[->] (C1)--node[fill=white,inner sep=.2mm] {$0$}(B3);\draw[->] (B3)--node[fill=white,inner sep=.2mm] {$1$}(C2);\draw[->] (C2)--node[fill=white,inner sep=.2mm] {$0$}(B4);
\draw[->] (A1)--node[fill=white,inner sep=.2mm] {$0$}(C1);\draw[->] (C1)--node[fill=white,inner sep=.2mm] {$1$}(A2);\draw[->] (A2)--node[fill=white,inner sep=.2mm] {$0$}(C2);

\draw[>->, very thick] (-3,1)--(0,1);

\node (B1') at (-13,0){$1$};
\node (B2') at (-11,0){$3$};
\node (B3') at (-7,0){$5$};
\node (B4') at (-3,0){$8$};
\node (A1') at (-12,2){$2$};
\node(A2') at (-6,2){$6$};
\node (C1') at (-9,-1){$4$};
\node (C2') at (-5,-1){$7$};
\draw[->] (B2')--(B1');\draw[->] (B3')--(B4');\draw[->](B3')--(B2');\draw[->] (A2')--(A1');\draw[->] (C2')--(C1');
\draw[->] (A1')--(B3');\draw[->] (A2')--(B4');
\draw[->] (C1')--(B3');\draw[->] (C2')--(B4');
\draw[->] (A1')--(C1');\draw[->] (A2')--(C2');
\draw[->>, very thick] (-13,1)--(-16,1);
\node (B1'') at (-22,0){$1$};
\node (B2'') at (-20,0){$3$};
\node (B3'') at (-16,0){$5$};

\node (A1'') at (-21,2){$2$};

\node (C1'') at (-18,-1){$4$};

\draw[->] (B2'')--(B1'');\draw[->](B3'')--(B2'');
\draw[->] (A1'')--(B3'');
\draw[->] (C1'')--(B3'');
\draw[->] (A1'')--(C1'');
\draw[loosely dotted, thick] (B1'')..controls (-18,.5)..(A1'');
\draw[loosely dotted, thick] (B2'')..controls (-17,-.3)..(C1'');

\node(A) at (-20,-3){$\underline{\Lambda}$};
\node (B) at (-8,-3){$\Lambda$};
\node (C) at (2,-3){$\End_\Pi(M_{\mathbf{w}})$};

\end{tikzpicture}}\]
where $\Lambda:=\End_{\gr \Pi}(M_{\mathbf{w}})$ and $\underline{\Lambda}=\underline{\End}_{\gr\Pi}(M_{\mathbf{w}})$.
Then we consider the following triangle functor $$\xymatrix{F:\Db(\underline{\Lambda})\ar[rr]^-{Res}& & \Db(\Lambda)\ar[rr]^{-\lten_{\Lambda}M_{\mathbf{w}}} && \Db(\Sub\Pi_w)}$$ induced by these algebra morphisms and check that it satisfies the hypotheses of Theorem \ref{univ prop}.
\end{example}




\subsection{Cohen Macaulay modules over isolated singularities}\label{subsection CM}
In this section we assume that the characteristic of $k$ is zero.

Let $S$ be the polynomial ring $k[x_0,x_1,\ldots,x_d]$ in $d+1$ variables. Let $G$ be the finite cyclic subgroup of order $n$ of $\SL_{d+1}(k)$ generated by $g=\rm{diag}(\zeta^{a_0},\ldots,\zeta^{a_d})$, where $\zeta$ is a primitive $n$-root of unity. We assume moreover that 
$a_0+a_1+\ldots a_d=n$. 

The group $G$ acts naturally on $S$.  We denote by $S^G$ the invariant algebra and by $S*G$ the skew group algebra. The algebra $S^G$ is a Gorenstein isolated singularity of Krull dimension $d+1$.

\begin{theorem}
In the setup above, the following assertions hold.
\begin{itemize}
\item[(a)] \cite{Aus78} The stable category $\underline{\CM}(S^G)$ of
maximal Cohen-Macaulay $S^G$-modules is a $d$-CY triangulated category.
\item[(b)] \cite{Iya07a} The $S^G$-module $S$ is a
$d$-cluster tilting object in $\underline{\CM}(S^G)$.
\item[(c)] \cite{Aus86}, \cite{Yos90} We have $\End_{S^G}(S)\simeq S*G$.
\end{itemize}
\end{theorem}

If we want to realize $\underline{CM}(S^G)$ as a cluster category, we have to find a natural $\ZZ$-grading on $S^G$.   If we set $\deg(x_i)=\frac{a_i}{n}$, we obtain a structure of $\frac{\mathbb{Z}}{n}$-graded algebra on $S=\bigoplus_{i\in\ZZ}S_{\frac{i}{n}}$. Then the invariant subring $S^G$ corresponds to the subalgebra of $S$ of the direct sum of homogeneous components of integral degrees:
$$S^G=\bigoplus_{i\in\ZZ}S_i.$$
Hence we obtain  a structure of $\ZZ$-graded algebra for $S^G$.

Now we define, for $j=0,\ldots,n-1$, the graded $S^G$-modules $T_j=\bigoplus_{i\in\ZZ}S_{i+\frac{j}{n}}$. Then we have $T_0=S^G$ and $T=\bigoplus_{j=0}^{n-1}T_j\simeq S$ as $S^G$-modules.

The following result is also proved in \cite{TVdB2} using different techniques.
\begin{theorem}[\cite{AIR}]
Let $S$ and $G$ be as above. Let $\underline{A}$ be the stable endomorphism algebra $\underline{A}:=\underline{\End}_{\Gr S^G}(T).$ Then the following assertions hold.
\begin{itemize}
\item[(a)] The algebra $\underline{A}$ is a finite dimensional, $\tau_d$-finite algebra of global dimension at most $d$. 
\item[(b)] There exists a triangle equivalence $\Cc_d(\underline{A})\simeq \underline{CM}(S^G)$ sending the $d$-cluster-tilting object $\pi\underline{A}$ to $T$.
\end{itemize}
\end{theorem}
 
\begin{example}
Let $d=3$ and $G$ be the subgroup generated by $\frac{1}{5}(1,2,2)$. Then the skew-group algebra $S*G\simeq\End_{S^G}(S)$ is presented by the McKay quiver
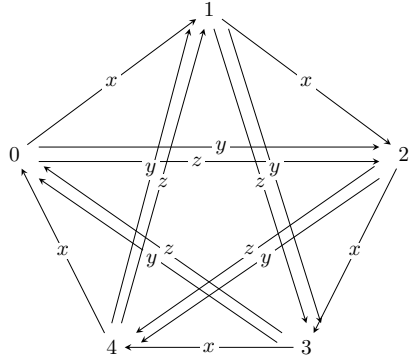
\[\scalebox{.8}{
\begin{tikzpicture}[>=stealth, scale=.8]
\node (P0) at (0,0){$0$}; \node (P1) at (4,3){$1$}; \node (P2) at
(8,0){$2$}; \node (P3) at (6,-4){$3$}; \node (P4) at (2,-4){$4$};

\draw [->] (0.5,0.15) -- node [fill=white,inner
sep=.5mm,xshift=2mm]{$y$} (7.5,0.15); \draw [->] (0.5,-0.15) -- node
[fill=white,inner sep=.5mm, xshift=-2mm]{$z$} (7.5,-0.15); \draw
[->] (7.4,-0.25) -- node [fill=white,inner sep=.5mm,
xshift=-1mm]{$z$} (2.5,-3.7); \draw [->] (7.5,-.5) -- node
[fill=white,inner sep=.5mm, xshift=1mm]{$y$} (2.6,-3.9); \draw [->]
(2,-3.5) -- node [fill=white,inner sep=.5mm, yshift=1mm]{$y$}
(3.6,2.6); \draw [->] (2.2,-3.5) -- node [fill=white,inner sep=.5mm,
yshift=-1mm]{$z$} (3.9,2.6); \draw [->] (4.4,2.6) -- node
[fill=white,inner sep=.5mm, yshift=1mm]{$y$} (6.3,-3.5); \draw [->]
(4.1,2.6) -- node [fill=white,inner sep=.5mm, yshift=-1mm]{$z$}
(6,-3.5); \draw [<-] (0.6,-0.25) -- node [fill=white,inner sep=.5mm,
xshift=1mm]{$z$} (5.5,-3.7); \draw [<-] (0.5,-.5) -- node
[fill=white,inner sep=.5mm, xshift=-1mm]{$y$} (5.4,-3.9);

\draw [->] (P0)-- node [fill=white,inner sep=.5mm]{$x$}(P1); \draw
[->] (P1)-- node [fill=white,inner sep=.5mm]{$x$}(P2); \draw [->]
(P2)-- node [fill=white,inner sep=.5mm]{$x$}(P3); \draw [->] (P3)--
node [fill=white,inner sep=.5mm]{$x$}(P4); \draw [->] (P4)-- node
[fill=white,inner sep=.5mm]{$x$}(P0);
\end{tikzpicture}}
\]
with the commutativity relations $xy=yx$, $yz=zy$, $zx=xz$. This is a Jacobian algebra associated with a non degenerate QP by \cite{TVdB1}. 
The vertex $j=0,\ldots,5$ corresponds to the summand $T_j$ of $T$.
With the grading defined as above, an arrow $i\to j$ has degree $0$ if $i<j$ and degree $1$ if $i>j$.
Hence the algebra $\End_{\Gr S^G}(T)$ is presented by
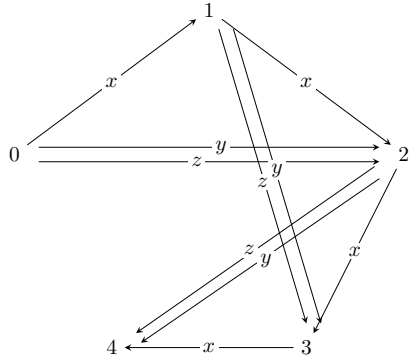
\[\scalebox{.8}{
\begin{tikzpicture}[>=stealth, scale=.8]
\node (P0) at (0,0){$0$}; \node (P1) at (4,3){$1$}; \node (P2) at
(8,0){$2$}; \node (P3) at (6,-4){$3$}; \node (P4) at (2,-4){$4$};

\draw [->] (0.5,0.15) -- node [fill=white,inner
sep=.5mm,xshift=2mm]{$y$} (7.5,0.15); \draw [->] (0.5,-0.15) -- node
[fill=white,inner sep=.5mm, xshift=-2mm]{$z$} (7.5,-0.15); \draw
[->] (7.4,-0.25) -- node [fill=white,inner sep=.5mm,
xshift=-1mm]{$z$} (2.5,-3.7); \draw [->] (7.5,-.5) -- node
[fill=white,inner sep=.5mm, xshift=1mm]{$y$} (2.6,-3.9);

\draw [->] (4.5,2.6) -- node [fill=white,inner sep=.5mm,
yshift=1mm]{$y$} (6.3,-3.5); \draw [->] (4.2,2.6) -- node
[fill=white,inner sep=.5mm, yshift=-1mm]{$z$} (6,-3.5);

\draw [->] (P0)-- node [fill=white,inner sep=.5mm]{$x$}(P1); \draw
[->] (P1)-- node [fill=white,inner sep=.5mm]{$x$}(P2); \draw [->]
(P2)-- node [fill=white,inner sep=.5mm]{$x$}(P3); \draw [->] (P3)--
node [fill=white,inner sep=.5mm]{$x$}(P4);

\end{tikzpicture}}
\]
with the commutativity relations. Therefore the stable endomorphism algebra $\underline{A}=\underline{\End}_{\Gr S^G}(T)\simeq \End_{\Gr S^G}(T_1\oplus\cdots\oplus T_4)$ is presented by the quiver
 \[\scalebox{.8}
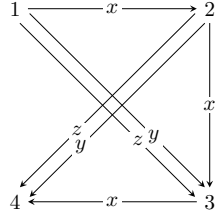
{
\begin{tikzpicture}[>=stealth, scale=.8]

\node (P1) at (2,0){$1$}; \node (P2) at (6,0){$2$}; \node (P3) at
(6,-4){$3$}; \node (P4) at (2,-4){$4$};

\draw [->] (5.7,-0.1) -- node [fill=white,inner sep=.5mm,
xshift=-5mm, yshift=-5mm]{$z$} (2.1,-3.7);
\draw [->] (5.9,-0.3) -- node [fill=white,inner sep=.5mm,
xshift=-6mm, yshift=-6mm]{$y$} (2.3,-3.9);


\draw [->] (2.3,-0.1) -- node [fill=white,inner sep=.5mm,
yshift=-6mm, xshift=6mm]{$y$} (5.9,-3.7); \draw [->] (2.1,-0.3) -- node
[fill=white,inner sep=.5mm, yshift=-5mm, xshift=5mm]{$z$} (5.7,-3.9);

\draw [->] (P1)-- node [fill=white,inner sep=.5mm]{$x$}(P2); \draw
[->] (P2)-- node [fill=white,inner sep=.5mm]{$x$}(P3); \draw [->]
(P3)-- node [fill=white,inner sep=.5mm]{$x$}(P4);

\end{tikzpicture}}\]
with the commutativity relations. By the previous theorem the category $\underline{\CM}(S^G)$ is triangle equivalent to the generalized cluster category $\Cc_2(\underline{A})$ .

\end{example}

 The proof of this theorem uses again Theorem \ref{univ prop}. One fundamental step in the proof consists in using the fact that the  skew-group algebra with the grading defined as above is bimodule $d$-Calabi-Yau of Gorenstein parameter 1 by \cite{BSW10}, and then in showing the following:
 
\begin{theorem}[\cite{AIR}]
Let $B$ be a $\ZZ$-graded algebra bimodule $d$-Calabi-Yau of Gorenstein parameter 1 such that $\dim_kB_i$ is finite for all $i$, and such that $B\in\per B^{\rm e}$. Then the DG algebra $\BPi_d(A)$ has its homology concentrated in degree $0$. Moreover, there is an isomorphism of $\ZZ$-graded algebras $\Pi_d(A)\simeq B.$
\end{theorem}

We end this section by asking the following intriguing questions:
\begin{question}
\begin{itemize}
\item[1.]  What can we say about the stable category $\underline{\CM}(S^G)$ when $G$ is not cyclic ? 
\item[2.] Is there an analogue of the categories $\Sub \Pi_{w}$ in higher CY-dimensions? 
\end{itemize}
\end{question}

\section{On the $\mathbb{Z}$-grading on the $3$-preprojective algebra}

Throughout the section $\Lambda$ is a finite dimensional algebra of global dimension at most two.

As we already saw in section 3, the $3$-preprojective algebra $\Pi_3(\Lambda)$ is isomorphic to the tensor algebra $\Talg_\Lambda \Ext^2_{\Lambda}(D\Lambda, \Lambda)$ (see \ref{equation isomorphism}) and thus is naturally positively $\mathbb{Z}$-graded as a tensor algebra. The algebra $\Lambda$ can easily be recovered from the graded algebra $\Pi_3(\Lambda)$ as its degree zero subalgebra. We already used this remark to prove equivalences between certain stable $2$-CY categories and cluster categories $\Cc_2(\Lambda)$ in subsections~\ref{subsection words} and~\ref{subsection CM}.

In this section, we study this grading, and outline applications in representation theory.

\subsection{Grading of $\Pi_3(\Lambda)$}

Let us describe more explicitly this grading. The algebra $\Pi_3(\Lambda)$ is isomorphic to the Jacobian algebra $\Jac(\bar{Q}_{\Lambda},\bar{W}_{\Lambda})$, where $\bar{Q}_\Lambda$ and $\bar{W}$ are obtained from the quiver $Q_\Lambda$  and the minimal relations of $\Lambda$ (cf. isomorphism \ref{equation jacobian}). This natural grading comes from a grading on $\bar{Q}_{\Lambda}$ defined as follows: 
\begin{itemize}
\item the arrows of $Q_{\Lambda}$ are of degree $0$, (indeed, they correspond to elements of $\Lambda$ which is the degree $0$ subalgebra of $\Pi_3(\Lambda)$); 
\item the new arrows, corresponding to minimal relations, are of degree $1$, (indeed, they correspond to elements of the bimodule $\Ext^2_{\Lambda}(D(\Lambda),\Lambda)$).
\end{itemize}
This grading on $\bar{Q}_\Lambda$ makes the potential $\bar{W}$ homogeneous of degree $1$. Hence the relations $\partial_a\bar{W}$ for $a\in Q_1$ are homogeneous.

\begin{example}
Let $\Lambda$ be as in Example \ref{example 3cycle}. The $3$-preprojective algebra $\Pi_3(\Lambda)$ is isomorphic, as a graded algebra, to the graded Jacobian algebra $\Jac(\bar{Q}_{\Lambda},W,d)$ with
$$\bar{Q}_{\Lambda}=\xymatrix@-0.3cm{ &2\ar[dl]_a& \\1\ar@<.5ex>[rr]^c\ar@<-.5ex>[rr]_{r_{ab}}&&3\ar[ul]_{b}}\quad d(a)=d(b)=d(c)=0,\quad d(r_{ab})=1\quad\textrm{and }W=r_{ab}ab.$$ 
\end{example}

Since we have an isomorphism of graded algebras $\Pi_3(\Lambda)\simeq \bigoplus_{p\in \ZZ}\Hom_{\Db(\Lambda)}(\Lambda,\SSS_2^{-p}\Lambda)$ (cf. \ref{equation isomorphism}), the subcategory $$\Uu_{\Lambda}:=\pi^{-1}(\Lambda)=\{ \SSS_{2}^{-p}\Lambda, \ p\in\ZZ\}\subset \Db(\Lambda)$$ provides a $\ZZ$-covering  of the graded algebra  $\Pi_3(\Lambda)$. 

Hence the study of the category $\Uu_\Lambda=\pi^{-1}(\Lambda)$ is one of the fundamental tools for the study of the grading of $\Pi_3(\Lambda)$.

\begin{theorem}[Thm 1.22 \cite{Iya08}, Prop 5.4.2 \cite{Ami08}]\label{orbit cluster tilting}
If $\Lambda$ is $\tau_2$-finite, then the subcategory $\Uu_\Lambda$ is a cluster-tilting subcategory of $\Db(\Lambda)$.
\end{theorem}

The next result shows that, moreover, the category $\Db(\Lambda)$ is determined by the category $\Uu_\Lambda$.
 
 \begin{theorem}[Recognition Theorem, Thm 3.5 \cite{AO10a}]\label{recognition theorem}
 Let $\Tt$ be an algebraic triangulated category, with Serre functor $\SSS$ and with a cluster-tilting subcategory $\Vv$. Assume that there exists a $\tau_2$-finite algebra $\Lambda$ and an equivalence $f:\Vv\simeq \Uu_\Lambda=\add\{\SSS_2^p\Lambda, p\in \mathbb{Z}\}$ such that $f\circ \SSS_2\simeq \SSS\circ f [-2]$. Then $\Tt$ is triangle equivalent to $\Db(\Lambda)$.
 \end{theorem} 
This result is the key step for the proof of all results presented in the next two sections.
\subsection{Mutation of graded QPs}
In this section, we assume that $\Lambda$ is a $\tau_2$-finite algebra of global dimension at most~2.

The observation that $\Pi_3(\Lambda)$ is a graded Jacobian algebra  leads us to study the notion of graded QP $(Q,W,d)$ with $W$ homogeneous of degree $1$. 
We can define the notion of reduction similarly to \cite{DWZ} (see \cite[Section 6]{AO10a}).

\begin{definition}[\cite{AO10a}, see also \cite{TVdB1}]
Let $(Q,W,d)$ be a graded QP such that $Q$ does not have any loops. Let $i$ be a vertex of $Q$. Then the \emph{left mutation} $\mu_i^L(Q,W,d)$ (resp. \emph{right mutation} $\mu_i^R(Q,W,d)$) of $(Q,W,d)$ at vertex $i$ is defined to be the reduction of the graded QP $(Q',W',d')$ constructed as follows:
\begin{itemize}
\item[(M1gr)] for each pair of arrows $\xymatrix{j\ar[r]^{a}& i\ar[r]^{b} & k}$, add an arrow $\xymatrix{j\ar[r]^{[ba]} & k}$ and put $d'([ba])=d(a)+d(b)$,
\item[(M2gr)] replace each arrow $\xymatrix{j\ar[r]^{a}& i}$ by an arrow $\xymatrix{j &i\ar[l]_{a^*}}$ and put $d'(a^*)=1-d(a)$ (resp. $d'(a^*)=-d(a)$),
\item[] replace each arrow $\xymatrix{ i\ar[r]^{b} & k}$ by an arrow $\xymatrix{i &k\ar[l]_{b^*}}$ and put $d'(b^*)=-d(b)$ (resp. $d'(b^*)=1-d(b)$).
\end{itemize}
All other arrows remain with the same degree and the potential $W'=[W]+W^*$ (as defined in \cite{DWZ}) and is again homogeneous of degree $1$.
\end{definition}

In the derived category $\Db(\Lambda)$, if $T$ is an object such that $\pi(T)\in \Cc_2(\Lambda)$ is cluster-tilting, then one can lift the exchange triangles of Theorem~\ref{theorem IY mutation} to triangles in $\Db(\Lambda)$. Then we obtain for any indecomposable direct summand $T_i$ of $T=T_0\oplus T_i$ the following triangles in $\Db(\Lambda)$ 
$$\xymatrix{T_i\ar[r] & B\ar[r]& T_i^L\ar[r] & T_i[1]} \quad \textrm{and}\quad \xymatrix{T_i^R\ar[r] & B'\ar[r]& T_i\ar[r] & T_i^R[1]}.$$ 
The objects $T_i^L$ and $T_i^R$ satisfy that $\pi(T_i^L)=\pi(T_i^R)=\pi(T_i)^*$ is the unique complement non isomorphic to $\pi(T_i)$ of the almost complete cluster-tilting object $\pi(T_0)$. We call the object $T_0\oplus T_i^L$ (resp. $T^R_i$) the \emph{left mutation} (resp. \emph{right mutation}) of $T$ at $T_i$.
 
Then we can formulate the graded analogue of Theorem \ref{BIRSm}, which  is a first motivation for the introduction of mutation of graded QP. 
\begin{theorem}[\cite{AO10a}]\label{gradedBIRSm}
Let $\Lambda$ and $T\in\Db(\Lambda)$ as above. Assume that there
exist a graded QP $(Q,W,d)$ with potential homogeneous of degree 1 such that we have an
isomorphism of graded algebras 
$$\xymatrix{\bigoplus_{p\in\ZZ}\Hom_{\Dd}(T,\SSS_2^{-p}T)\ar[r]^-\sim_-{\ZZ} & \Jac(Q,W,d).}$$
Let $T_i$ be an indecomposable summand of $T\simeq T_i\oplus T_0$ and assume that there are neither loops nor $2$-cycles incident to $i$ (corresponding to $T_i$) in the quiver of $\End_{\Cc}(T)$. 
Then there is an isomorphism of $\ZZ$-graded algebras
$$\xymatrix{\bigoplus_{p\in
  \mathbb{Z}}\Hom_\Dd(T_0\oplus T^L_i,\SSS_2^{-p}(T_0\oplus T^L_i))\ar[r]^-\sim_-{\ZZ}& \Jac(\mu_i^L(Q,W,d)).}$$  
\end{theorem}
\noindent
This result can be illustrated by the following picture
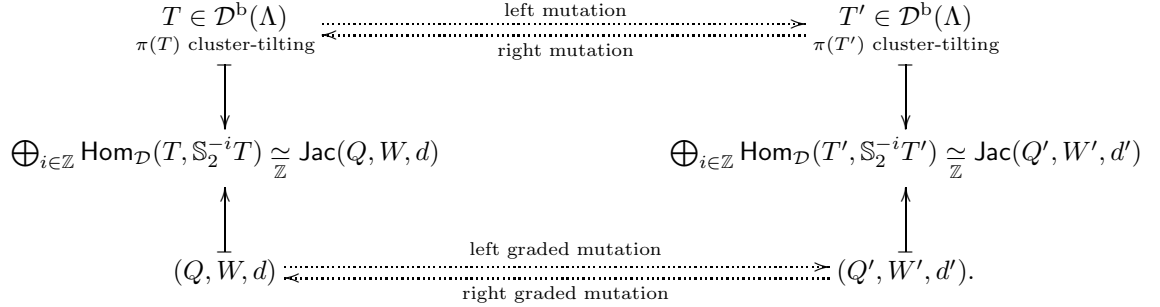
\[\xymatrix@C=1.3cm{\underset{\tiny \pi(T) \textrm{ cluster-tilting}}{T\in\Db(\Lambda)}\ar@<.5ex>@{..>}[rr]^{\tiny\textrm{left mutation}}\ar@{|->}[d] && \underset{\tiny\pi(T')\textrm{ cluster-tilting}}{T'\in\Db(\Lambda)}\ar@{|->}[d]\ar@<.5ex>@{..>}[ll]^{\tiny\textrm{right mutation}}\\ \bigoplus_{i\in\ZZ}\Hom_\Dd(T,\SSS_2^{-i}T)\underset{\ZZ}{\simeq} \Jac(Q,W,d) && \bigoplus_{i\in \ZZ}\Hom_\Dd(T',\SSS_2^{-i}T')\underset{\ZZ}{\simeq} \Jac(Q',W',d')\\ (Q,W,d)\ar@<.5ex>@{..>}[rr]^{\tiny\textrm{left graded mutation}}\ar@{|->}[u] && (Q',W',d')\ar@{|->}[u]  \ar@<.5ex>@{..>}[ll]^{\tiny\textrm{right graded mutation}}   .}\]

Our main application of the graded mutation is given by the following result, which gives a combinatorial criterion to see when two algebras of global dimension at most two are derived equivalent.
\begin{theorem}[\cite{AO10a}]\label{graded mutation and derived equivalence}
Let $\Lambda$ and $\Lambda '$ be two algebras of global dimension 2 , which are $\tau_2$-finite. Assume that one can pass from the graded QP $(\bar{Q},\bar{W},d)$ associated with $\Pi_3(\Lambda)$ to the graded QP $(\bar{Q}',\bar{W}',d')$ (up to graded right equivalence) associated with $\Pi_3(\Lambda')$ using a finite sequence of left and right mutations.  Then the algebras $\Lambda$ and $\Lambda'$ are derived equivalent.
\end{theorem}
This result can be seen as a generalization (in one direction) of Happel's Theorem (Thm \ref{Happel}). So we could ask what is the meaning of the hypothesis of $\tau_2$-finiteness for the algebras $\Lambda$ and $\Lambda'$ and wether this hypothesis is necessary or not. Indeed, cluster categories and cluster-tilting theory are hidden in the statement. However, the proof uses strongly Theorems \ref{orbit cluster tilting}, \ref{recognition theorem} and 
\ref{gradedBIRSm}, where the hypothesis of $\tau_2$-finiteness is fundamental.

\begin{example}
Let us illustrate this result with an example.
The following graded quivers with potential are linked by sequences of left and right graded mutations
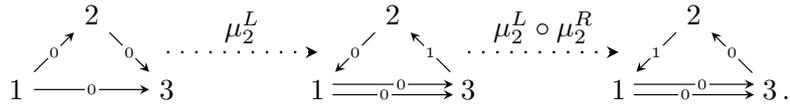
\[\scalebox{1.1}{
\begin{tikzpicture}[scale=.9,>=stealth]
\node (A1) at (0,0) {$1$};
\node (A2) at (1,1){$2$};
\node (A3) at (2,0){$3$};

\draw [->] (A1)--node [fill=white,inner sep=.2mm]{\tiny{0}}(A2);
\draw [->] (A2)--node [fill=white,inner sep=.2mm]{\tiny{0}}(A3);
\draw[->] (A1)--node [fill=white,inner sep=.2mm]{\tiny{0}}(A3);

\node (E1) at (4,0) {$1$};
\node (E2) at (5,1){$2$};
\node (E3) at (6,0){$3$};

\draw[loosely dotted, thick, ->] (2, 0.5) --node[yshift=3mm]{$\mu_2^L$} (4,0.5);
\draw[loosely dotted, thick, ->] (6, 0.5) --node[yshift=3mm]{$\mu_2^L\circ \mu^R_2$} (8,0.5);

\draw[->] (4.2,0.07)--node [fill=white,inner sep=.2mm,xshift=1mm]{\tiny{0}}(5.8,0.07);
\draw[->] (4.2,-0.07)--node [fill=white,inner sep=.2mm,xshift=-1mm]{\tiny{0}}(5.8,-0.07);
\draw[->] (E3)--node [fill=white,inner sep=.2mm]{\tiny{1}}(E2);
\draw[->] (E2)--node [fill=white,inner sep=.2mm]{\tiny{0}}(E1);

\node (F1) at (8,0) {$1$};
\node (F2) at (9,1){$2$};
\node (F3) at (10,0){$3$\,.\!\!\!};

\draw[->] (8.2,0.07)--node [fill=white,inner sep=.2mm,xshift=1mm]{\tiny{0}}(9.8,0.07);
\draw[->] (8.2,-0.07)--node [fill=white,inner sep=.2mm,xshift=-1mm]{\tiny{0}}(9.8,-0.07);
\draw[->] (F3)--node [fill=white,inner sep=.2mm]{\tiny{0}}(F2);
\draw [->] (F2)--node [fill=white,inner sep=.2mm]{\tiny{1}}(F1);

\end{tikzpicture}}\]

Here, the potential is $W=0$ for the first  quiver, and is $W=c$, where $c$ is one of the two 3-cycles, for the next two quivers.

Hence, using Theorem \ref{graded mutation and derived equivalence}, one deduces that the following algebras of global dimension at most~2 (presented by a quiver with relations) are in the same derived equivalence class. In this picture, a dotted line between two vertices means that a path between these two vertices is zero
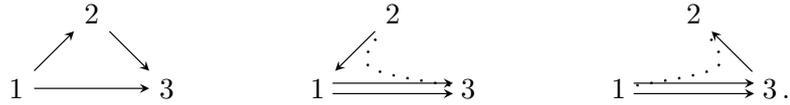
\[\scalebox{1.1}{
\begin{tikzpicture}[scale=.9,>=stealth]
\node (A1) at (0,0) {$1$};
\node (A2) at (1,1){$2$};
\node (A3) at (2,0){$3$};

\draw [->] (A1)--(A2);
\draw [->] (A2)--(A3);
\draw [->] (A1)--(A3);

\node (E1) at (8,0) {$1$};
\node (E2) at (9,1){$2$};
\node (E3) at (10,0){$3$\,.\!\!\!};

\draw[->] (8.2,0.07)--(9.8,0.07);
\draw[->] (8.2,-0.07)--(9.8,-0.07);
\draw[->] (E3)--(E2);
\draw [loosely dotted, thick] (E1)..controls(9.5,0.25)..(E2);

\node (F1) at (4,0) {$1$};
\node (F2) at (5,1){$2$};
\node (F3) at (6,0){$3$};

\draw[->] (4.2,0.07)--(5.8,0.07);
\draw[->] (4.2,-0.07)--(5.8,-0.07);
\draw[loosely dotted, thick] (F3)..controls (4.5, 0.25)..(F2);
\draw [->] (F2)--(F1);

\end{tikzpicture}}\]
\end{example}

\subsection{Application to cluster equivalence}
In this subsection, we describe another application of graded mutation to the notion of cluster equivalence. 
\begin{definition}
Two finite dimensional algebras $\Lambda$ and $\Lambda'$ of global dimension at most $2$ which are $\tau_2$-finite are called \emph{cluster equivalent}  if there exists a triangle equivalence $\Cc_2(\Lambda)\simeq \Cc_2(\Lambda')$ between their cluster categories.
\end{definition}

The notion of cluster equivalence seems to be a reasonable way to relate the homological algebras of two algebras of global dimension ~2 (see Theorem \ref{graded derived equivalence} below). This notion is more general than derived equivalence: two derived equivalent algebras of global dimension at most two are cluster equivalent, the converse is not true, as shown in the following example.

\begin{example}\label{example 3cycle bis} We take $\Lambda$ as in Example~\ref{example 3cycle}. It is a $\tau_2$-finite algebra of global dimension at most~$2$. The potential $\bar{W}=r_{ab}ab$ is rigid. The object $\pi(\Lambda)=\pi(P_1\oplus P_2\oplus P_3)$ is a cluster-tilting object with endomorphism algebra $\Pi_3(\Lambda)\simeq \Jac(\bar{Q}_\Lambda, \bar{W})$. By Theorem~\ref{BIRSm}, if we mutate $\pi(\Lambda)$ at $\pi(P_2)$, we obtain a new cluster-tilting object whose endomorphism algebra is  the Jacobian algebra $\Jac(Q,0)$, where $Q$ is the acyclic quiver $$\xymatrix@-0.5cm{ &&2\ar[dr]& \\Q=&1\ar[rr]\ar[ur]&&3.}$$
By Theorem \ref{KR recognition theorem}, the cluster category $\Cc_2(\Lambda)$ is equivalent to the acyclic category $\Cc_{Q}$. Therefore the algebra $\Lambda$ is cluster equivalent to the path algebra $kQ$. However, the algebra $\Lambda$ is not piecewise hereditary since its quiver $Q_{\Lambda}$ contains an oriented cycle. The algebras $\Lambda$ and $kQ$ are cluster equivalent, but not derived equivalent.
\end{example}

The next result gives a better understanding the notion of cluster equivalence. 
\begin{theorem}\label{graded derived equivalence}
Let $\Lambda$ and $\Lambda '$ be two $\tau_2$-finite algebras of global dimension at most~$2$. Denote by $(\bar{Q},\bar{W},d)$ (resp. $(\bar{Q}',\bar{W}',d')$) the graded QP associated with $\Lambda$ (resp. $\Lambda'$). Assume that there exists a sequence $\mu$ of mutations such that  \begin{itemize}
\item the QPs $(\bar{Q}',\bar{W}')$ and $\mu(\bar{Q},\bar{W})$ are right equivalent,
\item this right equivalence sends homogeneous elements of $\widehat{k\bar{Q'}}$ (w.r.t. the grading $d'$) to homogeneous elements of $\widehat{k\mu(\bar{Q})}$ (w.r.t. the grading $\mu^L(d)$). 
\end{itemize}
Then, the following statements hold
\begin{itemize}
\item this sequence $\mu$ of mutations induces a grading on $\Lambda$ and on $\Lambda'$, 
\item we have a triangle equivalence $\Db(\gr \Lambda)\simeq \Db(\gr\Lambda')$,
\item the algebras $\Lambda$ and $\Lambda'$ are cluster equivalent.
\end{itemize}
\end{theorem}
The proof uses a graded version of the Recognition Theorem \ref{recognition theorem}, and therefore the hypothesis of $\tau_2$-finiteness for the algebras $\Lambda$ and $\Lambda'$ is again strongly needed in the proof.
\begin{example}
We take $\Lambda$ as in Example~\ref{example 3cycle} and $\Lambda'=kQ$, where $Q$ is the acyclic quiver $\mu_2(\bar{Q}_\Lambda)$. If we do the left graded mutation at vertex $2$ of the graded QP $(\bar{Q}_{\Lambda},\bar{W}_{\Lambda})$ we obtain the graded quiver:

\[\scalebox{1.1}{
\begin{tikzpicture}[scale=.9,>=stealth]
\node (A1) at (8,0) {$1$};
\node (A2) at (9,1){$2$};
\node (A3) at (10,0){$3$\,.\!\!\!};

\draw [->] (A1)--node [fill=white,inner sep=.2mm]{\tiny{0}}(A2);
\draw [->] (A2)--node [fill=white,inner sep=.2mm]{\tiny{1}}(A3);
\draw[->] (A1)--node [fill=white,inner sep=.2mm]{\tiny{0}}(A3);

\node (E1) at (4,0) {$1$};
\node (E2) at (5,1){$2$};
\node (E3) at (6,0){$3$};

\draw[loosely dotted, thick, ->] (6,0.5)--node[swap,yshift=3mm] {$\mu_2^L$}(8,0.5);

\draw[->] (4.2,0.07)--node [fill=white,inner sep=.2mm,xshift=1mm]{\tiny{0}}(5.8,0.07);
\draw[->] (4.2,-0.07)--node [fill=white,inner sep=.2mm,xshift=-1mm]{\tiny{1}}(5.8,-0.07);
\draw[->] (E3)--node [fill=white,inner sep=.2mm]{\tiny{0}}(E2);
\draw[->] (E2)--node [fill=white,inner sep=.2mm]{\tiny{0}}(E1);

\end{tikzpicture}}\]

It induces a grading on the algebra $kQ$. Now, if we mutate the graded QP $(Q,0,0)$ at vertex $2$ using the right graded mutation (which is the inverse of the left graded mutation at $2$), we obtain the graded quiver:
\[\scalebox{1.1}{
\begin{tikzpicture}[scale=.9,>=stealth]
\node (A1) at (0,0) {$1$};
\node (A2) at (1,1){$2$};
\node (A3) at (2,0){$3$};

\draw [->] (A1)--node [fill=white,inner sep=.2mm]{\tiny{0}}(A2);
\draw [->] (A2)--node [fill=white,inner sep=.2mm]{\tiny{0}}(A3);
\draw[->] (A1)--node [fill=white,inner sep=.2mm]{\tiny{0}}(A3);

\node (E1) at (4,0) {$1$};
\node (E2) at (5,1){$2$};
\node (E3) at (6,0){$3$\,.\!\!\!};

\draw[loosely dotted, thick, ->] (2,0.5)--node[swap,yshift=3mm] {$\mu_2^R$}(4,0.5);

\draw[->] (4.2,0.07)--node [fill=white,inner sep=.2mm,xshift=1mm]{\tiny{0}}(5.8,0.07);
\draw[->] (4.2,-0.07)--node [fill=white,inner sep=.2mm,xshift=-1mm]{\tiny{0}}(5.8,-0.07);
\draw[->] (E3)--node [fill=white,inner sep=.2mm]{\tiny{1}}(E2);
\draw[->] (E2)--node [fill=white,inner sep=.2mm]{\tiny{0}}(E1);
\end{tikzpicture}}\]
It induces a grading on $\Lambda$. Theorem~\ref{graded derived equivalence} implies that the graded algebras
\[\scalebox{1.1}{
\begin{tikzpicture}[scale=.9,>=stealth]
\node (A1) at (0,0) {$1$};
\node (A2) at (1,1){$2$};
\node (A3) at (2,0){$3$};

\draw [->] (A1)--node [fill=white,inner sep=.2mm]{\tiny{0}}(A2);
\draw [->] (A2)--node [fill=white,inner sep=.2mm]{\tiny{1}}(A3);
\draw[->] (A1)--node [fill=white,inner sep=.2mm]{\tiny{0}}(A3);

\node (E1) at (4,0) {$1$};
\node (E2) at (5,1){$2$};
\node (E3) at (6,0){$3$};

\draw[->] (E1)--node [fill=white,inner sep=.2mm,]{\tiny{0}}(E3);
\draw[->] (E3)--node [fill=white,inner sep=.2mm]{\tiny{1}}(E2);
\draw[->] (E2)--node [fill=white,inner sep=.2mm]{\tiny{0}}(E1);
\draw[loosely dotted, thick] (E1)..controls (5,0.9)..(E3);
\end{tikzpicture}}\]
are graded derived equivalent.

\end{example}

Our hope was initially to understand and characterize the notion of cluster equivalence. Unfortunately, this seems to be a very difficult question. If two algebras $\Lambda$ and $\Lambda'$ are cluster equivalent, then we get a generalized cluster category with two canonical cluster-tilting objects. The first difficulty is that these two cluster-tilting objects are not known to be linked by a sequence of mutations. If they are, using left and right mutation of graded QP, the grading on $\Pi_3(\Lambda)$ induces a grading on $\Pi_3(\Lambda')$ a priori not equivalent to the natural one. But two $\ZZ$-gradings on an algebra does not give rise in general to a $\ZZ^2$-grading on it and  this fact is needed in the proof.  We call it \emph{compatibility condition} in \cite[Def. 8.5]{AO10a}. This condition ensures that there is a right equivalence between $(\bar{Q}',\bar{W}')$ and $\mu(\bar{Q},\bar{W})$ which sends homogeneous elements of $\widehat{k\bar{Q'}}$ (w.r.t. the grading $d'$) to homogeneous elements of $\widehat{k\mu(\bar{Q})}$ (w.r.t. the grading $\mu^L(d)$).

\subsection{Algebras of acyclic cluster type}

In this subsection, we focus on the case of algebras which are cluster equivalent to an hereditary algebra, and deduce properties of these algebras using results on acyclic cluster categories. All the results of this subsection are in \cite{AO10b}.

\begin{definition}
A $\tau_2$-finite algebra $\Lambda$ of global dimension at most 2 is said to be of \emph{acyclic cluster type} when its cluster category $\Cc_2(\Lambda)$ is equivalent to some cluster category $\Cc_Q$, where $Q$ is an acyclic quiver.
\end{definition}

This setup is especially nice since in this case, the two problems listed above do not appear. On one hand, all cluster-tilting objects in $\Cc_2(\Lambda)$ are IY-mutation equivalent (Thm \ref{HU theorem}).
Hence one gets a combinatorial characterization of algebras of acyclic cluster type.
\begin{proposition}  
An  algebra $\Lambda$ of global dimension at most $2$ is $\tau_2$-finite and of acyclic cluster type if and only if the quiver $\bar{Q}_\Lambda$ of $\Pi_3(\Lambda)$ is mutation acyclic and the potential $\bar{W}$ is rigid.
\end{proposition}

Note that, if the QP $(\bar{Q}_\Lambda,\bar{W})$ satisfies the hypothesis of the proposition above, it implies the $\tau_2$-finiteness for the algebra $\Lambda$. Indeed, the property of Jacobi-finiteness for a QP is stable under DWZ-mutation \cite{DWZ}.

\medskip
On the other hand, all gradings are compatible in the case of acyclic cluster type algebras. This point comes from the fact that, up to automorphism of algebras, two $\ZZ$-gradings on a path algebra, yield a $\ZZ^2$-grading on it. Hence we get a converse for Theorem \ref{graded mutation and derived equivalence} in the setup of algebras of acyclic cluster type. This is a generalization of Happel's result (Theorem \ref{Happel}).

\begin{theorem}[\cite{AO10b}]\label{generalization Happel}
Let $\Lambda$ and $\Lambda'$ be two algebras of global dimension $2$. Let $(\bar{Q},\bar{W},d)$ (resp. $(\bar{Q}',\bar{W}',d')$) be the graded QP associated with $\Pi_3(\Lambda)$ (resp. $\Pi_3(\Lambda')$).  Assume that the quiver $\bar{Q}$ is mutation acyclic, and that the potential $\bar{W}$ is rigid. Then $\Lambda$ and $\Lambda'$ are derived equivalent if and only if  one can pass from the graded QP $(\bar{Q},\bar{W},d)$ to the graded QP $(\bar{Q}',\bar{W}',d')$ (up to graded right equivalence) using a finite sequence of left and right mutations.  
\end{theorem}

If $Q$ is an acyclic quiver whose underlying graph is a tree, then one easily checks that one can pass from any grading on $Q$ to the trivial one using left and right mutations. Applying Theorem \ref{generalization Happel} to an algebra $\Lambda$ of cluster type $Q$, and to $\Lambda'=kQ$, one gets the following consequence. 

\begin{corollary}
Let $\Lambda$ be an algebra of cluster type $Q$, where $Q$ is a tree. Then the algebra $\Lambda$ is derived equivalent to  $kQ$.
\end{corollary}

Theorem~\ref{generalization Happel} can also be used to describe algebras of a given cluster type and to classify these algebras up to derived equivalence (see example below).

\begin{example}
As shown in \cite{AO10b}, there are 11 algebras of cluster type $\widetilde{A}_{2,2}$ (up to isomorphism). These algebras belong to two derived equivalence classes. 
In the first one we find the 8 piecewise hereditary algebras:

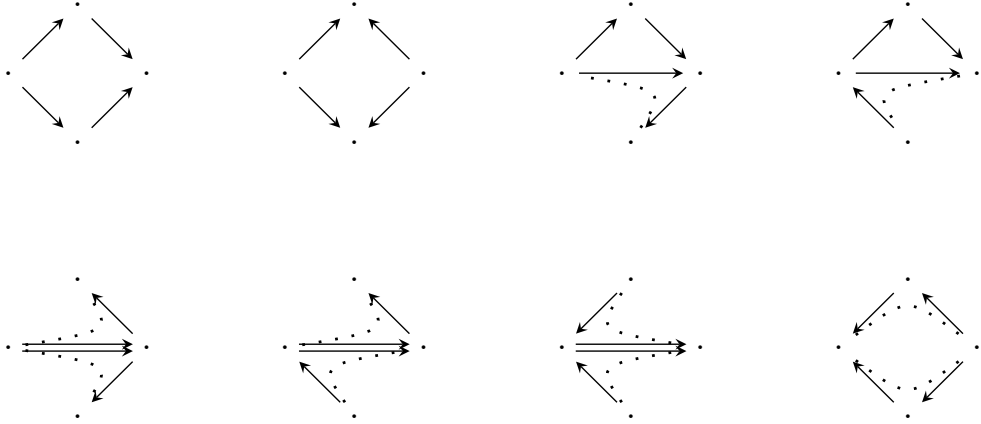
\[\scalebox{1.3}{
\begin{tikzpicture}[scale=.7,>=stealth]
\node (A1) at (0,0) {$.$};
\node (A2) at (1,1){$.$};
\node (A3) at (2,0){$.$};
\node (A4) at (1,-1){$.$};

\draw [->] (A1)--(A2);
\draw [->] (A2)--(A3);
\draw [->] (A1)--(A4);
\draw [->] (A4)--(A3);

\node (B1) at (4,0) {$.$};
\node (B2) at (5,1){$.$};
\node (B3) at (6,0){$.$};
\node (B4) at (5,-1){$.$};

\draw [->] (B1)--(B2);
\draw [->] (B3)--(B2);
\draw [->] (B1)--(B4);
\draw [->] (B3)--(B4);

\node (C1) at (8,0) {$.$};
\node (C2) at (9,1){$.$};
\node (C3) at (10,0){$.$};
\node (C4) at (9,-1){$.$};

\draw [->] (C1)--(C2);
\draw [->] (C2)--(C3);
\draw [->] (C3)--(C4); 
\draw [loosely dotted,thick] (C4)..controls (9.5,-0.25)..(C1);
\draw [->] (C1)--(C3);

\node (D1) at (12,0) {$.$};
\node (D2) at (13,1){$.$};
\node (D3) at (14,0){$.$};
\node (D4) at (13,-1){$.$};

\draw[->] (D1)--(D2);
\draw[->] (D2)--(D3);
\draw[->] (D1)--(D3);
\draw[->] (D4)--(D1);
\draw [loosely dotted, thick] (D3)..controls (12.5,-.25)..(D4);

\node (E1) at (0,-4) {$.$};
\node (E2) at (1,-3){$.$};
\node (E3) at (2,-4){$.$};
\node (E4) at (1,-5){$.$};

\draw[->] (0.2,-3.95)--(1.8,-3.95);
\draw[->] (0.2,-4.05)--(1.8,-4.05);
\draw[->] (E3)--(E2);
\draw[->] (E3)--(E4);
\draw [loosely dotted, thick] (E1)..controls(1.5,-3.75)..(E2);
\draw [loosely dotted, thick] (E1)..controls(1.5,-4.25)..(E4);

\node (F1) at (4,-4) {$.$};
\node (F2) at (5,-3){$.$};
\node (F3) at (6,-4){$.$};
\node (F4) at (5,-5){$.$};

\draw[->] (4.2,-3.95)--(5.8,-3.95);
\draw[->] (4.2,-4.05)--(5.8,-4.05);
\draw[->] (F3)--(F2);
\draw[->] (F4)--(F1);
\draw [loosely dotted, thick] (F1)..controls(5.5,-3.75)..(F2);
\draw [loosely dotted, thick] (F4)..controls(4.5,-4.25)..(F3);

\node (G1) at (8,-4) {$.$};
\node (G2) at (9,-3){$.$};
\node (G3) at (10,-4){$.$};
\node (G4) at (9,-5){$.$};

\draw[->] (8.2,-3.95)--(9.8,-3.95);
\draw[->] (8.2,-4.05)--(9.8,-4.05);
\draw[->] (G2)--(G1);
\draw[->] (G4)--(G1);
\draw [loosely dotted, thick] (G2)..controls (8.5,-3.75)..(G3);
\draw [loosely dotted, thick] (G4)..controls(8.5,-4.25)..(G3);

\node (H1) at (12,-4) {$.$};
\node (H2) at (13,-3){$.$};
\node (H3) at (14,-4){$.$};
\node (H4) at (13,-5){$.$};

\draw[->] (H2)--(H1);
\draw[->] (H3)--(H2);
\draw[->] (H4)--(H1);
\draw[->] (H3)--(H4);
\draw [loosely dotted, thick] (H1)..controls (13,-3.25)..(H3);
\draw [loosely dotted, thick] (H1)..controls (13,-4.75)..(H3);

\end{tikzpicture}}\]

The second equivalence class contains the 3 algebras which are not piecewise hereditary:

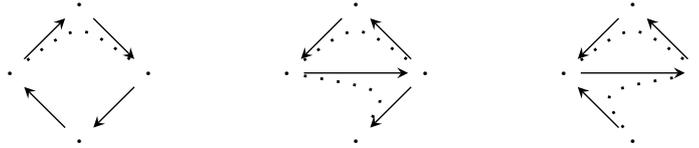
\[\scalebox{1.3}{
\begin{tikzpicture}[scale=0.7,>=stealth]
\node (A1) at (0,0) {$.$};
\node (A2) at (1,1){$.$};
\node (A3) at (2,0){$.$};
\node (A4) at (1,-1){$.$};

\draw [->] (A1)--(A2);
\draw [->] (A2)--(A3);
\draw [->] (A3)--(A4);
\draw [->] (A4)--(A1);
\draw [loosely dotted, thick]  (A1)..controls (1,0.75)..(A3);

\node (B1) at (4,0) {$.$};
\node (B2) at (5,1){$.$};
\node (B3) at (6,0){$.$};
\node (B4) at (5,-1){$.$};

\draw [->] (B2)--(B1);
\draw [->] (B3)--(B2);
\draw [->] (B1)--(B3);
\draw [->] (B3)--(B4);
\draw [loosely dotted, thick] (B1)..controls (5,0.75)..(B3);
\draw[loosely dotted, thick] (B1)..controls (5.5,-.25)..(B4);

\node (C1) at (8,0) {$.$};
\node (C2) at (9,1){$.$};
\node (C3) at (10,0){$.$};
\node (C4) at (9,-1){$.$};

\draw [->] (C2)--(C1);
\draw [->] (C3)--(C2);
\draw [->] (C1)--(C3); 
\draw [loosely dotted,thick] (C4)..controls (8.5,-0.25)..(C3);
\draw [->] (C4)--(C1);
\draw[loosely dotted, thick] (C1)..controls (9,0.75)..(C3);
\end{tikzpicture}}\]

\end{example}

\subsection{Gradable $\Pi_3{\Lambda}$-modules and piecewise hereditary algebras}
In the paper \cite{AO10c}, we study the functor $\pi:\Db(\Lambda)\to \Cc_2(\Lambda)$. In particular, we try to understand its image and to find conditions on the algebra $\Lambda$ to make it dense.

The following proposition gives another interpretation of the objects in the image of the functor $\pi$.
\begin{proposition}[\cite{AO10c}]
Let $\Lambda$ be a finite dimensional $\tau_2$-finite algebra. The following diagram commutes:
\[\xymatrix{\Db(\Lambda)\ar[rr]^{\pi_{\Lambda}}\ar@{->>}[d]_{\Hom_{\Dd}(\bigoplus_{i\in\mathbb{Z}}\SSS_2^{-i}\Lambda, -)} && \Cc_2(\Lambda)\ar@{->>}[d]^{\Hom_{\Cc}(\pi(\Lambda), -)}\\
\gr \Pi_3(\Lambda)\ar[rr] && \mod \Pi_3(\Lambda)\,,\!\!\!}\]
where the functor $\gr \Pi_3(\Lambda)\to \mod\Pi_3(\Lambda)$ is the forgetful functor. 
\end{proposition}

This means that the objects which are in the image of the triangle functor $\pi:\Db(\Lambda)\to \Cc_2(\Lambda)$ correspond to gradable modules when viewed in $\mod \Pi_3(\Lambda)$. As a consequence, we obtain that all rigid objects of $\Cc_2(\Lambda)$ are in the image of $\pi$~\cite[Cor 4.5]{AO10c}. Hence any cluster-tilting object $T$ in $\Cc_2(\Lambda)$ has a preimage in $\Db(\Lambda)$ and the subcategory $\pi^{-1}(\add T)\subset\Db(\Lambda)$ is a cluster-tilting subcategory of $\Db(\Lambda)$~\cite[Prop. 3.1]{AO10a}.

Following results of Gordon and Green \cite{GG} on graded algebras, one deduces results on the image of the functor $\pi$. 

\begin{corollary}[\cite{AO10c}]
The AR-quiver of the orbit category $\Db(\Lambda)/\SSS_2$ is the union of connected components of the AR-quiver of $\Cc_2(\Lambda)$.
\end{corollary}

\begin{example}
Let $\Lambda$ be as in Example \ref{example 3cycle}. We already saw in example \ref{example 3cycle bis} that it is an algebra of acyclic cluster type $Q$, where $Q$ is of type $\widetilde{A}_{2,1}$. In this example, we compute and compare the shape of the AR-quivers of the cluster category $\Cc_2(\Lambda)\simeq \Cc_Q=\Db(kQ)/\SSS_2$, and of the orbit category $\Db(\Lambda)/\SSS_2$.

\medskip
Since the cluster category $\Cc_Q$ is defined to be the orbit category of $\Db(kQ)/\SSS_2$, it is easy to compute the shape of its AR-quiver. The AR-quiver of $\Db(kQ)$ has three kinds of connected components (see \cite{Hap87}):
\begin{itemize}
\item A family, indexed by $\ZZ$, of connected components of type $\ZZ Q$, called \emph{transjective components}. They correspond to shifts of preprojective and preinjective $kQ$-modules.  
\item A family, indexed by $k^{\times}\times \ZZ$, of connected components of type $\ZZ A_{\infty}/\tau$ , called \emph{homogeneous tubes}.
\item  Two families, indexed by $\ZZ$,  one of tubes $\ZZ A_{\infty}/\tau^2$ of rank $2$,  and one of tubes $\ZZ A_{\infty}/\tau$ of rank 1, called \emph{exceptional tubes}.
\end{itemize}

Notice that in general, tubes of rank one are called \emph{homogeneous}, but in this example, one homogeneous can be considered as an exceptional one.

 The functor $\SSS_2\simeq \tau[-1]$ acts transitively on the family of components of type $\ZZ Q$, and sends a tube on its shift. Hence the shape of the AR-quiver of the category $\Cc_Q$ is the following:

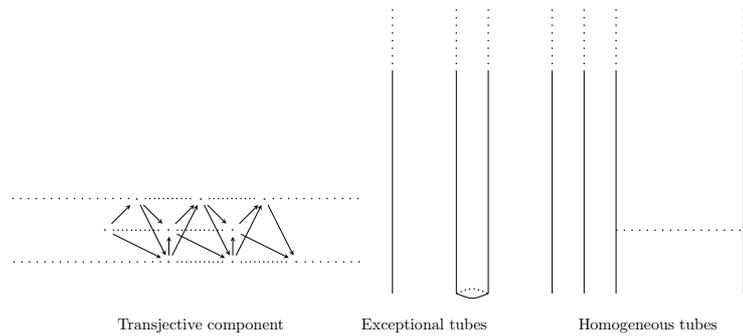
\[\scalebox{.6}{
\begin{tikzpicture}[scale=.7,>=stealth]

\node (A1) at (0,0) {$.$};
\node (A2) at (1,1) {$.$};
\node (A3) at (2,-1) {$.$};
\node (B1) at (2,0) {$.$};
\node (B2) at (3,1) {$.$};
\node (B3) at (4,-1) {$.$};
\node (C1) at (4,0) {$.$};
\node (C2) at (5,1) {$.$};
\node (C3) at (6,-1) {$.$};

\draw[->] (A1)--(A2);\draw[->] (A2)--(A3);\draw[->] (A1)--(A3);
\draw[->] (B1)--(B2);\draw[->] (B2)--(B3);\draw[->] (B1)--(B3);
\draw[->] (C1)--(C2);\draw[->] (C2)--(C3);\draw[->] (C1)--(C3);
\draw[->] (A2)--(B1);\draw[->] (A3)--(B2);\draw[->] (A3)--(B1);
\draw[->] (B2)--(C1);\draw[->] (B3)--(C2);\draw[->] (B3)--(C1);
\draw [dotted, thick] (A1)--(B1);\draw [dotted, thick] (B1)--(C1);\draw [dotted, thick] (A2)--(B2);\draw [dotted, thick] (B2)--(C2);\draw [dotted, thick] (A3)--(B3);\draw [dotted, thick] (B3)--(C3);

\draw[loosely dotted, thick] (A2)--(-3,1);\draw[loosely dotted, thick](A3)--(-3,-1); \draw[loosely dotted, thick] (C2)--(8,1); \draw[loosely dotted, thick] (C3)--(8,-1);

\draw (9,-2)--(9,5);
\draw[loosely dotted, thick] (9,5)--(9,7);
\draw (11,-2)--(11,5);
\draw[loosely dotted, thick] (11,5)--(11,7);
\draw (12,-2)--(12,5);
\draw[loosely dotted, thick] (12,5)--(12,7);
\draw (11,-2)..controls (11.5,-2.2)..(12,-2);\draw[dotted, thick] (11,-2)..controls (11.5,-1.8)..(12,-2);

\draw (14,-2)--(14,5);
\draw[loosely dotted, thick] (14,5)--(14,7);
\draw (15,-2)--(15,5);
\draw[loosely dotted, thick] (15,5)--(15,7);
\draw (16,-2)--(16,5);
\draw[loosely dotted, thick] (16,5)--(16,7);
\draw (20,-2)--(20,5);
\draw[loosely dotted, thick] (20,5)--(20,7);
\draw[loosely dotted, thick] (16,0)--(20,0);

\node (P1) at (3,-3){Transjective component};
\node (P2) at (10,-3){Exceptional tubes};
\node (P3) at (17,-3){Homogeneous tubes};

\end{tikzpicture}}\]

On the other hand, the AR-quiver of the derived category $\Db(\Lambda)$ has three connected components. One is of type $\ZZ A_{\infty}$ and contains the projective indecomposable $P_2={\bsm 2\\3\\1\\2\esm}$:
\[\scalebox{.8}{
\begin{tikzpicture}[scale=.8,>=stealth]

\node (P1) at (0,0) {${\bsm 2\\3\\1\\2\esm}[-1]$};
\node (P2) at (2,0) {\frame{${\bsm &2&\\&3&\\&1&\\&2&\esm}$}};
\node (P3) at (4,0) {${\bsm 2\\3\\1\\2\esm}[1]$};
\node (P4) at (6,0) {\frame{${\bsm 2\\3\\1\\2\esm}[2]=\SSS_2^{-1}{\bsm 2\\3\\1\\2\esm}$}};
\node (P5) at (8,0) {${\bsm 2\\3\\1\\2\esm}[3]$};

\node (Q1) at (0,2) {};
\node (Q2) at (2,2) {};
\node (Q3) at (4,2) {};
\node (Q4) at (6,2) {};
\node (Q5) at (8,2) {};

\node (R1) at (-1,1) {};
\node (R2) at (1,1) {};
\node (R3) at (3,1) {};
\node (R4) at (5,1) {};
\node (R5) at (7,1) {};
\node (R6) at (9,1) {};

\draw [->] (R1)--(Q1);
\draw [->] (R1)--(P1);
\draw [->] (R2)--(Q2);
\draw [->] (R2)--(P2);
\draw [->] (R3)--(Q3);
\draw [->] (R3)--(P3);
\draw [->] (R4)--(Q4);
\draw [->] (R4)--(P4);
\draw [->] (R5)--(Q5);
\draw [->] (R5)--(P5);

\draw [->] (P1)--(R2);
\draw [->] (Q1)--(R2);
\draw [->] (P2)--(R3);
\draw [->] (Q2)--(R3);
\draw [->] (Q3)--(R4);
\draw [->] (P3)--(R4);
\draw [->] (P4)--(R5);
\draw [->] (P5)--(R6);
\draw [->] (Q4)--(R5);
\draw [->] (Q5)--(R6);

\draw [loosely dotted, thick] (-3,1)--(R1);
\draw [loosely dotted, thick] (0,4)--(Q1);
\draw [loosely dotted, thick] (12,1)--(R6);
\draw [loosely dotted, thick] (8,4)--(Q5);
\end{tikzpicture}}
\]

One is of type $\ZZ A_{\infty}$ and contains the simple $\Lambda$-modules $S_1={\bsm 1\esm}$ and $S_3={\bsm 3\esm}$, and the module ${\bsm3\\1\esm}$: 
\[\scalebox{.8}{
\begin{tikzpicture}[scale=.8,>=stealth]

\node (P1) at (0,0) {\frame{${\bsm 3\esm}[1]$}};
\node (P2) at (2,0) {\frame{${\bsm&&\\ &1&\\&&\esm}$}};
\node (P3) at (4,0) {\frame{${\bsm 3\esm}=\SSS_2({\bsm 1 \esm})$}};
\node (P4) at (6,0) {\frame{${\bsm 1\esm}[-1]$}};
\node (P5) at (8,0) {\frame{${\bsm 3\esm}[-1]$}};

\node (Q1) at (0,2) {};
\node (Q2) at (2,2) {};
\node (Q3) at (4,2) {};
\node (Q4) at (6,2) {};
\node (Q5) at (8,2) {};

\node (R1) at (-1,1) {};
\node (R2) at (1,1) {};
\node (R3) at (3,1) {${\bsm3\\1\esm}$};
\node (R4) at (5,1) {};
\node (R5) at (7,1) {};
\node (R6) at (9,1) {};

\draw [->] (R1)--(Q1);
\draw [->] (R1)--(P1);
\draw [->] (R2)--(Q2);
\draw [->] (R2)--(P2);
\draw [->] (R3)--(Q3);
\draw [->] (R3)--(P3);
\draw [->] (R4)--(Q4);
\draw [->] (R4)--(P4);
\draw [->] (R5)--(Q5);
\draw [->] (R5)--(P5);

\draw [->] (P1)--(R2);
\draw [->] (Q1)--(R2);
\draw [->] (P2)--(R3);
\draw [->] (Q2)--(R3);
\draw [->] (Q3)--(R4);
\draw [->] (P3)--(R4);
\draw [->] (P4)--(R5);
\draw [->] (P5)--(R6);
\draw [->] (Q4)--(R5);
\draw [->] (Q5)--(R6);

\draw [loosely dotted, thick] (-3,1)--(R1);
\draw [loosely dotted, thick] (0,4)--(Q1);
\draw [loosely dotted, thick] (12,1)--(R6);
\draw [loosely dotted, thick] (8,4)--(Q5);
\end{tikzpicture}}
\]

The third one is of type $\ZZ A^{\infty}_{\infty}$ and contains the other indecomposable $\Lambda$-modules:
\[
\scalebox{.8}{
\begin{tikzpicture}[scale=.7,>=stealth]
\node (A0) at (0,0){.};
\node (A1) at (4,0){.};
\node (A2) at (8,0){.};
\node (A3) at (12,0){.};
\node (a1) at (2,1){.};
\node (a2) at (6,1){.};
\node (a3) at (10,1){\frame{$\ \SSS_2^{-2}({\bsm 2\esm})\ $}};
\draw[loosely dotted, thick] (A0)--(-2,0);
\draw [->](A0)--(a1);\draw [->](A1)--(a2);\draw [->](A2)--(a3);
\draw [->](a1)--(A1);\draw [->](a2)--(A2);\draw [->](a3)--(A3);
\draw[loosely dotted, thick](A0)--(0,-1);\draw[loosely dotted, thick](A1)--(4,-1);\draw[loosely dotted, thick](A2)--(8,-1);\draw[loosely dotted, thick](A3)--(12,-1);\draw[loosely dotted, thick](A3)--(14,0);

\node (B0) at (0,2){.};
\node (B1) at (4,2){.};
\node (B2) at (8,2){.};
\node (B3) at (12,2){.};
\node (b1) at (2,3){${\bsm 2\\3\esm}[1]$};
\node (b2) at (6,3){.};
\node (b3) at (10,3){.};
\draw[loosely dotted, thick](B0)--(-2,2);
\draw[loosely dotted, thick](B3)--(14,2);
\draw [->](B0)--(a1);\draw [->](B1)--(a2);\draw [->](B2)--(a3);
\draw [->](a1)--(B1);\draw [->](a2)--(B2);\draw [->](a3)--(B3);

\draw [->](B0)--(b1);\draw [->](B1)--(b2);\draw [->](B2)--(b3);
\draw [->](b1)--(B1);\draw [->](b2)--(B2);\draw [->](b3)--(B3);
\node (C0) at (0,4){.};
\node (C1) at (4,4){${\bsm 2\esm}[1]$};
\node (C2) at (8,4){\frame{$\ \SSS_2^{-1}({\bsm2\esm})\ $}};
\node (C3) at (12,4){.};
\node (c1) at (2,5){${\bsm 2\\3\\1\esm}$};
\node (c2) at (6,5){${\bsm1\\2\esm}[1]$};
\node (c3) at (10,5){.};
\draw[loosely dotted, thick](C0)--(-2,4);
\draw[loosely dotted, thick](C3)--(14,4);
\draw [->](C0)--(b1);\draw [->](C1)--(b2);\draw [->](C2)--(b3);
\draw [->](b1)--(C1);\draw [->](b2)--(C2);\draw [->](b3)--(C3);
\draw [->](C0)--(c1);\draw [->](C1)--(c2);\draw [->](C2)--(c3);
\draw [->](c1)--(C1);\draw [->](c2)--(C2);\draw [->](c3)--(C3);
\node (D0) at (0,6){.};
\node (D1) at (4,6){${\bsm 2\\3\esm}$};
\node (D2) at (8,6){${\bsm3\\1\\2\esm}[1]$};
\node (D3) at (12,6){.};
\node (d1) at (2,7){.};
\node (d2) at (6,7){\frame{${\bsm &&\\& &\\&2&\\&&\esm}$}};
\node (d3) at (10,7){.};
\draw[loosely dotted, thick](D0)--(-2,6);
\draw[loosely dotted, thick](D3)--(14,6);
\draw [->](D0)--(c1);\draw [->](D1)--(c2);\draw [->](D2)--(c3);
\draw [->](c1)--(D1);\draw [->](c2)--(D2);\draw [->](c3)--(D3);
\draw [->](D0)--(d1);\draw [->](D1)--(d2);\draw [->](D2)--(d3);
\draw [->](d1)--(D1);\draw [->](d2)--(D2);\draw [->](d3)--(D3);
\node (E0) at (0,8){.};
\node (E1) at (4,8){.};
\node (E2) at (8,8){${\bsm 1\\2\esm}$};
\node (E3) at (12,8){.};
\node (e1) at (2,9){.};
\node (e2) at (6,9){.};
\node (e3) at (10,9){${\bsm 3\\1\\2\esm}$};
\draw[loosely dotted, thick](E0)--(-2,8);
\draw[loosely dotted, thick](E3)--(14,8);
\draw [->](E0)--(d1);\draw [->](E1)--(d2);\draw [->](E2)--(d3);
\draw [->](d1)--(E1);\draw [->](d2)--(E2);\draw [->](d3)--(E3);

\draw [->](E0)--(e1);\draw [->](E1)--(e2);\draw [->](E2)--(e3);
\draw [->](e1)--(E1);\draw [->](e2)--(E2);\draw [->](e3)--(E3);
\node (F0) at (0,10){.};
\node (F1) at (4,10){\frame{$\ \SSS_2({\bsm2\esm})\ $}};
\node (F2) at (8,10){.};
\node (F3) at (12,10){.};
\node (f1) at (2,11){.};
\node (f2) at (6,11){.};
\node (f3) at (10,11){.};
\draw[loosely dotted, thick](F0)--(-2,10);
\draw[loosely dotted, thick](F3)--(14,10);
\draw [->](F0)--(e1);\draw [->](F1)--(e2);\draw [->](F2)--(e3);
\draw [->](e1)--(F1);\draw [->](e2)--(F2);\draw [->](e3)--(F3);
\draw [->](F0)--(f1);\draw [->](F1)--(f2);\draw [->](F2)--(f3);
\draw [->](f1)--(F1);\draw [->](f2)--(F2);\draw [->](f3)--(F3);
\node (G0) at (0,12){.};
\node (G1) at (4,12){.};
\node (G2) at (8,12){.};
\node (G3) at (12,12){.};
\draw[loosely dotted, thick](G0)--(-2,12);
\draw[loosely dotted, thick](G3)--(14,12);
\draw[loosely dotted, thick](G0)--(0,13);
\draw[loosely dotted, thick](G1)--(4,13);\draw[loosely dotted, thick](G2)--(8,13);
\draw[loosely dotted, thick](G3)--(12,13);
\draw [->](G0)--(f1);\draw [->](G1)--(f2);\draw [->](G2)--(f3);
\draw [->](f1)--(G1);\draw [->](f2)--(G2);\draw [->](f3)--(G3);

\end{tikzpicture}}\]

The framed objects are part of the same $\SSS_2$-orbit. Hence we see that the AR-quiver of the orbit category $\Db(\Lambda)/\SSS_2$ has three connected components, one is tubular $\ZZ A_{\infty}/\tau^2$ of rank 2, one is tubular $\ZZ A_{\infty}/\tau$ of rank 1, and the third one is of type $\ZZ Q$. Hence the image of the functor $\pi:\Db(\Lambda)\to \Cc_2(\Lambda)$ is the union of the transjective component with the exceptional tubes. 

\end{example}

We end this section by stating a conjecture of a criterion on the density of the functor $\pi$.

\begin{conjecture}\label{AOconjecture}
Let $\Lambda$ be a $\tau_2$-finite algebra of global dimension at most $2$. The algebra $\Lambda$ is piecewise hereditary if and only if the functor $\Db(\Lambda)\rightarrow \Cc_2(\Lambda)$ is dense. 

Or equivalently, the algebra $\Lambda$ is piecewise hereditary if and only if any $\Pi_3(\Lambda)$-module is gradable.
\end{conjecture}
The ``only if'' part of this conjecture holds by Keller's result on triangulated orbit categories \cite{Kel05}.

This last result summarizes cases where Conjecture \ref{AOconjecture} holds.
\begin{theorem}[\cite{AO10b}]
Let $\Lambda$ be a $\tau_2$-finite algebra of global dimension at most 2. Then Conjecture~\ref{AOconjecture} holds if we are in one of the following cases:
\begin{itemize}
\item the quiver of $\Lambda$ contains an oriented cycle;
\item there exists an object $X$ in $\Db(\Lambda)$ such that $\mathbb{S}^{a}X=X[b]$ for some $a\neq b\in \ZZ$.
\end{itemize}
\end{theorem}
In the example treated above, the quiver $Q$ contains an oriented cycle, and also we have ${\bsm 1\esm}[3]\simeq \SSS^2({\bsm 1\esm})$ and $\SSS({\bsm 2\\3\\1\\2\esm})={\bsm 2\\3\\1\\2\esm}[0]$.


\def\cprime{$'$}
\frenchspacing

\end{document}